\documentclass[a4paper,article,10pt]{article}
\usepackage[english]{babel}
\usepackage[reqno]{amsmath}% AMS-LaTeX-first 
\usepackage{amssymb,amsthm,stmaryrd}
\theoremstyle{plain}
\newtheorem{theorem}{Theorem}
\newtheorem{proposition}[theorem]{Proposition}

\theoremstyle{definition}

\newtheorem{remark}[theorem]{Remark}

\usepackage[sort,compress]{cite} % bibliography
\usepackage[]{graphicx}               % graphics
\usepackage{yhmath} % for wider hats

\allowdisplaybreaks%[1]

\usepackage{xcolor}
\usepackage{ulem}

\newcommand{\bm}[1]{\text{\boldmath $#1$\unboldmath}}
\newcommand{\abs}[1]{\lvert#1\rvert}
\newcommand{\norm}[1]{\lVert#1\rVert}
\newcommand{\bddot}{\operatorname{\bm{:}}}
\newcommand{\vect}[1]{\mathbf{#1}}
\newcommand{\mat}[1]{\mathbf{#1}}
\newcommand{\node}[1]{\mathrm{#1}}
\newcommand{\Div}{{\bm{\nabla}{\cdot}}}
\newcommand{\grad}{\bm{\nabla}}
\newcommand{\curl}{\operatorname{curl}}
\newcommand{\defo}{\bm{\nabla}^{\texttt{{s}}}}
\newcommand{\RR}{\mathbb{R}}

\renewcommand{\SS}{\mathbb{S}}
\newcommand{\Ga}[1]{\Gamma_{\!\!#1}}
\newcommand{\sobo}[1][1]{\ensuremath{\mathcal{H}^{#1\!}}}
\newcommand{\Htrace}{\ensuremath{\mathcal{H}^{\frac{1}{2}}}}
\newcommand{\eltwo}{\ensuremath{\mathcal{L}_{2_{\!}}}}
\newcommand{\hDiv}[1]{\ensuremath{\mathcal{H}(\operatorname{div};{{#1}})}}
\newcommand{\Hdiv}{H(\operatorname{div})}
\newcommand{\Assem}{\mbox{\textsf{\textbf{\Large A}}}}

\newcommand{\nen}  {\ensuremath{\texttt{n}_{\texttt{en}}}}
\newcommand{\nfn}     {\texttt{n}_{\texttt{fn}}}

\newcommand{\nsd}  {\ensuremath{\texttt{n}_{\texttt{sd}}}}
\newcommand{\msd}  {\ensuremath{\texttt{m}_{\texttt{sd}}}}
\newcommand{\numel}{\ensuremath{\texttt{n}_{\texttt{el}}}}
\newcommand{\numfa}{\ensuremath{\texttt{n}_{\texttt{fa}}}}

\newcommand{\Insd}{\mat{I}_{\nsd\!}}
\newcommand{\Imsd}{\mat{I}_{\msd}}

\DeclareMathOperator{\tr}{tr}

\newcommand{\bu}{\bm{u}}
\newcommand{\bhu}{\bm{\hat{u}}}
\newcommand{\bv}{\bm{v}}

\newcommand{\bsigma}{\bm{\sigma}}
\newcommand{\bomega}{\bm{\omega}}

\newcommand{\bn}{\bm{n}}
\newcommand{\btau}{\bm{\tau}}

\newcommand{\bt}{\bm{t}}
\newcommand{\bx}{\bm{x}}
\newcommand{\bL}  {\bm{L}}
\newcommand{\bG}  {\bm{G}}

\newcommand{\bw}{\bm{w}}
\newcommand{\bhw}{\bm{\widehat{w}}}
\newcommand{\ba}{\bm{a}}
\newcommand{\bha}{\bm{\hat{a}}}

\newcommand{\VhHat}{\ensuremath{\mathcal{\hat{V}}^h}}
\newcommand{\Vh}{\ensuremath{\mathcal{V}^h}}

\newcommand{\nDeg}{\ensuremath{k}}
\newcommand{\Pk}{\ensuremath{\mathcal{P}^{\nDeg}}}

\newcommand{\hv}{\hat{v}}
\newcommand{\hu}{\hat{u}}

\newcommand{\gradS}{\bm{\nabla}_{\!\!\texttt{\footnotesize s}}}
\newcommand{\gradW}{\bm{\nabla}_{\!\texttt{\footnotesize w}}}
\newcommand{\stressV}{\bm{\sigma}_{\!\texttt{\footnotesize  v}}}
\newcommand{\strainV}{\bm{e}_{\texttt{\footnotesize  v}}}
\newcommand{\bD}{\mat{D}}
\newcommand{\bDHalf}{\bD^{1/2}}

\newcommand{\bE}{\mat{E}}
\newcommand{\bN}{\mat{N}}
\newcommand{\bT}{\mat{T}}

\newcommand{\nrr}  {\ensuremath{\texttt{n}_{\texttt{rr}}}}

\newcommand{\Ke}{\mat{K}_{\! e}}
\newcommand{\Be}{\mat{B}_{\! e}}
\newcommand{\Ce}{\mat{C}_{\! e}}
\newcommand{\De}{\mat{D}_{\! e}}
\newcommand{\Se}{\mat{S}_{\! e}}
\newcommand{\inv}[1]{#1^{-1}}
\newcommand{\pinv}[1]{#1^{\dagger}}

\newcommand{\nodaluV}{\vect{u}}
\newcommand{\nodalpV}{{\node{p}}}
\newcommand{\nodaluhV}{\vect{\hat{u}}}
\newcommand{\nodalLV}{\vect{L}}
\newcommand{\bxi}{\bm{\xi}}
\newcommand{\bet}{\bm{\eta}}

\newcommand{\nipe}{\texttt{n}_{\texttt{ip}}^{\texttt{e}}}
\newcommand{\sumge}{\sum_{\texttt{g}=1}^{\nipe}}

\newcommand{\bxige}{\bxi^\texttt{e}_\texttt{g}}
\newcommand{\wge}{w^\texttt{e}_\texttt{g}}
\newcommand{\bJ}{\mat{J}}

\newcommand{\nipf}{\texttt{n}_{\texttt{ip}}^{\texttt{f}}}
\newcommand{\sumgf}{\sum_{\texttt{g}=1}^{\nipf}}

\newcommand{\bxigf}{\bxi^\texttt{f}_\texttt{g}}
\newcommand{\wgf}{w^\texttt{f}_\texttt{g}}

\newcommand{\Nmat}{\bm{\mathcal{N}}}
\newcommand{\Mmat}{\bm{\mathcal{M}}}
\newcommand{\Qmat}{\bm{\mathcal{Q}}}
\newcommand{\NmatHat}{\bm{\widehat{\mathcal{N}}}}

\newcommand{\jump}[1]{\llbracket #1\rrbracket}
\newcommand{\bigjump}[1]{\bigl\llbracket #1\bigr\rrbracket}

\usepackage{scalerel,stackengine}
\stackMath
\newcommand\reallywidehat[1]{%
\savestack{\tmpbox}{\stretchto{%
  \scaleto{%
    \scalerel*[\widthof{\ensuremath{#1}}]{\kern-.6pt\bigwedge\kern-.6pt}%
    {\rule[-\textheight/2]{1ex}{\textheight}}%WIDTH-LIMITED BIG WEDGE
  }{\textheight}% 
}{0.5ex}}%
\stackon[1pt]{#1}{\tmpbox}%
}

\usepackage{fancyhdr}
\usepackage{lastpage} %FOR HEADER IN FIRST PAGE

\begin{document}
%________________________________________________________________________
\title{Tutorial on Hybridizable Discontinuous Galerkin (HDG) Formulation for Incompressible Flow Problems}

\author{
\renewcommand{\thefootnote}{\arabic{footnote}}
			  M. Giacomini\footnotemark[1], \
			  R. Sevilla\footnotemark[2], \ and
             A. Huerta\footnotemark[1]
}

\date{}
%________________________________________________________________________
\maketitle

\renewcommand{\thefootnote}{\arabic{footnote}}

\footnotetext[1]{Laboratori de C\`alcul Num\`eric (LaC\`aN), ETS de Ingenieros de Caminos, Canales y Puertos, Universitat Polit\`ecnica de Catalunya, Barcelona, Spain}
\footnotetext[2]{Zienkiewicz Centre for Computational Engineering, College of Engineering, Swansea University, Wales, UK
}

%________________________________________________________________________
\begin{abstract}
A hybridizable discontinuous Galerkin (HDG) formulation of the linearized incompressible Navier-Stokes equations, known as Oseen equations, is presented. 
The Cauchy stress formulation is considered and the symmetry of the stress tensor and the mixed variable, namely the scaled strain-rate tensor, is enforced pointwise via Voigt notation.
Using equal-order polynomial approximations of degree $k$ for all variables, HDG provides a stable discretization. Moreover, owing to Voigt notation, optimal convergence of order $k+1$ is obtained for velocity, pressure and strain-rate tensor and a local postprocessing strategy is devised to construct an approximation of the velocity superconverging with order $k+2$, even for low-order polynomial approximations.
A tutorial for the numerical solution of incompressible flow problems using HDG is presented, with special emphasis on the technical details required for its implementation.
\end{abstract}

%________________________________________________________________________
%\begin{keywords}
%hybridizable discontinuous Galerkin, incompressible flows, Stokes flow, incompressible Navier-Stokes, Cauchy stress formulation, Voigt notation
%\end{keywords}

%________________________________________________________________________
\section{Introduction}

Computational engineering has always been concerned with the solution of equations of mathematical physics. Development of robust, accurate and efficient techniques to approximate solutions of these problems is still an important area of research. In recent years, hybrid discretization methods have gained popularity, in particular, for flow problems. This Chapter presents the extension to incompressible flows (Stokes and Navier-Stokes) of the primer\ \cite{RS-SH:16}, which was restricted to the Poisson (thermal) problem. 

Hybrid methods have been proposed for some time. In fact, already\ \cite{Ciarlet:2002} describes a hybrid method as ``any finite element method based on a formulation where one unknown is a function, or some of its derivatives, on the set $\Omega$, and the other unknown is the trace of some of the derivatives of the same function, or the trace of the function itself, along the boundaries of the set". In fact,\ \cite{Raviart-Thomas:77} propose a discontinuous Galerkin (DG) technique such that the continuity constrain is eliminated from the finite element space and imposed by means of Lagrange multipliers on the inter-element boundaries. 
Stemming from this work, several hybrid methods have been developed using both primal (see\ \cite{Egger-EW-12b,Oikawa-15,Ern-DPE-15}) and mixed formulations  (see\ \cite{Jay-CG:04,Jay-CG:05,Jay-CG:05-GAMM,Jay-CGL:09}). The latter family of techniques is nowadays known as \emph{hybridizable discontinuous Galerkin} (HDG) method. For a literature review on hybrid discretization methods and their recent developments, the interested reader is referred to \cite{Cockburn-Encyclopedia} or \cite{MG-GS:19}.

Since its introduction, HDG has been objective of intensive research and has been applied to a large number of problems in different areas, including fluid mechanics  (see\ \cite{Jay-CG:09,Nguyen-CNP:10,JP-PNC:10-compressible,JP-NPC:10-incompressible,Nguyen-NPC:10,Nguyen-NPC:11-NS,MG-GKSH:18}), wave propagation (see\ \cite{Nguyen-NPC:11-Maxwell,Nguyen-NPC:11-acoustics,GG-GFH:13}) and solid mechanics (see\ \cite{Cockburn-SCS:09-elasticity,Cockburn-KLC:15-NonLinearElasticity,RS-SGKH:18,HDG-NEFEM-Elas}), to name but a few.

This Chapter starts from the HDG method originally proposed by \cite{Jay-CGL:09} and, following \cite{RS-SH:16}, provides a tutorial for the implementation of an HDG formulation for incompressible flow problems. Note that HDG features a mixed formulation and a hybrid variable, which is the trace of the primal one. The  \emph{hybridization}\ \cite{Fraeijs-65} (also known as \emph{static condensation} for primal formulations\ \cite{Guyan-65}) allows to reduce the number of the globally-coupled degrees of freedom of the problem\ (see\ \cite{Cockburn-CDGRS:09,Cockburn-KSC:11,GG-GMFH:13,AH-HARP:13}). Moreover, the superconvergent properties of HDG in elliptic problems allow to define an efficient and inexpensive error indicator to drive degree adaptivity procedures, not feasible in a standard CG approach, see for instance\ \cite{GG-GFH:13,GG-GFH:14} or \cite{RS-SH:18}.

Moreover, the HDG formulation presented in this Chapter approximates the Cauchy stress formulation of incompressible flow equations. In particular, Voigt notation allows to easily enforce the symmetry of second-order tensors pointwise. Optimal convergence of order $k+1$ is obtained for velocity, pressure and strain-rate tensor, even for low-order polynomial approximations. Moreover,  the local postprocessing strategy is adapted to construct an approximation of the velocity superconverging with order $k+2$, even for low-order polynomial approximations.

%This work presents a tutorial on HDG for the numerical solution of incompressible flow problems. 
The remainder of this Chapter is organized as follows.
Section\ \ref{sc:ProbSta} introduces the problem statement for incompressible flows. The implementation of HDG for the linearized incompressible Navier-Stokes, known as Oseen equations, is presented in Section\ \ref{sc:Oseen}. Section\ \ref{sc:NumEx} presents some numerical results to validate the optimal convergence properties of HDG for Stokes, Oseen and incompressible Navier-Stokes equations. Eventually, two appendices,\ \ref{ap:Symmetry} and\ \ref{ap:Implementation}, present the proof of the symmetry of the matrix of the HDG global problem for Stokes equations and the technical details required for implementation, respectively.

%________________________________________________________________________
\section{Incompressible flows: problem statement}\label{sc:ProbSta}
%________________________________________________________________________
%\subsection{Strong forms}\label{sc:NS}
In compact form, the steady-state incompressible Navier-Stokes problem in the open bounded computational domain $\Omega\in\RR^{\nsd}$ with boundary $\partial\Omega$ and $\nsd$ the number of spatial dimensions reads
\begin{equation} \label{eq:NS}
 \left\{\begin{aligned}
 -\grad\cdot(2\nu\defo{\bu} - p \Insd) + \grad {\cdot} (\bu \otimes \bu) &= \bm{s}       &&\text{in $\Omega$,}\\
   \grad\cdot\bu &= 0  &&\text{in $\Omega$,}\\
   \bu &= \bu_D  &&\text{on $\Ga{D}$,}\\
    \bigl((2\nu\defo{\bu} - p \Insd ) - (\bu \otimes \bu)\bigr)\,\bn &= \bt        &&\text{on $\Ga{N}$,}\\
 \end{aligned}\right.
\end{equation}
where the pair $(\bu,p)$ represents the unknown velocity and pressure fields, $\nu>0$ is the kinematic viscosity of the fluid, $\bn$ is the outward unit normal vector to the corresponding boundary, in this case $\Ga{N}$, $\bm{s}$ is the volumetric source term, and $\defo{\bu}$ is the strain-rate tensor, that is, the symmetric part of the gradient of velocity with $\defo{} := \bigl( \grad + \grad^T \bigr)/2$. Consequently, $\bsigma := 2\nu\defo{\bu} - p \Insd$ is the stress tensor. Recall that $[\grad\bu]_{ij}=\partial u_i/\partial x_j$.

The boundary $\partial\Omega$ is composed of two disjoint parts, the Dirichlet portion $\Ga{D}$, where the value $\bu_D$ of the velocity is imposed, and the Neumann one $\Ga{N}$. Formally, $\partial\Omega = \Ga{D} \cup \Ga{N}$, $\Ga{D} \cap \Ga{N} = \emptyset$. 
On Neumann boundaries, two typical options are found. On the one hand, material surfaces (i.e.\ $\bu\cdot\bn=0$ and, thus, $(\bu \otimes \bu)\bn=\bm{0}$) where a traction $\bt$ is applied. For this reason many fluid references define Neumann boundary conditions as $(2\nu\defo{\bu} - p \Insd)\,\bn = \bt $. On the other hand, artificial boundaries where the convection term cannot be neglected. This is typical when implementing Neumann boundaries in synthetic problems. 

\begin{remark}[Cauchy stress vs.\ velocity-pressure formulation]
It is standard for incompressible flow problems to use the pointwise solenoidal property of the velocity to modify the momentum equations. Then, the Navier-Stokes problem is rewritten as
\begin{equation*}
 \left\{\begin{aligned}
 -\grad\cdot(\nu\grad{\bu} - p \Insd) + \grad {\cdot} (\bu \otimes \bu) &= \bm{s}       &&\text{in $\Omega$,}\\
   \grad\cdot\bu &= 0  &&\text{in $\Omega$,}\\
   \bu &= \bu_D  &&\text{on $\Ga{D}$,}\\
   \bigl((\nu\grad{\bu} - p \Insd) - (\bu \otimes \bu)\bigr)\,\bn &= \bt        &&\text{on $\Ga{N}$.}\\
 \end{aligned}\right.
\end{equation*}
However, as noted in\ \cite[Section 6.5]{Donea-Huerta}, from a mechanical viewpoint the two formulations are \emph{not identical} because, in general,
\begin{equation*}
  \bsigma\,\bn= (2\nu\defo{\bu} - p \Insd )\,\bn \neq (\nu\grad{\bu} - p \Insd)\,\bn .
\end{equation*}
In fact, in a velocity-pressure formulation $\bt$ is a \emph{pseudo-traction} imposed on the Neumann boundary $\Ga{N}$, unless a Robin-type boundary condition is imposed, namely
\begin{equation*}
  \bigl((\nu\grad{\bu} - p \Insd) - (\bu \otimes \bu)\bigr)\,\bn = \bt - (\nu\grad{\bu}^T) \bn ,
\end{equation*}
which, obviously, does not correspond to the natural boundary condition associated with the operator in the partial differential equation.
\end{remark}

As usual in computational mechanics, when modeling requires it, along the same portion of a boundary it is also typical to impose Dirichlet and Neumann boundary conditions at once but along orthogonal directions. For instance, this is the case of a perfect-slip boundary, which is the limiting case of those discussed in Remark\ \ref{rm:Slip}, and of a fully-developed outflow boundary in Remark\ \ref{rm:Outflow}.
\begin{remark}[Slip/friction boundary condition]\label{rm:Slip}
Another family of physical boundary conditions that can be considered for the Navier-Stokes equations includes slip/friction boundaries, see\ \cite{Volker-02}, which, in fact, correspond to Robin-type boundary conditions. Slip boundaries are denoted by $\Ga{S}$ and need to verify: $\partial\Omega=\Ga{D}\cup\Ga{N}\cup\Ga{S}$, being $\Ga{D}$, $\Ga{N}$ and $\Ga{S}$ disjoint by pairs. 
Slip boundary conditions are:
%
%\begin{equation*}
% \left\{\begin{aligned}
%          \bu\cdot\bn    + \alpha\bn\cdot \bigl((2\nu\defo{\bu} - p \Insd ) \AH{- (\bu \otimes \bu)}\bigr)\,\bn   &=0 &&\text{on $\Ga{S}$,}\\
% \beta \bu\cdot\bt_k +        \bt_k\cdot \bigl((2\nu\defo{\bu} - p \Insd ) \AH{- (\bu \otimes \bu)}\bigr)\,\bn   &=0 &&\text{for $k=1,\dotsc ,\nsd{-}1$, on $\Ga{S}$,}\\
% \end{aligned}\right.
%\end{equation*}
\begin{equation*}
 \left\{\begin{aligned}
          \bu\cdot\bn    + \alpha\bn\cdot \bsigma \bn   &=0 &&\text{on $\Ga{S}$,}\\
 \beta \bu\cdot\bt_k +        \bt_k\cdot \bsigma \bn   &=0 &&\text{for $k=1,\dotsc ,\nsd{-}1$, on $\Ga{S}$,}\\
 \end{aligned}\right.
\end{equation*}
where $\alpha$ and $\beta$ are the penetration and friction coefficients respectively, whereas $\bn$ is the outward unit normal to $\Ga{S}$ and the tangential vectors $\bt_k$, for $k = 1,\dotsc ,\nsd {-}1$, are such that $\{\bn,\bt_1,\dotsc ,\bt_{\nsd-1}\}$ form an orthonormal system of vectors.
 
Note that symmetry-type boundary conditions can be seen as a particular case of slip boundary conditions. That is, along the plane of symmetry (in 3D and axis of symmetry in 2D) a non-penetration ($\alpha{=}0$) and perfect-slip ($\beta{=}0$) boundary condition needs to be imposed.  Obviously, this plane of symmetry can be oriented arbitrarily in the domain $\Omega$.
\end{remark}
\begin{remark}[Outflow boundary condition]\label{rm:Outflow}
Outflow surfaces are of great interest in the simulation of engineering problems, especially for internal flows.
Common choices in the literature are represented by traction-free conditions
\begin{subequations}
\begin{equation}\label{eq:outTraction}
 	(2\nu\defo{\bu} - p \Insd )\,\bn =  \bm{0}
\end{equation}
or homogeneous Neumann conditions,
\begin{equation}\label{eq:outNeumann}
	\bigl((2\nu\defo{\bu} - p \Insd ) - (\bu \otimes \bu)\bigr)\,\bn =  \bm{0} .
\end{equation}
Nonetheless, as will be shown in Section\ \ref{sc:NS-NumEx}, these conditions introduce perturbations in the flow near the outlet and do not meet the goal of outflow boundaries, that is, to model a flow exiting the domain undisturbed.
Outflow conditions for fully-developed flows can be devised as Robin-type boundary conditions  (see\ \cite{vandeVosse-VHOBGSWSB-03}) similarly to the perfect-slip boundary conditions described in Remark\ \ref{rm:Slip}. 
Given an outflow boundary $\Ga{O}$ such that $\partial\Omega=\Ga{D}\cup\Ga{N}\cup\Ga{S}\cup\Ga{O}$, being $\Ga{D}$, $\Ga{N}$, $\Ga{S}$ and $\Ga{O}$ disjoint by pairs, the corresponding conditions are
\begin{equation}\label{eq:outCorrect}
 \left\{\begin{aligned}
          \bu\cdot\bt_k   &=0 &&\text{for $k=1,\dotsc ,\nsd{-}1$, on $\Ga{O}$,}\\
		  \bn\cdot \bsigma \bn   &=0 &&\text{on $\Ga{O}$,}\\
 \end{aligned}\right.
\end{equation}
\end{subequations}
where $\bn$ is the outward unit normal to $\Ga{O}$ and the tangential vectors $\bt_k$, for $k = 1,\dotsc ,\nsd {-}1$, are such that $\{\bn,\bt_1,\dotsc ,\bt_{\nsd-1}\}$ form an orthonormal system of vectors.
\end{remark}

To further simplify the presentation, the linearized incompressible Navier-Stokes equations are considered. The so-called Oseen equations, with linear convection driven by the solenoidal field $\ba$, are
\begin{equation} \label{eq:Oseen}
 \left\{\begin{aligned}
 -\grad\cdot(2\nu\defo{\bu} - p \Insd) + \grad {\cdot} (\bu \otimes \ba) &= \bm{s}       &&\text{in $\Omega$,}\\
  \grad\cdot\bu &= 0  &&\text{in $\Omega$,}\\
  \bu &= \bu_D  &&\text{on $\Ga{D}$,}\\
  \bigl((2\nu\defo{\bu} - p \Insd ) - (\bu \otimes \ba)\bigr)\,\bn &= \bt        &&\text{on $\Ga{N}$,}\\
 \end{aligned}\right.
\end{equation}
where $\ba$ coincides with $\bu$ when confronted with Navier-Stokes.
In Einstein notation and for $i=1,\dotsc ,\nsd$, the strong form of the Oseen equations reads as
\begin{equation*}
 \left\{\begin{aligned}
  -\frac{\partial}{ \partial x_j}\biggl( \nu \Bigl(\frac{\partial u_i}{\partial x_j} + \frac{\partial u_j}{\partial x_i}\Bigr) - p\,\delta_{ij} \biggr) 
  +\frac{\partial}{\partial x_j}\bigl(u_i\, a_j \bigr) &= s_i       &&\text{in $\Omega$,}\\
                                          \frac{\partial u_j}{\partial x_j} &= 0  &&\text{in $\Omega$,}\\
                                          u_i &= u_{D,i}  &&\text{on $\Ga{D}$,}\\
      \nu\Bigl(\frac{\partial u_i}{\partial x_j} + \frac{\partial u_j}{\partial x_i}\Bigr) n_j - p n_i - u_i a_j n_j &= t_i       &&\text{on $\Ga{N}$.}\\
 \end{aligned}\right.
\end{equation*}

Finally, it is important to recall that, when inertial effects are negligible (viz.\ in micro-fluidics), incompressible flows are modeled by means of the Stokes problem. That is,
\begin{equation} \label{eq:Stokes}
 \left\{\begin{aligned}
 -\grad\cdot(2\nu\defo{\bu} - p \Insd) &= \bm{s}       &&\text{in $\Omega$,}\\
  \grad\cdot\bu &= 0  &&\text{in $\Omega$,}\\
  \bu &= \bu_D  &&\text{on $\Ga{D}$,}\\
  (2\nu\defo{\bu} - p \Insd )\,\bn &= \bt        &&\text{on $\Ga{N}$,}\\
 \end{aligned}\right.
\end{equation}
which is a second-order elliptic equation. Stokes problem can be viewed as the particular case of the Oseen equations\ \eqref{eq:Oseen}, for $\ba=\bm{0}$. Thus, in the remainder of this Chapter, the Oseen equations will be detailed. Navier-Stokes and Stokes problems will be recovered after replacing $\ba$ by $\bu$ or $\bm{0}$, respectively.

It is worth noticing that the divergence-free equation in Navier-Stokes, Oseen and Stokes problems induces a compatibility condition on the velocity field, namely
\begin{equation} \label{eq:compatibilityCondition}
\int_{\Ga{D}} \bu_D \cdot \bn\, d\Gamma + \int_{\Ga{N}} \bu \cdot \bn\, d\Gamma = 0.
\end{equation}
In particular, if only Dirichlet boundary conditions are considered (i.e.\ $\Ga{N} = \emptyset$), an additional constraint on the pressure needs to be imposed to avoid its indeterminacy. It is common for hybrid formulations (see for instance\ \cite{Jay-CG:09,Nguyen-CNP:10,Cockburn-CS:14}) to impose the mean pressure on the boundary, namely
\begin{equation} \label{eq:constraintDirichlet}
\int_{\partial\Omega} p\, d\Gamma = 0.
\end{equation}

%________________________________________________________________________
\section{HDG method for Oseen flows}\label{sc:Oseen}
%________________________________________________________________________
\subsection{Functional and discrete approximation setting}\label{sc:FunSet}

The HDG method relies on a mixed hybrid formulation of the problem under analysis.
Thus, assume that $\Omega$ is partitioned in $\numel$ disjoint subdomains $\Omega_e$,
\begin{equation*}
 \overline{\Omega} =  \bigcup_{e=1}^{\numel} \overline{\Omega}_e, \;\text{with }
 \Omega_i \cap \Omega_j = \emptyset \text{ for } i\neq j ,
\end{equation*}
with boundaries $\partial\Omega_e$, which define an internal interface, also known as \emph{internal skeleton}, $\Gamma$
\begin{equation*}%\label{eq:Gamma}
 \Gamma := \Big[ \bigcup_{e=1}^{\numel} \partial\Omega_e \Big]\setminus\partial\Omega .
\end{equation*}
Moreover, in what follows, the classical $\eltwo$ inner products for vector-valued functions are considered for a generic domain $D \subseteq \Omega \subset \RR^{\nsd}$ and a generic line/surface $S \subset\Gamma\cup\partial\Omega$, that is, 
\begin{equation*}%\label{eq:innerScalar}
 (\bu,\bw)_{D} := \int_{D} \bu \cdot \bw \, d\Omega  , \qquad 
 \langle \bu, \bw \rangle_{S} := \int_{S} \bu \cdot \bw \, d\Gamma .
\end{equation*}
The $\eltwo$ inner product for tensor-valued functions will also be used on $D$, namely
\begin{equation*}
 (\bL,\bG)_{D} := \int_{D} \bL \bddot \bG \, d\Omega  \text{, i.e.\ in Einstein notation: $\int_{D} [\bL]_{ij} [\bG]_{ij} \, d\Omega$.}
\end{equation*}
In the following subsections the Sobolev space $\sobo(D)$, $D\subseteq\Omega$, of $\eltwo(D)$ functions whose gradient also is an $\eltwo(D)$ function, will be used systematically. The space of $\eltwo(D)$ symmetric tensors $\SS$ of order $\nsd$ with $\eltwo(D)$ row-wise divergence is denoted by $\left[\hDiv{D};\SS\right]$. Moreover, the trace space $\Htrace(\partial D)$, being the space of the restriction of $\sobo(D)$ functions to the boundary $\partial D$, is introduced. Note also that for functions on $S\subset\Gamma\cup\partial\Omega$, the typical $\eltwo(S)$ space is employed. 

Moreover, the following discrete functional spaces are considered according to the notation introduced in\ \cite{RS-SH:16}
\begin{subequations}\label{eq:HDG-Spaces}
	\begin{align} 
	\Vh(\Omega) & := \left\{ v \in \eltwo(\Omega) : v \vert_{\Omega_e} \in \Pk(\Omega_e) \;\forall\Omega_e \, , \, e=1,\dotsc ,\numel \right\} , \label{eq:spaceScalarElem} \\
	\VhHat(S) & := \left\{ \hv \in \eltwo(S) : \hv\vert_{\Gamma_i}\in \Pk(\Gamma_i)	\;\forall\Gamma_i\subset S\subseteq\Gamma\cup\partial\Omega \right\}, \label{eq:spaceScalarFace} %\\
%	\Vh_\star(\Omega) & := \left\{ v \in \eltwo(\Omega) : v \vert_{\Omega_e}\in \mathcal{P}^{k+1}(\Omega_e) \;\forall\Omega_e \, , \, e=1,\dotsc ,\numel \right\}, \label{eq:Vstar}
	\end{align}
\end{subequations}
where $\Pk(\Omega_e)$ and $\Pk(\Gamma_i)$ are the spaces of polynomial functions of complete degree at most $k$ in $\Omega_e$ and on $\Gamma_i$, respectively. 
%

%________________________________________________________________________
\subsection{Strong forms of the local and global problems}\label{sc:O-Strong}

As noted earlier, the HDG method relies on a mixed hybrid formulation. Following the rationale described by \cite{RS-SH:16}, the Oseen problem reported in \eqref{eq:Oseen} is rewritten as a system of first-order equations as
\begin{equation} \label{eq:OseenBrokenFirstOrder}
\hspace{-0.5em}
\left\{
\hspace{-0.25em}
\begin{aligned}
  \bL + \sqrt{2\nu} \defo{\bu}                 &= \bm{0}    &&
\hspace{-0.5em}
  \text{in $\Omega_e$, and for $e=1,\dotsc ,\numel$,}\\	
  \grad{\cdot}(\sqrt{2\nu} \bL + p \Insd) + \grad {\cdot} (\bu \otimes \ba) &= \bm{s}    &&
\hspace{-0.5em}
  \text{in $\Omega_e$, and for $e=1,\dotsc ,\numel$,}\\
  \grad\cdot\bu                                      &= 0            &&
\hspace{-0.5em}
  \text{in $\Omega_e$, and for $e=1,\dotsc ,\numel$,}\\
  \bu                                                      &= \bu_D    &&
\hspace{-0.5em}
  \text{on $\Ga{D}$,}\\
  \bigl((\sqrt{2\nu} \bL + p \Insd) + (\bu \otimes \ba)\bigr)\,\bn            &= -\bt         &&
\hspace{-0.5em}
  \text{on $\Ga{N}$,}\\
  \jump{\bu \otimes \bn}                        &= \bm{0}   &&
\hspace{-0.5em}
  \text{on $\Gamma$,}\\
  \bigjump{\bigl((\sqrt{2\nu} \bL + p \Insd)+ (\bu \otimes \ba)\bigr)\,\bn} &= \bm{0}  &&
\hspace{-0.5em}
  \text{on $\Gamma$,}\\
\end{aligned} \right.
\end{equation}
where $\bL= - \sqrt{2\nu} \defo\bu$ is the mixed variable representing the \emph{scaled} strain-rate or second-order velocity deformation tensor. The last two equations are the so-called \emph{transmission conditions}. They impose, respectively, continuity of velocity and normal flux (normal component of the stress --i.e.\ diffusive flux-- and of the convective flux), across the interior faces. The \emph{jump} $\jump{\cdot}$ operator has been introduced following the definition by \cite{AdM-MFH:08}, such that, along each portion of the interface $\Gamma$ it sums the values from the element on the left and on the right, say $\Omega_l$ and $\Omega_r$, namely
\begin{equation*}
\jump{\odot} = \odot_l + \odot_r.
\end{equation*}
Note that the above definition of the jump operator always involves the outward unit normal to a surface, say $\jump{ \odot \bn}$ . Thus, at the interface between elements $\Omega_l$ and $\Omega_r$, this definition implies $\jump{ \odot \bn}= \odot_l \bn_l + \odot_r \bn_r$ where $\bn_l$ and $\bn_r$ are the outward unit normals to $\partial\Omega_l$ and $\partial\Omega_r$, respectively. Moreover, recall that $\bn_l=-\bn_r$ along their interface.

\begin{remark}
Recall that for incompressible flow problems with purely Dirichlet boundary conditions (i.e.\ $\Ga{N} = \emptyset$) an additional constraint is required to avoid indeterminacy of pressure. A common choice is to impose an arbitrary mean value of the pressure on the boundary (usually zero), that is
$ \langle p, 1\rangle_{\partial\Omega} = \texttt{Cst} $.
\end{remark}

%The HDG formulation solves the problem in two phases, see the seminal contribution in\ \cite{Jay-CGL:09} and the subsequent series of papers by Cockburn and coworkers\ \cite{Jay-CG:09,Nguyen-NPC:09-LinearConDif,Nguyen-NPC:09-NonLinearConDif,Nguyen-NPC:10,Nguyen-CNP:10,JP-PNC:10-compressible,JP-NPC:10-incompressible,Nguyen-NPC:11-NS,Jay-CGNPS:11} where the HDG formulation for flow problems has been theoretically and numerically analyzed. 
%
Starting from the mixed formulation on the broken computational domain, see \eqref{eq:OseenBrokenFirstOrder}, HDG features two stages.
First, a set of $\numel$ local problems is introduced to define element-by-element $(\bL_e,\bu_e,p_e)=(\bL,\bu,p)$ for all $\bx\in\Omega_e\subset\Omega$ in terms of a novel independent variable $\bhu$, namely
\begin{subequations}\label{eq:O-StrongLocal}
\begin{equation} \label{eq:O-StrongLocal-inde}
\left\{\begin{aligned}
  \bL_e + \sqrt{2\nu} \defo{\bu_e}                  &= \bm{0}   &&\text{in $\Omega_e$,} \\	
  \grad {\cdot}(\sqrt{2\nu} \bL_e + p_e \Insd) + \grad {\cdot} (\bu_e \otimes \ba) &= \bm{s}   &&\text{in $\Omega_e$,} \\
  \grad\cdot \bu_e                                          &= 0           &&\text{in $\Omega_e$,} \\
  \bu_e                                                           &= \bu_D   &&\text{on $\partial\Omega_e \cap \Ga{D}$,}\\
  \bu_e                                                           &= \bhu     &&\text{on $\partial\Omega_e \setminus \Ga{D}$,}\\
\end{aligned} \right.
\end{equation}
where $\bhu$ represents the trace of the velocity on the mesh skeleton $\Gamma \cup \Ga{N}$.
Note that \eqref{eq:O-StrongLocal-inde} constitutes a purely Dirichlet boundary value problem. As previously observed, an additional constraint needs to be added to remove the indeterminacy of the pressure, namely
\begin{equation}\label{eq:constraintLoc}
\langle p_e, 1 \rangle_{\partial\Omega_e} = \rho_e ,
\end{equation}
\end{subequations}
where $\rho_e$ denotes the scaled mean pressure on the boundary of the element $\Omega_e$.
Hence, the local problem defined by \eqref{eq:O-StrongLocal} provides $(\bL_e,\bu_e,p_e)$, for $e=1,\dotsc ,\numel$, in terms of the global unknowns $\bhu$ and $\bm{\rho}=(\rho_1\dotsc,\rho_{\numel})^T$.

The second stage computes the trace of the velocity $\bhu$ and the scaled mean pressure $\bm{\rho}$ on the element boundaries by solving the global problem accounting for the following transmission conditions and the Neumann boundary conditions
\begin{subequations}\label{eq:O-StrongGlobal}
\begin{equation} 
\left\{\begin{aligned}
%  \jump{\bu \otimes \bn}                        &= \bm{0}  &&\text{on $\Gamma$,}\\
  \bigjump{\bigl((\sqrt{2\nu} \bL + p \Insd)+ (\bu \otimes \ba)\bigr)\,\bn} &= \bm{0}  &&\text{on $\Gamma$,}\\
            \bigl((\sqrt{2\nu} \bL + p \Insd)+ (\bu \otimes \ba)\bigr)\,\bn  &= -\bt       &&\text{on $\Ga{N}$.}\\
\end{aligned} \right.
\end{equation}
Note that the continuity of velocity on $\Gamma$, $\jump{\bu \otimes \bn}=\bm{0}$, is not explicitly written because it is automatically satisfied. This is due to the Dirichlet boundary condition $\bu_e = \bhu$ imposed in every local problem and to the unique definition of the hybrid variable $\bhu$ on each edge/face of the mesh skeleton.
Moreover, the divergence-free condition in the local problem induces the following compatibility condition for each element $\Omega_e, \ e=1,\dotsc ,\numel$
\begin{equation}\label{eq:divergenceFreeConstraint}
     \langle \bhu \cdot \bn_e , 1 \rangle_{\partial \Omega_e \setminus \Ga{D}} 
  + \langle \bu_D \cdot \bn_e, 1 \rangle_{\partial \Omega_e\cap \Ga{D}} = 0 ,
\end{equation}
\end{subequations}
which is utilized to close the global problem.

\begin{remark}[Neumann local problems] 
As detailed by \cite{RS-SH:16}, an alternative HDG formulation is obtained if the Neumann boundary condition is imposed in the local problem. In this case, the trace of the velocity, $\bhu$, is only defined on $\Gamma$ rather than on $\Gamma \cup \Ga{N}$. This leads to a marginally smaller discrete global problem. For the sake of clarity, the standard HDG formulation with Neumann boundary conditions imposed in the global problem is only considered here.
\end{remark}

For a complete introduction to HDG, the interested reader is referred the seminal contribution by \cite{Jay-CG:09} and to the subsequent series of papers by Cockburn and coworkers (\cite{Jay-CGL:09,Nguyen-CNP:10,Jay-CGNPS:11,Nguyen-NPC:09-LinearConDif,Nguyen-NPC:09-NonLinearConDif,JP-NPC:10-incompressible,Nguyen-NPC:10,Nguyen-NPC:11-NS,JP-PNC:10-compressible}) where the HDG formulation for flow problems has been theoretically and numerically analyzed.

%________________________________________________________________________
\subsection{Weak forms of the local and global problems}\label{sc:O-Weak}

For each element $\Omega_e, \ e=1,\ldots,\numel$, the weak formulation of \eqref{eq:O-StrongLocal} is as follows: given $\bu_D$ on $\Ga{D}$ and $\bhu$ on $\Gamma \cup \Ga{N}$, find $(\bL_e, \bu_e, p_e) \in \left[\hDiv{\Omega_e};\SS\right] \times \left[\sobo(\Omega_e)\right]^{\nsd} \times \sobo(\Omega_e)$ that satisfies
\begin{equation*}\left\{
\begin{aligned} 
  {-}&\bigl( \bG, \bL_e \bigr)_{\Omega_e}  
  + \bigl( \Div(\sqrt{2\nu}\bG), \bu_e \bigr)_{\Omega_e}
  \\[-0.5ex] &\hspace{65pt}
  = \langle \bG\,\bn_e, \sqrt{2\nu}\,\bu_D \rangle_{\partial\Omega_e\cap\Ga{D}} 
  + \langle \bG\,\bn_e, \sqrt{2\nu}\,\bhu \rangle_{\partial \Omega_e \setminus \Ga{D} } ,
\\[1ex]
  &\bigl( \bw, \grad{\cdot}( \sqrt{2\nu} \bL_e) \bigr)_{\Omega_e}  
  {+} \bigl( \bw, \grad p_e \bigr)_{\Omega_e}  
  \\[-0.5ex] &\hspace{65pt}
  {+} \bigl\langle \bw, (\reallywidehat{\sqrt{2\nu} \bL_e {+} p_e \Insd})\,\bn_e{-}(\sqrt{2\nu} \bL_e {+} p_e \Insd)\,\bn_e \bigr\rangle_{\partial\Omega_e}
  \\[-0.5ex] &\hspace{75pt}
  {-} \bigl( \grad\bw , \bu_e \otimes \ba \bigr)_{\Omega_e} 
  {+} \langle \bw,(\widehat{\bu_e \otimes \ba})  \bn_e \rangle_{\partial\Omega_e} 
  = \bigl( \bw, \bm{s} \bigr)_{\Omega_e} ,
\\[1ex]
   &\bigl( \grad q, \bu_e \bigr)_{\Omega_e}  
 = \langle q, \bu_D \cdot \bn_e \rangle_{\partial \Omega_e \cap \Ga{D}}
 + \langle q, \bhu \cdot \bn_e \rangle_{\partial \Omega_e \setminus \Ga{D} } ,
\\
 &\langle p_e, 1 \rangle_{\partial\Omega_e}   = \rho_e ,
\end{aligned}\right.
\end{equation*}
for all $(\bG, \bw, q) \in \left[\hDiv{\Omega_e};\SS\right] \times \left[\sobo(\Omega_e)\right]^{\nsd} \times \sobo(\Omega_e)$, where, as defined in Section \ref{sc:FunSet}, $\left[\hDiv{\Omega_e};\SS\right]$ is the space of square-integrable symmetric tensors $\SS$ of order $\nsd$ on $\Omega_e$ with square-integrable row-wise divergence. 
Note that the variational form of the momentum equation above is obtained after integrating by parts the diffusive part of the flux twice, whereas the convective term is integrated by parts only once. This, as noted in\ \cite{RS-SH:16}, preserves the symmetry of the local problem for the diffusive operator and, consequently, for the Stokes problem ($\ba=\bm{0}$).

The numerical trace of the diffusive flux is defined in \cite{Jay-CGL:09} and \cite{RS-SH:16},
\begin{subequations}\label{eq:NumFlux}
\begin{equation} \label{eq:traceDiffusion}
  \hspace{-1em}
  (\reallywidehat{\sqrt{2\nu} \bL_e {+} p_e \Insd})\,\bn_e {:=} 
  \hspace{-0.25em}
    \begin{cases}
    (\sqrt{2\nu} \bL_e {+} p_e \Insd)\,\bn_e {+} \tau^d (\bu_e {-} \bu_D) &\hspace{-0.75em} \text{on $\partial\Omega_e\cap\Ga{D}$,} \\
    (\sqrt{2\nu} \bL_e {+} p_e \Insd)\,\bn_e {+} \tau^d (\bu_e {-} \bhu) &\hspace{-0.75em} \text{elsewhere.}  
  \end{cases}
\end{equation}
For the convection flux, several alternatives are possible, and they can be written in general form as 
\begin{equation} \label{eq:traceAdvection}
 (\widehat{\bu_e \otimes \ba})\bn_e := \begin{cases}
 (\bm{\tilde{u}} \otimes \bha)\bn_e + \tau^a(\bu_e-\bu_D) & \text{on $\partial\Omega_e\cap\Ga{D}$,} \\
 (\bm{\tilde{u}} \otimes \bha)\bn_e + \tau^a(\bu_e-\bhu) & \text{elsewhere,}  
 \end{cases}
\end{equation}
\end{subequations}
where $\bha$ is the trace of the convective field evaluated on the mesh edges/faces (and it is unique, as $\bhu$, on the interior faces) whereas $\bm{\tilde{u}}$ needs to be appropriately defined. The option $\bm{\tilde{u}}=\bu_e$ is rarely used because, in general, $\bm{\tilde{u}}$ can be seen as an intermediate state typically used in exact Riemann solvers (see\ \cite[Section 6.1]{Hesthaven}). 
Here, as usually done in HDG (see\ \cite{Nguyen-NPC:09-LinearConDif,Nguyen-NPC:09-NonLinearConDif,Nguyen-NPC:11-NS}) $\bm{\tilde{u}}$ is chosen such that $\bm{\tilde{u}} {=} \bu_D$ on $\Ga{D}$ and $\bm{\tilde{u}} {=} \bhu$, elsewhere. 
The coefficients $\tau^d$ and $\tau^a$ are stabilization parameters that play a crucial role on the stability, accuracy and convergence properties of the resulting HDG method, see Remark\ \ref{rm:StabPara}.

\begin{remark}[Stabilization paramaters]\label{rm:StabPara}
The influence of the stabilization parameters on the well-posedness of the HDG method has been studied extensively (see for instance\ \cite{Jay-CGL:09,Cockburn-CDG:08,Nguyen-NPC:10}).
As noted in those references, for second-order elliptic problems, $\tau^d$ is of order one \emph{for dimensionless problem}; thus, here $\tau^d$ is proportional to the viscosity, namely
\begin{equation}\label{eq:tauDif}
  \tau^d = \kappa {\nu}/{\ell}  ,
\end{equation}
where $\ell$ is the characteristic size of the problem under analysis and $\kappa >0$ is a scaling factor. 
It is worth noticing that the purely diffusive (i.e.\ Stokes) problem is not very sensitive to the stabilization parameter $\tau^d$, as will be shown later in the numerical experiments.

Stabilization of the convective term emanates from the extensive literature on DG methods for hyperbolic problems. Obviously, in presence of both diffusion and convection phenomena, an appropriate choice of the stabilization parameters $\tau^d$ and $\tau^a$ is critical to guarantee stability and optimal convergence of the HDG method, as extensively analyzed by \cite{Cockburn-CDGRS:09} and \cite{Cesmelioglu-CCNP:13,Cesmelioglu-CCQ:17}. Typically for the convection term, a characteristic velocity field of the fluid is considered to determine $\tau^a$, namely
\begin{subequations}\label{eq:tauConv}
\begin{equation}\label{eq:tauConvElement}
  \tau^a = \beta \norm{\ba}_2  \quad \text{or} \quad \tau^a = \beta \norm{\ba}_\infty ,
\end{equation}
where the above norms may be defined either locally on a single element $\Omega_e$ or globally on the domain $\Omega$ and $\beta >0$ is a positive constant independent on the Reynolds number.
Note that alternative choices for the stabilization parameters $\tau^d$ and $\tau^a$ are possible, e.g. depending on the spatial coordinate $\bx$. 
More precisely, global definitions of the parameter $\tau^a$ may be obtained by considering the maximum of the expressions reported in \eqref{eq:tauConvElement}, over all the nodes in the computational mesh, that is
\begin{equation}\label{eq:tauConvUnique}
  \tau^a = \beta \max_{\bx \in \mathcal{N}_{\Omega}} \norm{\ba(\bx)}_2  \quad \text{or} \quad \tau^a = \beta  \max_{\bx \in \mathcal{N}_{\Omega}} \norm{\ba(\bx)}_\infty ,
\end{equation}
where $\mathcal{N}_{\Omega}$ is the set of nodes of the computational mesh associated with the domain $\Omega$.
The possibility of defining different values of the stabilization on each face of $\partial\Omega_e$ is also of great interest, e.g. the expression proposed by \cite{Cesmelioglu-CCQ:17} for each face $\Gamma_e$, namely
\begin{equation}
  \tau^a \vert_{_{\Gamma_e}} = \beta \max \{ \bha {\cdot} \bn_e , 0 \} .
\end{equation}
\end{subequations}

Nonetheless, $\tau^d{+}\tau^a {-} \bha {\cdot} \bn_e$ needs to be nonnegative on all the faces of the element and positive at least on one, see\ \cite{Cockburn-CDGRS:09}.
Thus, for any pair $(\tau^d,\tau^a)$, the following admissibility condition is introduced
\begin{equation*}
  \min_{\bx \in \partial\Omega_e} \{ \tau^d{+}\tau^a {-} \bha {\cdot} \bn_e \} \geq \gamma > 0 .
\end{equation*}
\end{remark}

Introducing the definition of the numerical traces in \eqref{eq:NumFlux} into the momentum equation leads to the weak form of the local problem: for $e=1,\dotsc,\numel$, find $(\bL_e, \bu_e, p_e) \in \left[\hDiv{\Omega_e};\SS\right] \times \left[\sobo(\Omega_e)\right]^{\nsd} \times \sobo(\Omega_e)$ that satisfies
\begin{equation}\label{eq:StokesWeakLocal}
\hspace{-1em}\left\{\begin{aligned} 
  {-} &\bigl( \bG, \bL_e \bigr)_{\Omega_e}  
 {+} \bigl( \Div(\sqrt{2\nu}\bG), \bu_e \bigr)_{\Omega_e}
 \\[-0.75ex] &\hspace{55pt}
 =   \langle \bG\,\bn_e, \sqrt{2\nu}\,\bu_D \rangle_{\partial\Omega_e\cap\Ga{D}} 
 {+} \langle \bG\,\bn_e, \sqrt{2\nu}\,\bhu \rangle_{\partial \Omega_e \setminus \Ga{D} }  ,
\\
  &\bigl( \bw, \grad{\cdot}( \sqrt{2\nu} \bL_e) \bigr)_{\Omega_e}  
 {+} \bigl( \bw, \grad p_e \bigr)_{\Omega_e} 
  {-} \bigl( \grad\bw , \bu_e \otimes \ba \bigr)_{\Omega_e} 
 {+} \langle \bw , \tau \bu_e \rangle_{\partial\Omega_e} 
 \\[-0.75ex] &\hspace{55pt}
 = \bigl( \bw, \bm{s} \bigr)_{\Omega_e} 
 \\[-0.75ex] &\hspace{65pt}
  {+} \langle \bw, (\tau {-} \bha{\cdot}\bn_e) \bu_D \rangle_{\partial\Omega_e \cap \Ga{D}}
  {+} \langle \bw, (\tau {-} \bha{\cdot}\bn_e) \bhu \rangle_{\partial\Omega_e \setminus \Ga{D}} ,    
\\
   &\bigl( \grad q, \bu_e \bigr)_{\Omega_e}  
 = \langle q, \bu_D \cdot \bn_e \rangle_{\partial \Omega_e \cap \Ga{D}}
 {+} \langle q, \bhu \cdot \bn_e \rangle_{\partial \Omega_e \setminus \Ga{D} } ,
\\
   &\langle p_e, 1 \rangle_{\partial\Omega_e}   = \rho_e ,
\end{aligned}\right.
\end{equation}
for all $(\bG, \bw, q) \in \left[\hDiv{\Omega_e};\SS\right] \times \left[\sobo(\Omega_e)\right]^{\nsd} \times \sobo(\Omega_e)$ and where the stabilization parameter is defined as $\tau = \tau^d + \tau^a$. Note that, as expected, the previous problem provides $(\bL_e,\bu_e,p_e)$, for $e=1,\dotsc,\numel$, in terms of the global unknowns $\bhu$ and $\bm{\rho}=(\rho_1\dotsc,\rho_{\numel})^T$.

For the global problem, the weak formulation equivalent to \eqref{eq:O-StrongGlobal} is: find $\bhu \in \left[\Htrace(\Gamma \cup \Ga{N})\right]^{\nsd}$ and $\bm{\rho} \in \mathbb{R}^{\numel}$ that satisfies 
\begin{equation*}
   \left\{\begin{aligned}
  \sum_{e=1}^{\numel} \Big\{ 
        \bigl\langle\bhw , (\reallywidehat{\sqrt{2\nu}\bL_e{+}p_e\Insd})\,\bn_e
                   {+} (\widehat{\bu_e {\otimes} \ba})\,\bn_e \bigr\rangle_{\partial\Omega_e\setminus\partial\Omega} 
  \hspace{100pt} & \\[-1.5ex] 
  {+} \bigl\langle \bhw , (\reallywidehat{\sqrt{2\nu}\bL_e{+}p_e\Insd})\,\bn_e
                   {+} (\widehat{\bu_e {\otimes} \ba})\,\bn_e {+} \bt \bigr\rangle_{\partial \Omega_e \cap \Ga{N}}  \Big\}
  {=} 0 , &
  \\
       \langle \bhu \cdot \bn_e , 1 \rangle_{\partial \Omega_e \setminus \Ga{D} } 
  = - \langle \bu_D  \cdot \bn_e , 1 \rangle_{\partial \Omega_e \cap \Ga{D}} \qquad \text{for $e=1,\dotsc,\numel$}, &
  \end{aligned}\right.
\end{equation*}
for all $\bhw \in \left[\eltwo(\Gamma \cup \Ga{N})\right]^{\nsd}$. 

Replacing in the transmission equations the numerical traces defined in \eqref{eq:NumFlux}, the variational form of the global problem is obtained. Namely, find $\bhu \in \left[\Htrace(\Gamma \cup \Ga{N})\right]^{\nsd}$ and $\bm{\rho} \in \mathbb{R}^{\numel}$ such that, for all $\bhw \in \left[\eltwo(\Gamma \cup \Ga{N})\right]^{\nsd}$, it holds
\begin{equation} \label{eq:OseenWeakGlobal}
\left\{\begin{aligned}
  \sum_{e=1}^{\numel} \Big\{ 
     \langle\bhw , (\sqrt{2\nu}\bL_e{+}p_e\Insd)\,\bn_e \rangle_{\partial \Omega_e\setminus\Ga{D}}    
  + \langle\bhw , \tau\bu_e \rangle_{\partial \Omega_e\setminus\Ga{D}}    
  \hspace{50pt} & \\[-2ex]
  -  \langle \bhw , \tau \bhu \rangle_{\partial \Omega_e\cap\Gamma}
  -  \langle \bhw , (\tau {-} \bha {\cdot} \bn_e)\bhu \rangle_{\partial \Omega_e\cap\Ga{N}}  \Big\}
  \hspace{20pt}&  \\[-1.5ex]
  = {-}\sum_{e=1}^{\numel} \langle \bhw , \bt \rangle_{\partial \Omega_e \cap \Ga{N}} & ,
  \\
       \langle \bhu \cdot \bn_e , 1  \rangle_{\partial \Omega_e \setminus \Ga{D} }  
  =   -\langle \bu_D  \cdot \bn_e , 1 \rangle_{\partial \Omega_e \cap \Ga{D}} \qquad \text{for $e=1,\dotsc,\numel$} & .
\end{aligned}\right.
\end{equation}
Remark \ref{rm:Hint} shows how the first equation in \eqref{eq:OseenWeakGlobal} is obtained exploiting the uniqueness of $\bhu$ and $\bha$ on the internal skeleton, recall that for Stokes it is trivial ($\bha=\bm{0}$) and for Navier-Stokes it follows from $\bha=\bhu$ on $\Gamma\cup\Ga{N}$. 

\begin{remark}\label{rm:Hint}
In order to obtain the first equation in \eqref{eq:OseenWeakGlobal} the usual choice $\bm{\tilde{u}}= \bhu$ is employed in the definition of the convection numerical flux, see \eqref{eq:traceAdvection}. Thus, 
\begin{equation*}
  \begin{aligned}
    & \sum_{e=1}^{\numel}\langle \bhw,  (\widehat{\bu_e \otimes \ba })\bn_e \rangle_{\partial\Omega_e\setminus \Ga{D}}  
    \\[-1.5ex] & \hspace{15pt}
    = \sum_{e=1}^{\numel}\langle \bhw,  (\bhu \otimes \bha)\bn_e +  \tau^a(\bu_e {-} \bhu) \rangle_{\partial\Omega_e\setminus \Ga{D}}  
    \\[-0.75ex] & \hspace{15pt}
    = \sum_{e=1}^{\numel}\langle \bhw,  \tau^a\bu_e - (\tau^a{-}\bha{\cdot}\bn_e) \bhu  \rangle_{\partial\Omega_e\setminus \Ga{D}}
    \\[-0.75ex] & \hspace{15pt}
    = \sum_{e=1}^{\numel}\Bigl\{
                                            \langle \bhw,  \tau^a\bu_e \rangle_{\partial\Omega_e\setminus \Ga{D}} 
                                          {-} \langle \bhw,  (\tau^a{-}\bha{\cdot}\bn_e) \bhu  \rangle_{\partial\Omega_e\cap\Gamma} 
    \\[-1.5ex] & \hspace{190pt}
                                          {-} \langle \bhw,  (\tau^a{-}\bha{\cdot}\bn_e) \bhu  \rangle_{\partial\Omega_e\cap\Ga{N}} 
                                            \Bigr\}
    \\[-0.75ex] & \hspace{15pt}
    = \sum_{e=1}^{\numel}\Bigl\{
                                            \langle \bhw,  \tau^a\bu_e \rangle_{\partial\Omega_e\setminus \Ga{D}} 
                                          {-} \langle \bhw,  \tau^a \bhu  \rangle_{\partial\Omega_e\cap\Gamma} 
                                          {-} \langle \bhw,  (\tau^a{-}\bha{\cdot}\bn_e) \bhu  \rangle_{\partial\Omega_e\cap\Ga{N}} 
                                            \Bigr\} , 
  \end{aligned}
\end{equation*}
where the last line follows from the uniqueness of $\bhw$, $\bhu$ and $\bha$ on the internal skeleton $\Gamma$ and from the relationship $\bn_l = - \bn_r$ between the outward unit normal vectors on an internal edge/face shared by two neighboring elements $\Omega_l$ and $\Omega_r$.
\end{remark}

%________________________________________________________________________
\subsection{Discrete forms and the resulting linear system}\label{sc:O-Eqs}

The discrete functional spaces defined in \eqref{eq:HDG-Spaces} are used in this section to construct the block matrices involved in the discretization of the HDG local and global problems. Moreover, following the rationale proposed in\ \cite{MG-GKSH:18} for the Stokes equations, the HDG-Voigt formulation is introduced to discretize the Cauchy stress form of the incompressible flow equations under analysis.
The advantage of this formulation compared with classical HDG approaches lies in the easy implementation of the space $\SS$ of symmetric tensors of order $\nsd$ which thus allows to retrieve optimal convergence and superconvergence properties even for low-order approximations.

%________________________________________________________________________
\subsubsection{Voigt notation for symmetric second-order tensors}\label{sc:Voigt}

A symmetric second-order tensor is written in Voigt notation by storing its diagonal and off-diagonal components in vector form, after an appropriate rearrangement.
Exploiting the symmetry, only $\msd = \nsd(\nsd+1)/2$ components (i.e.\ three in 2D and six in 3D) are stored. For the strain-rate tensor $\defo{\bu}$, the following column vector $\strainV \in \RR^{\msd}$ is obtained 
\begin{equation} \label{eq:strainVoigt}
  \strainV := \begin{cases}
                    \bigl[e_{11} ,\; e_{22} ,\; e_{12} \bigr]^T                                               &\text{in 2D,} \\
                    \bigl[e_{11} ,\; e_{22} ,\; e_{33} ,\; e_{12} ,\; e_{13} ,\; e_{23} \bigr]^T &\text{in 3D,} 
  \end{cases}
\end{equation}
where the arrangement proposed by\ \cite{FishBelytschko2007} has been utilized and the components are defined as
\begin{equation} \label{eq:NormalShear}
 e_{ij} := \frac{\partial u_i}{\partial x_j} + (1-\delta_{ij}) \frac{\partial u_j}{\partial x_i}, \quad \text{for $i,j = 1,\dotsc,\nsd$ and $i\le j$ ,}
\end{equation}
being $\delta_{ij}$ the classical Kronecker delta.
Moreover, the strain-rate tensor can be written as $\strainV = \gradS \bu$ by introducing the $\msd \times \nsd$ matrix 
\begin{equation} \label{eq:symmGrad}
  \gradS := \begin{cases}
                  \begin{bmatrix}
                    \partial/\partial x_1 & 0 & \partial/\partial x_2 \\
                    0 & \partial/\partial x_2 & \partial/\partial x_1
                  \end{bmatrix}^T                                                      &\hspace{-1em}\text{in 2D,} \\[1.5em]
                  \begin{bmatrix}
                    \partial/\partial x_1 & 0 & 0 & \partial/\partial x_2 & \partial/\partial x_3 & 0 \\
                    0 & \partial/\partial x_2 & 0 & \partial/\partial x_1 & 0 & \partial/\partial x_3 \\
                    0 & 0 & \partial/\partial x_3 & 0 & \partial/\partial x_1 & \partial/\partial x_2
                  \end{bmatrix}^T                                                      &\hspace{-1em}\text{in 3D.} 
                  \end{cases}
\end{equation}

\begin{remark}
The components of the strain-rate tensor $\defo{\bu}$ can be retrieved from its Voigt counterpart by multiplying the off-diagonal terms $e_{ij}, \ i \neq j$ by a factor $1/2$, namely
\begin{equation} \label{eq:strainClassicVoigt}
  \defo{\bu} :=\begin{cases}
                      \begin{bmatrix}
                        e_{11} & e_{12}/2 \\
                        e_{12}/2 & e_{22}
                      \end{bmatrix}                                 &\text{in 2D,} \\[1.5em]
                      \begin{bmatrix}
                        e_{11}    & e_{12}/2 & e_{13}/2 \\
                        e_{12}/2 & e_{22}    & e_{23}/2 \\
                        e_{13}/2 & e_{23}/2 & e_{33}
                      \end{bmatrix}                                 &\text{in 3D.} 
                      \end{cases}
\end{equation}
\end{remark}

Moreover, recall the following definitions introduced in\ \cite{MG-GKSH:18} for Voigt notation.
First, the vorticity $\bomega := \grad \times \bu$ is expressed using Voigt notation as $\bomega = \gradW \bu$, where the $\nrr \times \nsd$ matrix $\gradW$, with $\nrr=\nsd(\nsd-1)/2$ number of rigid body rotations (i.e.\ one in 2D and three in 3D), has the form
\begin{equation} \label{eq:curlVoigt}
\gradW :=\begin{cases}
\bigl[-\partial/\partial x_2 ,\; \partial/\partial x_1 \bigr]
&\text{in 2D,} \\
\begin{bmatrix}
0 & -\partial/\partial x_3 & \partial/\partial x_2 \\
\partial/\partial x_3 & 0 & -\partial/\partial x_1 \\
-\partial/\partial x_2 & \partial/\partial x_1 & 0
\end{bmatrix}
&\text{in 3D.} 
\end{cases}
\end{equation}
\begin{remark}[Vorticity in 2D]
Recall that in 2D the $\curl$ of a vector $\bv = [v_1, v_2]^T$ is a scalar quantity.
By setting $v_3 = 0$, the resulting vector $\bv = [v_1, v_2, 0]^T$ is embedded in the three dimensional space $\RR^3$.
Thus, $\grad \times \bv$ may be interpreted as a vector pointing along $x_3$ and with magnitude equal to $ -({\partial v_1}/{\partial x_2}) + ({\partial v_2}/{\partial x_1})$.
\end{remark}

Stokes' law is expressed in Voigt notation as $\stressV = \bD \gradS \bu -\bE p$, where the vector $\bE \in \RR^{\msd}$ and the matrix $\bD \in \RR^{\msd \times \msd}$ are defined as
\begin{equation} \label{eq:EDVoigt}
  \bE :=\begin{cases}
              \bigl[1 ,\; 1 ,\; 0 \bigr]^T                               &\hspace{-0.5em}\text{in 2D,} \\[1ex]
              \bigl[1 ,\; 1 ,\; 1 ,\; 0 ,\; 0 ,\; 0 \bigr]^T          &\hspace{-0.5em}\text{in 3D.} 
            \end{cases}
            \;
  \bD :=\begin{cases}
              \begin{bmatrix}
                2\nu \Insd & \bm{0}_{\nsd \times 1} \\
                \bm{0}_{\nsd \times 1}^T & \nu
              \end{bmatrix}                                                        &\hspace{-0.5em}\text{in 2D,} \\[1em]
              \begin{bmatrix}
                2\nu \Insd & \bm{0}_{\nsd} \\
                \bm{0}_{\nsd} & \nu \Insd
              \end{bmatrix}                                                        &\hspace{-0.5em}\text{in 3D.} 
              \end{cases}
\end{equation}

Similarly, traction boundary conditions on $\Ga{N}$ are rewritten as $\bN^T \stressV = \bt$, where $\bN$ is the $\msd \times \nsd$ matrix
\begin{equation} \label{eq:normalVoigt}
\bN :=\begin{cases}
            \begin{bmatrix}
              n_1 & 0 & n_2 \\
              0 & n_2 & n_1
            \end{bmatrix}^T                        &\text{in 2D,} \\[1em]
            \begin{bmatrix}
              n_1 & 0 & 0 & n_2 & n_3 & 0\\
              0 & n_2 & 0 & n_1 & 0 & n_3 \\
              0 & 0 & n_3 & 0 & n_1 & n_2
            \end{bmatrix}^T                        &\text{in 3D,} 
            \end{cases}
\end{equation}
which accounts for the normal direction to the boundary.

For the sake of completeness, the matrix $\bT \in \RR^{\nrr \times \nsd}$ denoting the tangential direction $\btau$ to a line/surface is introduced
\begin{equation} \label{eq:tangentVoigt}
\bT :=\begin{cases}
\bigl[-n_2 ,\; n_1 \bigr]
&\text{in 2D,} \\
\begin{bmatrix}
0 & -n_3 & n_2 \\
n_3 & 0 & -n_1 \\
-n_2 & n_1 & 0
\end{bmatrix}
&\text{in 3D,} 
\end{cases}
\end{equation}
and the projection of a vector $\bu$ along the such direction is given by $\bu \cdot \btau = \bT \bu$.

Voigt notation is thus exploited to rewrite the Oseen equations in\ \eqref{eq:Oseen}. First, recall that the divergence operator applied to a symmetric tensor is equal to the transpose of the matrix $\gradS$ introduced in\ \eqref{eq:symmGrad}.
Moreover, the vector $\bE$ in\ \eqref{eq:EDVoigt} is utilized to compactly express the trace of a symmetric tensor required in the incompressibility equation, namely $\grad\cdot\bu = \tr(\defo{\bu}) = \bE^T \gradS \bu = 0$.
Finally, the formulation of the Oseen problem equivalent to \eqref{eq:Oseen} using Voigt notation is obtained by the combining the above matrix equations
\begin{equation} \label{eq:StokesVoigt}
  \left\{\begin{aligned}
  -\gradS^T ( \bD \gradS \bu -\bE p ) + \grad \cdot (\bu \otimes \ba) &= \bm{s} &&\text{in $\Omega$,}\\
   \bE^T \gradS \bu &= 0        &&\text{in $\Omega$,}\\
   \bu &= \bu_D  &&\text{on $\Ga{D}$,}\\
   \bN^T ( \bD \gradS \bu -\bE \, p ) - (\bu \otimes \ba)\, \bn &= \bt         &&\text{on $\Ga{N}$.}
  \end{aligned}\right.
\end{equation}
Of course, as previously observed, replacing $\ba$ by $\bu$ or $\bm{0}$ leads to the Navier-Stokes, see\ \eqref{eq:NS}, and the Stokes, see\ \eqref{eq:Stokes}, problems, respectively.

%________________________________________________________________________
\subsubsection{HDG-Voigt strong forms}\label{sc:O-Voigt-Strong}

The HDG local and global problems introduced in Section\ \ref{sc:O-Strong} are now rewritten using Voigt notation.
The local problem in\ \eqref{eq:O-StrongLocal} becomes
\begin{equation} \label{eq:O-StrongLocal-Voigt}
\left\{\begin{aligned}
\bL_e + \bDHalf \gradS \bu_e &= \bm{0}    &&\text{in $\Omega_e$,}\\	
\gradS^T \bDHalf \bL_e + \gradS^T \bE \, p_e  + \grad \cdot (\bu \otimes \ba) &= \bm{s}          &&\text{in $\Omega_e$,}\\
\bE^T \gradS \bu_e &= 0           &&\text{in $\Omega_e$,}\\
\bu_e &= \bu_D     &&\text{on $\partial\Omega_e \cap \Ga{D}$,}\\
\bu_e &= \bhu  &&\text{on $\partial\Omega_e \setminus \Ga{D}$,}\\
\langle p_e, 1 \rangle_{\partial\Omega_e} &= \rho_e ,
\end{aligned} \right.
\end{equation}
where the additional constraint described in\ \eqref{eq:constraintLoc} to remove the indeterminacy of pressure is included.
Note that $\bL_e$ is the restriction to element $\Omega_e$ of the HDG mixed variable introduced to represent the strain-rate tensor in Voigt notation scaled via the matrix $\bDHalf$, featuring the square root of the eigenvalues of the diagonal matrix $\bD$, see Equation\ \eqref{eq:EDVoigt}.

Similarly, the global problem using Voigt notation, see\ \eqref{eq:O-StrongGlobal}, is
\begin{equation} \label{eq:O-StrongGlobal-Voigt}
  \hspace{-1.5em}
  \left\{\hspace{-0.25em}\begin{aligned}
%          \jump{\bN^T (\bDHalf \bL + \bE \,p) + (\bu \otimes \ba) \bn} &= \bm{0}  &&\text{on $\Gamma$,}\\
          \bigjump{\bN^T (\bDHalf \bL + \bE \,p) + (\bu \otimes \ba) \bn} &= \bm{0}  &&\text{on $\Gamma$,}\\
          \bN^T (\bDHalf \bL + \bE \,p) + (\bu \otimes \ba) \bn &= -\bt  &&\text{on $\Ga{N}$,}\\
          \langle \bE^T \bN_e \bhu , 1 \rangle_{\partial \Omega_e \setminus \Ga{D}} 
          + \langle \bE^T \bN_e \bu_D , 1 \rangle_{\partial \Omega_e\cap \Ga{D}} &= 0 &&\text{for $e=1,\dotsc ,\numel$},
          \end{aligned} \right.
\end{equation}
where the compatibility condition to close the global problem, see\ \eqref{eq:divergenceFreeConstraint}, is also written in Voigt notation.

Moreover, by exploiting Voigt notation, the definition of the trace of the diffusive numerical flux, see\ \eqref{eq:traceDiffusion}, becomes
\begin{equation} \label{eq:traceDiffusion-Voigt}
  \hspace{-1em}
  \bN_e^T (\reallywidehat{\bDHalf \bL_e {+} \bE \, p_e}) {:=} 
    \begin{cases}
    \bN_e^T ( \bDHalf \bL_e {+} \bE \, p_e ) + \tau^d (\bu_e {-} \bu_D) &\hspace{-0.75em}\text{on $\partial\Omega_e\cap\Ga{D}$,} \\
    \bN_e^T ( \bDHalf \bL_e {+} \bE \, p_e ) + \tau^d (\bu_e {-} \bhu)    &\hspace{-0.75em}\text{elsewhere.}  
    \end{cases}
  \hspace{-0.5em}
\end{equation}
Of course, the trace of the convective numerical flux, since it follows the standard tensorial notation, is the one defined in\ \eqref{eq:traceAdvection}, with $\bm{\tilde{u}}=\bhu$.

%________________________________________________________________________
\subsubsection{HDG-Voigt discrete weak forms}\label{sc:O-Voigt-Weak}

Following the derivation in Section\ \ref{sc:O-Weak}, the HDG formulation of the Oseen equations with Voigt notation is obtained.
The discrete weak formulation of the local problems proposed in\ \eqref{eq:O-StrongLocal-Voigt} is as follows: for $e=1,\dotsc ,\numel$, given $\bu_D$ on $\Ga{D}$ and $\bhu$ on $\Gamma\cup\Ga{N}$, find $(\bL_e ,\bu_e,p_e) \in [\Vh(\Omega_e)]^{\msd} \times [\Vh(\Omega_e)]^{\nsd} \times \Vh(\Omega_e)$ such that
\begin{subequations}\label{eq:O-WeakLocal-Voigt}
\begin{gather}
  \label{eq:O-Loc-L}
  \begin{aligned}
  - (\bv,\bL_e)_{\Omega_e} + (\gradS^T \bDHalf \bv, \bu_e)_{\Omega_e}
    \hspace{150pt} & \\[-0.5ex]
  = \langle \bN_e^T \bDHalf \bv , \bu_D\rangle_{\partial\Omega_e\cap\Ga{D}} 
  + \langle \bN_e^T \bDHalf \bv , \bhu \rangle_{\partial\Omega_e\setminus\Ga{D}} , &
    \end{aligned}
  \\[1ex]
  \label{eq:O-Loc-u}
  \begin{aligned}
    \bigl(\bw, \gradS^T \bDHalf \bL_e \bigr)_{\Omega_e} {+} \bigl(\bw, \gradS^T \bE \, p_e \bigr)_{\Omega_e}
    \hspace{130pt} & \\[-0.5ex]
    {-} (\grad \bw, \bu_e \otimes\ba)_{\Omega_e} 
    {+} \langle \bw , \tau \bu_e \rangle_{\partial\Omega_e}
    \hspace{100pt} & \\[-0.5ex]
    = (\bw,\bm{s} )_{\Omega_e}
    \hspace{180pt} & \\[-0.5ex]
    {+} \langle \bw, (\tau {-} \bha{\cdot}\bn_e) \bu_D \rangle_{\partial\Omega_e \cap \Ga{D}} 
    {+} \langle \bw, (\tau {-} \bha{\cdot}\bn_e) \bhu \rangle_{\partial\Omega_e \setminus \Ga{D}} , &  
    \end{aligned}
  \\[1ex]
  \label{eq:O-Loc-p}
  (\gradS^T \bE \, q, \bu_e)_{\Omega_e}  
  = \langle \bN_e^T \bE\, q ,  \bu_D \rangle_{\partial\Omega_e\cap\Ga{D}} 
  + \langle \bN_e^T \bE\, q ,  \bhu \rangle_{\partial\Omega_e\setminus\Ga{D}} , 
  \\[1ex]
  \label{eq:O-Loc-rho}
  \langle p_e, 1 \rangle_{\partial\Omega_e} = \rho_e  , 
\end{gather}
\end{subequations}
for all $(\bv ,\bw,q) \in [\Vh(\Omega_e)]^{\msd} \times [\Vh(\Omega_e)]^{\nsd} \times \Vh(\Omega_e)$, where, as done previously, the stabilization parameter is defined as $\tau = \tau^d + \tau^a$.

The global problem computes the hybrid variable $\bhu\in[\VhHat(\Gamma\cup\Ga{N})]^{\nsd}$ and the scaled mean pressure on each element boundary $\bm{\rho}\in\RR^{\numel}$.
The global problem is obtained starting from\ \eqref{eq:O-StrongGlobal-Voigt}, recall also Remark\ \ref{rm:Hint}, and the definitions of the numerical flux given in\ \eqref{eq:traceDiffusion-Voigt} and\ \eqref{eq:traceAdvection}: find $\bhu\in[\VhHat(\Gamma\cup\Ga{N})]^{\nsd}$ and $\bm{\rho}\in\RR^{\numel}$, such that 
\begin{subequations}\label{eq:O-WeakGlobal-Voigt}
\begin{align}
  \label{eq:O-Glob-Transmission}
  \begin{aligned}
    \sum_{e=1}^{\numel}\Bigl\{
         \langle \bhw,  \bN_e^T \bigl(\bDHalf \bL_e {+} \bE \, p_e\bigr) \rangle_{\partial\Omega_e\setminus\Ga{D}} 
        {+} \langle \bhw, \tau\, \bu_e \rangle_{\partial\Omega_e\setminus\Ga{D}} 
        \hspace{60pt} &\\[-3ex]
        {-} \langle \bhw, \tau \bhu \rangle_{\partial\Omega_e\cap\Gamma} 
        {-} \langle \bhw, (\tau {-} \bha {\cdot} \bn_e) \bhu \rangle_{\partial\Omega_e\cap\Ga{N}} \Bigr\}             
        = -\sum_{e=1}^{\numel} \langle \bhw, \bt \rangle_{\partial\Omega_e\cap\Ga{N}} ,  & \\
  \end{aligned} &
  \\
  \label{eq:O-Glob-Compatibility}
  \langle \bE^T \bN_e \bhu, 1 \rangle_{\partial \Omega_e \setminus \Ga{D}} 
     = - \langle \bE^T \bN_e \bu_D, 1 \rangle_{\partial \Omega_e\cap \Ga{D}} \; \text{ for $e=1,\dotsc,\numel$,} &
\end{align}
\end{subequations}
 for all $\bhw\in[\VhHat(\Gamma\cup\Ga{N})]^{\nsd}$.

%________________________________________________________________________
\subsubsection{HDG linear system}\label{sc:O-Voigt-Eqs}

The discretization of the weak form of the local problem given by\ \eqref{eq:O-WeakLocal-Voigt} using an isoparametric formulation for the primal,  mixed and hybrid variables leads to a linear system with the following structure
\begin{equation} \label{eq:localProblemSystem}
  \begin{bmatrix}
    \mat{A}_{LL}    & \mat{A}_{Lu} & \mat{0}              & \mat{0}                   \\
    \mat{A}_{Lu}^T & \mat{A}_{uu} & \mat{A}_{pu}^T & \mat{0}                   \\
    \mat{0}             & \mat{A}_{pu} & \mat{0}              & \mat{a}_{\rho p}^T \\
    \mat{0}             & \mat{0}          & \mat{a}_{\rho p}& 0                   \\
  \end{bmatrix}_{\! e}
  \hspace{-0.5ex}
\begin{Bmatrix}
   \vect{L}_e \\
   \vect{u}_e \\
   \vect{p}_e \\
   \zeta 
  \end{Bmatrix} 
  {=} 
  \begin{Bmatrix}
    \vect{f}_L  \\
    \vect{f}_u  \\
    \vect{f}_p \\
    0 
  \end{Bmatrix}_{\! e} 
  \hspace{-0.5ex}
  {+}
  \begin{bmatrix}
    \mat{A}_{L\hu} \\
    \mat{A}_{u\hu} \\
    \mat{A}_{p\hu} \\
    \mat{0}
  \end{bmatrix}_{\! e}
  \hspace{-0.5ex}
  \vect{\hu}_e
  {+}
  \begin{Bmatrix}
    \vect{0} \\
    \vect{0} \\
    \vect{0} \\
    1
  \end{Bmatrix}_{\! e}
  \hspace{-0.5ex}
  \rho_e,
\end{equation}
for $e=1,\dotsc ,\numel$. It is worth noticing that the last equation in\ \eqref{eq:O-WeakLocal-Voigt} is the restriction \eqref{eq:constraintLoc} imposed via the Lagrange multiplier $\zeta$ in the system of equations above. Recalling the dimensions of the unknown variables, the solution of this system implies inverting a matrix of dimension $\bigl((\msd + \nsd + 1)\nen+1\bigr)$ for each element in the mesh, being $\nen$ the number of element nodes of $\Omega_e$. Note that $\nen$ is directly related to the degree $k$ of polynomial approximations. Table \ref{tb:neqLocal} shows the dimension of the local problem stated in \eqref{eq:localProblemSystem} for Lagrange elements on simplexes and parallelepipeds for 2D and 3D and different orders of approximation.  
\begin{table} 
\caption{Dimension of the local problem}
\centering
\begin{tabular}{lrrrrrrrr}\hline
 Order of interpolation &   $1$ &     $2$ &     $3$ &          $4$ &          $5$ &           $7$ &            $9$ &          $11$ \\ \hline
 Simplexes \\ \hline
 2D                              & $19$ &   $37$ &   $61$ &        $91$ &      $127$ &      $217$  &        $331$ &        $469$  \\
 3D                              & $41$ & $101$ & $201$ &      $351$ &      $561$ & $1\, 201$  &   $2\, 201$ &   $3\, 641$  \\ \hline\hline
 Parallelepipeds \\ \hline
 2D                              & $25$ &   $55$ &   $97$ &      $151$ &      $217$ &      $385$  &        $601$ &        $865$  \\
 3D                              & $81$ & $271$ & $641$ & $1\, 251$ & $2\, 161$ & $5\, 121$  & $10\, 001$ & $17\, 281$  \\ \hline
\end{tabular}	
\label{tb:neqLocal}
\end{table}

Similarly, the following system of equations is obtained for the global problem
\begin{equation}\label{eq:globalProblemSystem}
  \begin{aligned}
  \sum_{e=1}^{\numel}\Big\{
    \begin{bmatrix} \mat{A}_{L \hu}^T & \mat{A}_{\hu u} & \mat{A}_{p\hu}^T \end{bmatrix}_{e}
    \begin{Bmatrix} \vect{L}_e \\ \vect{u}_e \\ \vect{p}_e \end{Bmatrix}
   +
    [\mat{A}_{\hu\hu}]_e \, \vect{\hu}_e \Big\}
  &= 
  \sum_{i=e}^{\numel} [\vect{f}_{\hu}]_e, \\
  \vect{1}^T \, [\mat{A}_{p\hu}]_e \vect{\hu}_e &= - \vect{1}^T \, [\vect{f}_p]_e
  %[\vect{f}_{\rho}]_e
  \end{aligned}
\end{equation}

The expressions of the matrices and vectors appearing in\ \eqref{eq:localProblemSystem} and \eqref{eq:globalProblemSystem} are detailed in Appendix\ \ref{ap:Implementation}.

After replacing the solution of the local problem\ \eqref{eq:localProblemSystem} in\ \eqref{eq:globalProblemSystem}, the global problem becomes
\begin{equation}\label{eq:globalProblemSystemFinal}
  \begin{bmatrix}\widehat{\mat{K}} & \mat{G} \\
                   \mat{G}^T           & \mat{0}  \end{bmatrix}
  \begin{Bmatrix}\vect{\hu} \\ \bm{\rho} \end{Bmatrix} 
  =
  \begin{Bmatrix}\vect{\hat{f}}_{\hu} \\ \vect{\hat{f}}_{\rho} \end{Bmatrix},
\end{equation}
where
\begin{subequations}\label{eq:Dglobal-matvect}
\begin{gather}
  \widehat{\mat{K}} {=} \Assem_{e=1}^{\numel}
      \bigl[\begin{array}{@{}c@{\,}c@{\,}c@{\,}c@{}} \mat{A}_{L \hu}^T & \mat{A}_{\hu u} & \mat{A}_{p\hu}^T & \mat{0}\end{array}\bigr]_{e}
  \hspace{-0.5ex}
      \left[\begin{array}{@{}c@{\,}c@{\,}c@{\,}c@{}}
        \mat{A}_{LL}    & \mat{A}_{Lu} & \mat{0}              & \mat{0}                   \\
        \mat{A}_{Lu}^T & \mat{A}_{uu} & \mat{A}_{pu}^T & \mat{0}                   \\
        \mat{0}             & \mat{A}_{pu} & \mat{0}              & \mat{a}_{\rho p}^T \\
        \mat{0}             & \mat{0}          & \mat{a}_{\rho p}& 0                   \\
      \end{array}\right]_{\! e}^{-1}
  \hspace{-0.5ex}
      \begin{bmatrix}
        \mat{A}_{L\hu} \\
        \mat{A}_{u\hu} \\
        \mat{A}_{p\hu} \\
        \mat{0}
      \end{bmatrix}_{\! e}
  \hspace{-0.5ex}
      {+}
      [\mat{A}_{\hu\hu}]_e ,\\
%
%  \mat{G} = \Assem_{e=1}^{\numel}
%      \begin{bmatrix} \mat{A}_{L \hu}^T & \mat{A}_{\hu u} & \mat{A}_{p\hu}^T & \mat{0}\end{bmatrix}_{e}
%      \begin{bmatrix}
%        \mat{A}_{LL}    & \mat{A}_{Lu} & \mat{0}              & \mat{0}                   \\
%        \mat{A}_{Lu}^T & \mat{A}_{uu} & \mat{A}_{pu}^T & \mat{0}                   \\
%        \mat{0}             & \mat{A}_{pu} & \mat{0}              & \mat{a}_{\rho p}^T \\
%        \mat{0}             & \mat{0}          & \mat{a}_{\rho p}& 0                   \\
%      \end{bmatrix}_{\! e}^{-1}
%    \begin{Bmatrix}
%      \vect{0} \\
%      \vect{0} \\
%      \vect{0} \\
%      1
%    \end{Bmatrix}_{\! e} ,
%    \label{eq:eqG} \\
%
  \mat{G}^T = \begin{bmatrix} \vect{1}^T  \, [\mat{A}_{p\hu}]_1 \\ \vect{1}^T  \, [\mat{A}_{p\hu}]_2 \\ \dotsb \\ \vect{1}^T  \, [\mat{A}_{p\hu}]_{\numel}\end{bmatrix}, 
      \label{eq:eqGT} \\
  \vect{\hat{f}}_{\hu}{=}\Assem_{e=1}^{\numel} [\vect{f}_{\hu}]_e
  {-}\bigl[\begin{array}{@{}c@{\,}c@{\,}c@{\,}c@{}} \mat{A}_{L \hu}^T & \mat{A}_{\hu u} & \mat{A}_{p\hu}^T & \mat{0}\end{array}\bigr]_{e}
  \hspace{-0.5ex}
      \left[\begin{array}{@{}c@{\,}c@{\,}c@{\,}c@{}}
        \mat{A}_{LL}    & \mat{A}_{Lu} & \mat{0}              & \mat{0}                   \\
        \mat{A}_{Lu}^T & \mat{A}_{uu} & \mat{A}_{pu}^T & \mat{0}                   \\
        \mat{0}             & \mat{A}_{pu} & \mat{0}              & \mat{a}_{\rho p}^T \\
        \mat{0}             & \mat{0}          & \mat{a}_{\rho p}& 0                   \\
      \end{array}\right]_{\! e}^{-1}
  \hspace{-0.5ex}
      \begin{Bmatrix}
        \vect{f}_L  \\
        \vect{f}_u  \\
        \vect{f}_p \\
        0 
      \end{Bmatrix}_{\! e} , \\
  \vect{\hat{f}}_{\rho}  =   - \begin{bmatrix} \vect{1}^T  \, [\vect{f}_p]_1 \\ \vect{1}^T  \, [\vect{f}_p]_2 \\ \dotsb \\ \vect{1}^T  \, [\vect{f}_p]_{\numel}\end{bmatrix} .
\end{gather}
\end{subequations}

It is worth noticing that the resulting system in\ \eqref{eq:globalProblemSystemFinal} has a saddle-point structure, as it is classical in the discretization of incompressible flow problems (see\ \cite[Section 6.5]{Donea-Huerta}). The proof of the symmetry of the rectangular off-diagonal blocks in\ \eqref{eq:globalProblemSystemFinal} is presented in Appendix\ \ref{ap:Symmetry}.

%________________________________________________________________________
\subsection{Local postprocess of the primal variable}\label{sc:PostProcess}

HDG methods feature the possibility of constructing a superconvergent approximation of the primal variable via a local postprocess performed element-by-element.
Several strategies have been proposed in the literature, either to recover a globally $\Hdiv$-conforming and pointwise solenoidal velocity field, see\ \cite{Jay-CGNPS:11,Nguyen-CNP:10} or simply to improve the approximation of the velocity field (see\ \cite{Nguyen-NPC:10} or\ \cite{RS-SH:16}) and then, derive an error indicator for degree adaptive procedures (see\ \cite{GG-GFH:14} and\ \cite{RS-SH:18}).
The former approach stem from the Brezzi-Douglas-Marini (BDM) projection operator, see\ \cite{brezzi1991mixed}, whereas the latter is inspired by the work of  \cite{Stenberg-90}.

The superconvergent property of the HDG postprocessed solution depends, on the one hand, on the optimal convergence of order $k+1$ of the mixed variable and, on the other hand, on a procedure to resolve the indeterminacy due to rigid body motions.
This is especially critical when the Cauchy stress formulation of incompressible flow problems is considered since a loss of optimal convergence and superconvergence is experienced for low-order polynomial approximations, as noted by \cite{Nguyen-CNP:10}.
Cockburn and coworkers, see\ \cite{Fu-CFS-17,Fu-CF-17,Fu-CF-17b}, proposed the $\bm{M}$-decomposition to enrich the discrete space of the mixed variable. Whereas \cite{Shi-QS-16} suggested an HDG formulation with different polynomial degrees of approximation for primal, mixed and hybrid variables.
Using a pointwise symmetric mixed variable and equal-order of approximation for all the variables, HDG-Voigt formulation has proved its capability to retrieve optimal convergence of the mixed variable and superconvergence of the postprocessed primal one, as shown by \cite{RS-SGKH:18} and \cite{MG-GKSH:18}.

To determine the postprocessed velocity field $\bu_e^\star$ in each element $\Omega_e$, a local problem to be solved element-by-element is devised.
Applying the divergence operator to the first equation in\ \eqref{eq:O-StrongLocal-inde}, for each element $\Omega_e, \ e = 1,\ldots,\numel$, it follows
\begin{subequations}\label{eq:Postprocess-Strong}
\begin{equation}
\grad {\cdot} \left(\sqrt{2\nu} \defo{\bu_e^\star} \right) = - \grad {\cdot} \bL_e ,
\end{equation}
whereas boundary conditions enforcing equilibrated fluxes on $\partial\Omega_e$ are imposed, namely
\begin{equation}%\label{eq:PostprocessBoundary}
\left( \sqrt{2\nu} \defo{\bu_e^\star} \right) \bn = - \bL_e \bn .
\end{equation}
\end{subequations}

The elliptic problem\ \eqref{eq:Postprocess-Strong} admits a solution up to rigid motions, that is $\nsd$ translations (i.e.\ two in 2D and three in 3D) and $\nrr$ rotations (i.e.\ one in 2D and three in 3D).
First, the indeterminacy associated with the $\nsd$ rigid translational modes is handled by means of a constraint on the mean value of the velocity, namely
\begin{equation} \label{eq:PostprocessMean}
( \bu^\star_e , 1 )_{\Omega_e} = ( \bu_e , 1)_{\Omega_e} .
\end{equation} 
Second, a condition on the $\curl$ of the velocity (i.e.\ the vorticity) is introduced to resolve the $\nrr$ rigid rotational modes
\begin{equation} \label{eq:PostprocessRot}
( \grad \times \bu^\star_e , 1 )_{\Omega_e}  =  \langle \bu_D \cdot \btau , 1 \rangle_{\partial \Omega_e \cap \Gamma_D} + \langle \bhu \cdot \btau , 1 \rangle_{\partial \Omega_e \setminus \Gamma_D} ,
\end{equation}
where $\btau$ is the tangential direction to the boundary $\partial\Omega_e$ and the right-hand side follows from Stokes' theorem and the Dirichlet conditions imposed in the local problem\ \eqref{eq:O-StrongLocal-inde}.

As previously done for the local and global problems, the discrete form of the postprocessing procedure is obtained by exploiting Voigt notation to represent the involved symmetric second-order tensors, that is
\begin{equation} \label{eq:postprocessVoigt}
\left\{\begin{aligned}
\gradS^T \bDHalf \gradS \bu_e^\star  &= - \gradS^T \bL_e       &&\text{in $\Omega_e$,}\\	
\bN_e^T \bDHalf \gradS \bu_e^\star &= - \bN_e^T \bL_e        &&\text{on $\partial\Omega_e$.}\\
\end{aligned}\right.
\end{equation}
Similarly, condition\ \eqref{eq:PostprocessRot} is equivalent to
\begin{equation} \label{eq:postprocessRotVoigt}
( \gradW \bu^\star_e , 1 )_{\Omega_e} = \langle \bT \bu_D , 1 \rangle_{\partial \Omega_e \cap \Gamma_D} + \langle \bT \bhu , 1 \rangle_{\partial \Omega_e \setminus \Gamma_D},
\end{equation}
where $\gradW$ and $\bT$ are defined in\ \eqref{eq:curlVoigt} and\ \eqref{eq:tangentVoigt}, repsectively.

Thus, following\ \cite{MG-GKSH:18}, the local postprocess needed to compute a superconvergent velocity in the HDG-Voigt formulation is as follows.
%
%	\Vh_\star(\Omega) & := \left\{ v \in \eltwo(\Omega) : v \vert_{\Omega_e}\in \mathcal{P}^{k+1}(\Omega_e) \;\forall\Omega_e \, , \, e=1,\dotsc ,\numel \right\}, \label{eq:Vstar}
%
First, consider a space for the velocity components with, at least, one degree more than the previous one defined in\ \eqref{eq:spaceScalarElem}, namely:
\begin{equation*}
  \Vh_\star(\Omega)  := \left\{ v \in \eltwo(\Omega) : v \vert_{\Omega_e}\in \mathcal{P}^{k+1}(\Omega_e) \;\forall\Omega_e \, , \, e=1,\dotsc ,\numel \right\},
\end{equation*}
where $\mathcal{P}^{k+1}(\Omega_e)$ denotes the space of polynomial functions of complete degree at most $k+1$ in $\Omega_e$.
Then, for each element $\Omega_e, \ e=1,\dotsc,\numel$, the discrete local postprocess is: given $\bL_e$, $\bu_e$ and $\bhu$ solutions of \eqref{eq:O-WeakLocal-Voigt} and \eqref{eq:O-WeakGlobal-Voigt} with optimal rate of converge $k+1$, find the velocity $\bu_e^\star \in \left[ \Vh_\star(\Omega_e) \right]^{\nsd}$ such that
\begin{equation} \label{eq:postprocessDiscrete}
\left\{\begin{aligned}
- ( \gradS \bv^\star, \bDHalf \gradS \bu_e^\star )_{\Omega_e}  &= ( \gradS \bv^\star, \bL_e )_{\Omega_e} , \\
( \bu^\star_e , 1 )_{\Omega_e} &= ( \bu_e , 1)_{\Omega_e} , \\
( \gradW \bu^\star_e , 1 )_{\Omega_e} &= \langle \bT \bu_D , 1 \rangle_{\partial \Omega_e \cap \Gamma_D} + \langle \bT \bhu , 1 \rangle_{\partial \Omega_e \setminus \Gamma_D} ,
\end{aligned}\right.
\end{equation}
for all $\bv^\star \in \left[ \Vh_\star(\Omega_e) \right]^{\nsd}$.
It is worth emphasizing that in the first equation in\ \eqref{eq:postprocessDiscrete}, the boundary terms appearing after integration by parts are naturally equilibrated by the boundary condition on $\partial\Omega_e$, see\ \eqref{eq:postprocessVoigt}. Moreover, it is also important to note that the right-hand-sides of the last two equations in\ \eqref{eq:postprocessDiscrete} converge with orders larger than $k+1$. Finally, the rate of convergence of $\bu_e^\star$ in diffusion dominated areas is $k+2$.

%________________________________________________________________________
\section{Numerical examples}\label{sc:NumEx}
%________________________________________________________________________
\subsection{Stokes flow}\label{sc:S-NumEx}

The first numerical example involves the solution of the so-called Wang flow\ \cite{wang1991exact} in  $\Omega = [0,1]^2$. This problem has a known analytical velocity field given by
\begin{equation}
\bu(\bx) =
\begin{Bmatrix} 
2 a x_2 - b \lambda \cos(\lambda x_1) \exp\{-\lambda x_2\} \\[1ex] 
b \lambda \sin(\lambda x_1) \exp\{-\lambda x_2\}
\end{Bmatrix}.
\end{equation}
The value of the parameters is selected as $a=b=1$ and $\lambda=10$ and the kinematic viscosity $\nu$ is taken equal to 1. The pressure is selected to be uniformly zero in the domain. 

Neumann boundary conditions are imposed on the bottom part of the domain, $\Gamma_N = \{(x_1,x_2) \in \Omega \; | \; x_2=0\}$, whereas Dirichlet boundary conditions are imposed on $\Gamma_D = \partial \Omega \setminus \Gamma_N$.

\begin{figure}[!tb]
	\centering
	\includegraphics[width=0.45\textwidth]{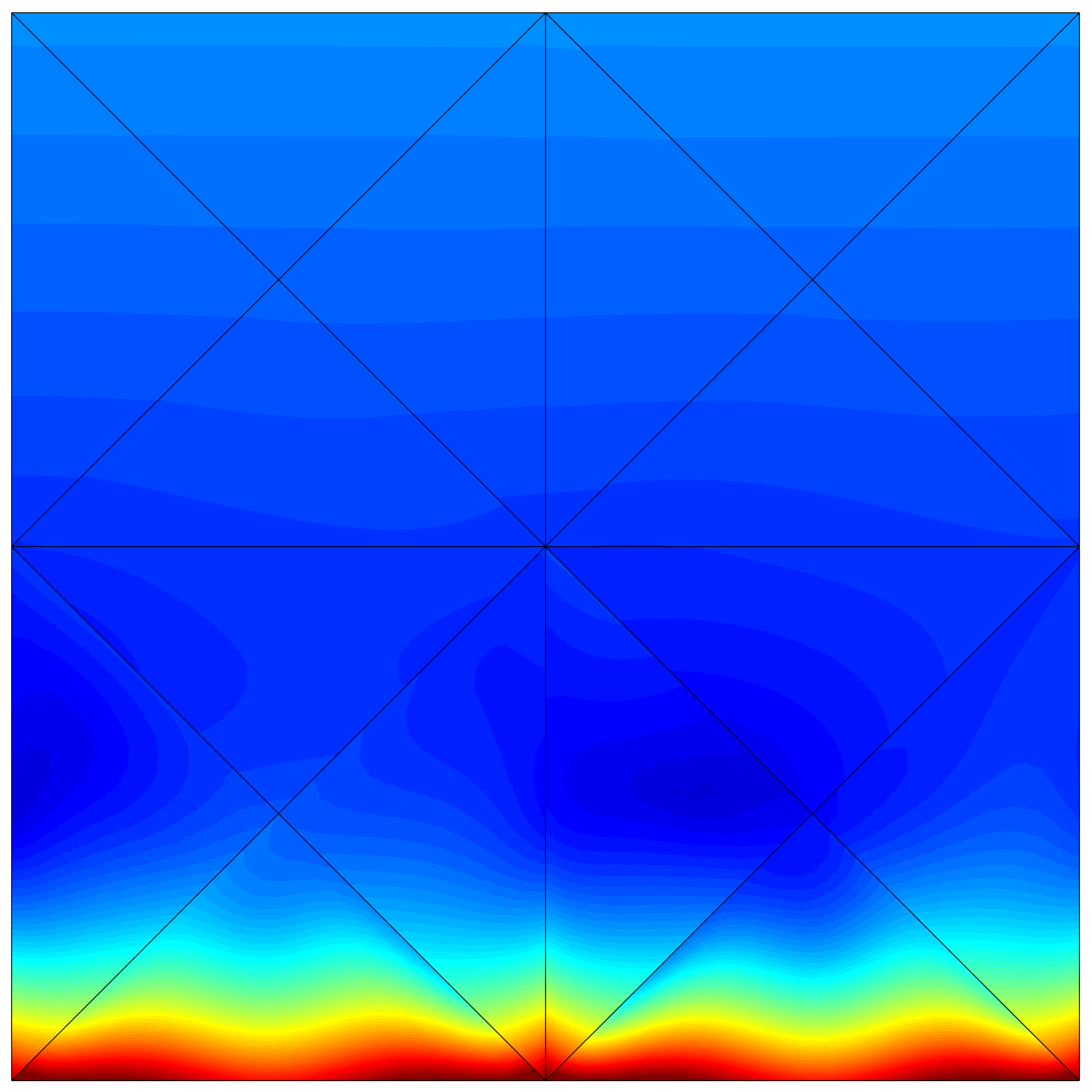}
	\includegraphics[width=0.45\textwidth]{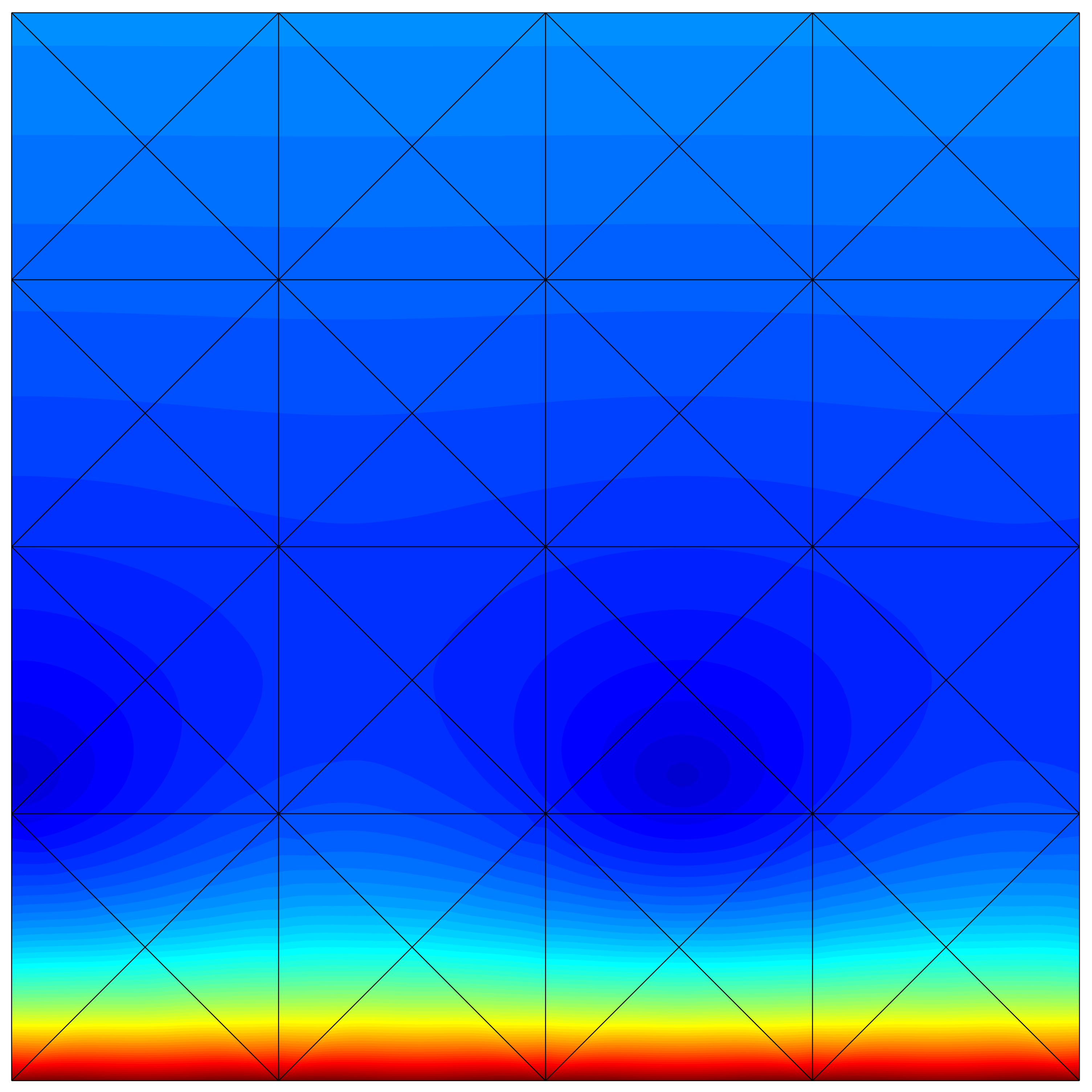}
	\caption{Stokes equations - Wang flow: magnitude of the velocity field computed with the HDG-Voigt method using a degree of approximation $k=3$ on two triangular meshes.}
	\label{fig:StokesWang_sol}
\end{figure}

A sequence of uniform triangular meshes is generated. Figure\ \ref{fig:StokesWang_sol} shows the magnitude of the velocity computed in the first two uniform triangular meshes with a degree of approximation $k=3$.
The results clearly illustrate the gain in accuracy provided by a single level of mesh refinement.

\begin{figure}[!tb]
	\centering
	\includegraphics[width=0.49\textwidth]{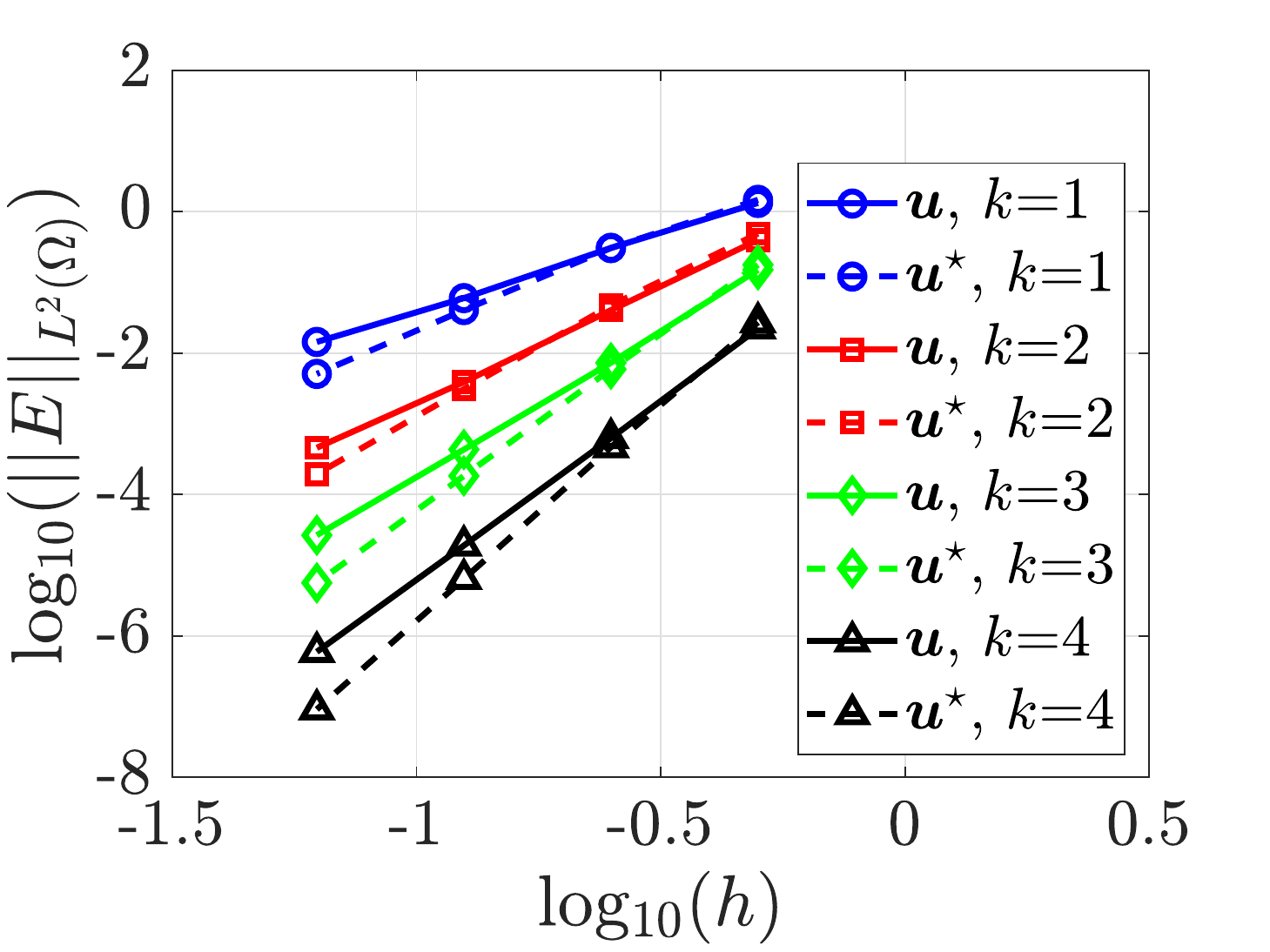}
	\includegraphics[width=0.49\textwidth]{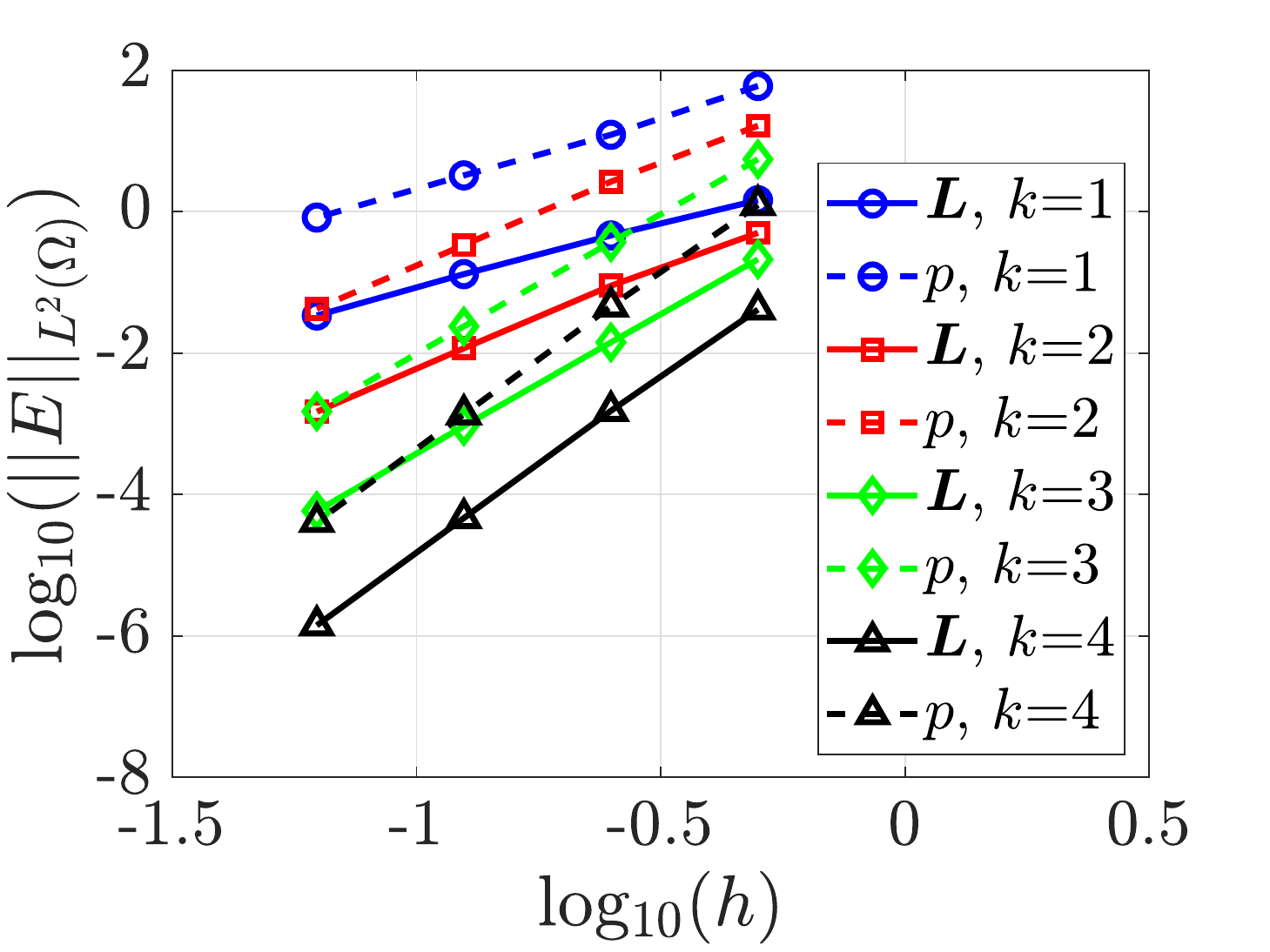}
	\caption{Stokes equations - Wang flow: error of the velocity ($\bu$), postprocessed velocity ($\bu^\star$), scaled strain-rate tensor ($\bL$) and pressure ($p$) in the \eltwo($\Omega$) norm as a function of the characteristic element size $h$.}
	\label{fig:StokesWang_hConv}
\end{figure}

To validate the implementation of the HDG-Voigt formulation, a mesh convergence study is performed. Figure\ \ref{fig:StokesWang_hConv} shows the relative error in the \eltwo($\Omega$) norm as a function of the characteristic element size $h$ for the velocity ($\bu$), postprocessed velocity ($\bu^\star$), scaled strain-rate tensor ($\bL$) and pressure ($p$).
Optimal rate of convergence, equal to $k+1$, is observed for the velocity, scaled strain-rate tensor and pressure and optimal rate, equal to $k+2$, is observed for the postprocessed velocity. In all the examples, the stabilization parameter is taken as $\tau^d=\kappa \nu/\ell$, with $\kappa=3$.

It is worth emphasizing that the optimal rate of convergence is observed even for linear and quadratic approximations. This is in contrast with the suboptimal convergence of the mixed variable, and a loss of superconvergence of the postprocessed velocity, using low-order approximations with the classical HDG equal-order approximation for the Cauchy formulation as shown by\ \cite{Nguyen-CNP:10}.

\begin{figure}[!tb]
	\centering
	\includegraphics[width=0.49\textwidth]{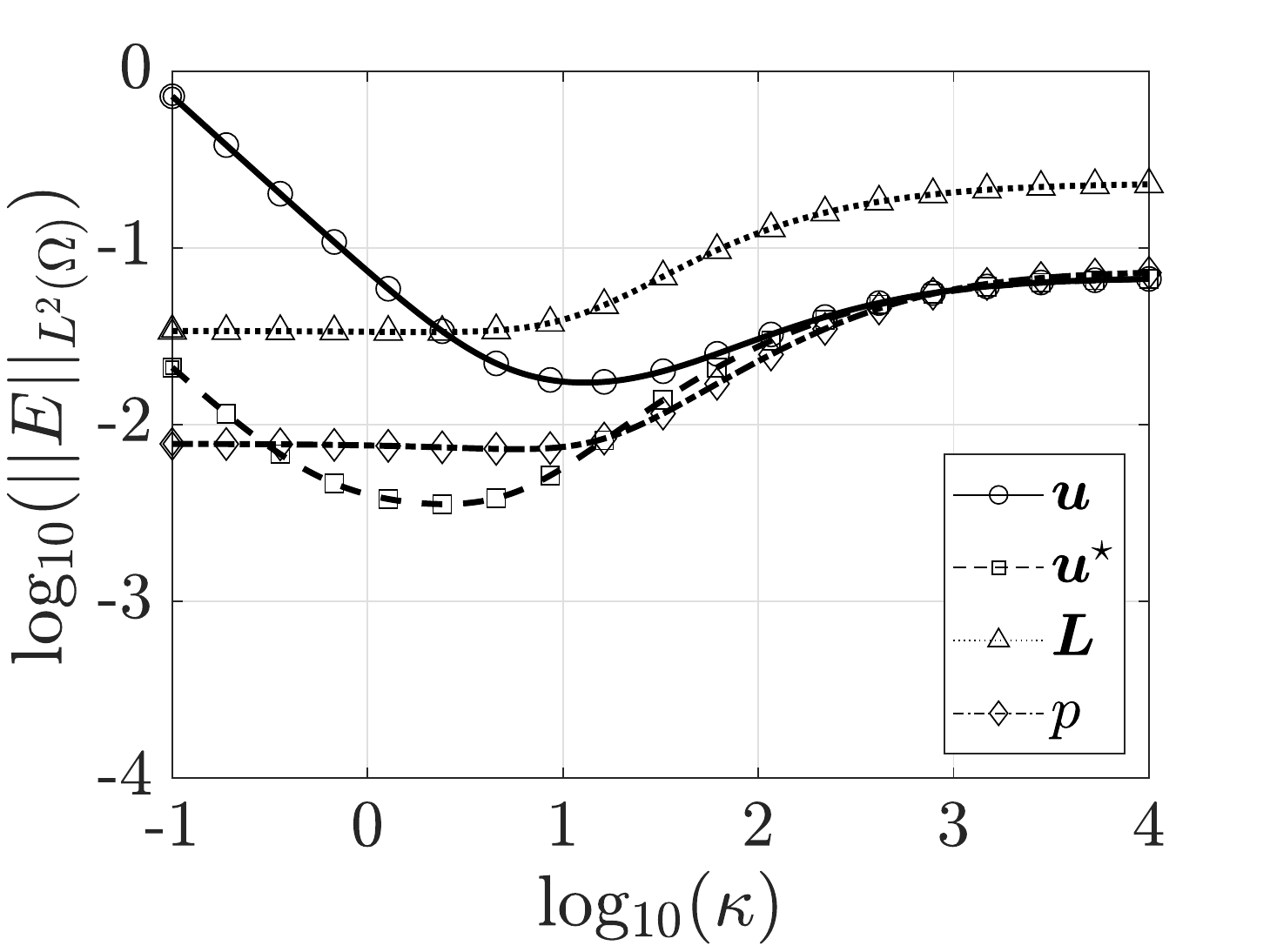}
	\includegraphics[width=0.49\textwidth]{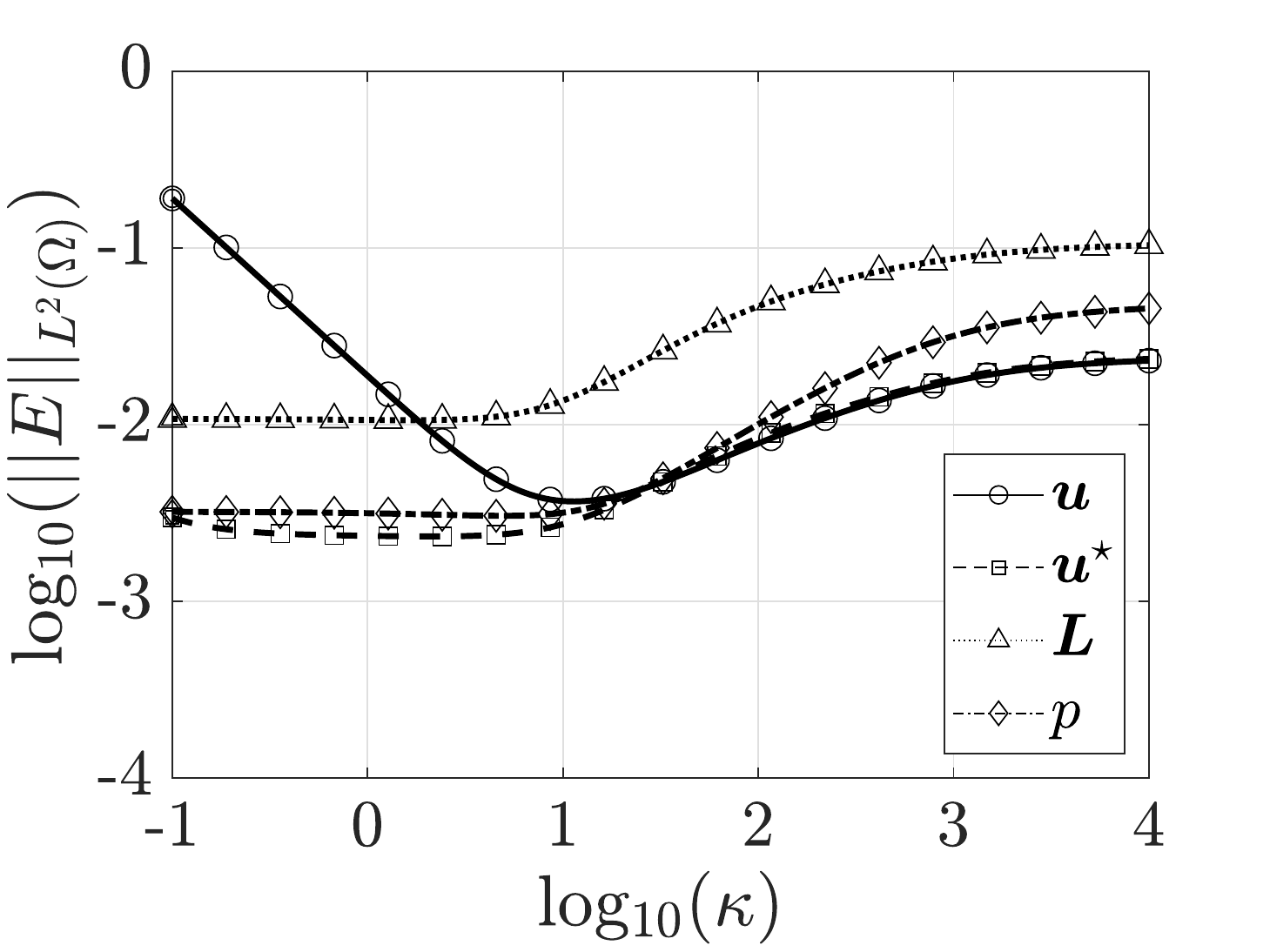}
	\caption{Stokes equations - Wang flow: error of the velocity ($\bu$), postprocessed velocity ($\bu^\star$), scaled strain-rate tensor ($\bL$) and pressure ($p$) in the \eltwo($\Omega$) norm as a function of the parameter $\kappa$ used to define $\tau^d$. Third mesh with $k=1$ (left) and second mesh with $k=2$ (right).}
	\label{fig:StokesWang_tau}
\end{figure}

The influence of the stabilization parameter $\kappa$ in the accuracy of the HDG-Voigt formulation is studied numerically. Figure\ \ref{fig:StokesWang_tau} shows the relative error in the \eltwo($\Omega$) norm, as a function of the parameter $\kappa$ used to define $\tau^d=\kappa \nu/\ell$. The error for the velocity, postprocessed velocity, scaled strain-rate tensor and pressure is considered.
The results are displayed for two different levels of mesh refinement and two different degrees of approximation. In both cases, it can be observed that there is an optimal value of the stabilization parameter that provides the maximum accuracy for the velocity and the postprocessed velocity, corresponding to $\kappa \in [1,10]$. The error for the scaled strain-rate tensor and the pressure is less sensitive to the value of $\kappa$. It can be observed that low values of the stabilization parameter substantially decrease the accuracy of the velocity field whereas large values produce less accurate results for all variables.

For more numerical examples in two and three dimensions and using different element types, the reader is referred to\ \cite{MG-GKSH:18}, where the HDG-Voigt formulation with a strong imposition of the symmetry of the stress tensor via Voigt notation was originally introduced.

%________________________________________________________________________
\subsection{Oseen flow}\label{sc:O-NumEx}

The numerical study of the Oseen problem involves the solution of the so-called Kovasznay flow \cite{kovasznay1948laminar} in $\Omega = [0,1]^2$. This problem has a known analytical solution given by
\begin{equation}
  \hspace{-0.5ex}
  \bu(\bx) {=} \begin{Bmatrix} 
                  1 {-} \exp(2\lambda x_1) \cos\bigl((4x_2{-}1)\pi \bigr)  \\[0.5ex]
                  (\lambda/2\pi) \exp(2\lambda x_1) \sin\bigl((4x_2{-}1)\pi \bigr)
                  \end{Bmatrix}\hspace{-0.25ex}, \,
  p(\bx) {=} {-} \frac{1}{2} \exp(4 \lambda x_1) {+} C,
\end{equation}
where $\lambda := Re/2 {-} \sqrt{Re^2/4 {+} 4\pi^2}$, $C {=} \bigl[ 1 {+} \exp(4 \lambda) {-} (1/2\lambda)\bigl(1 {-} \exp(4 \lambda)\bigr) \bigr]/8$. The Reynolds number is taken as $Re=100$ and the convection velocity field $\ba$ is taken as the exact velocity.

As done in the previous example, Neumann boundary conditions are imposed on the bottom part of the domain, $\Gamma_N = \{(x_1,x_2) \in \Omega \; | \; x_2=0\}$, whereas Dirichlet boundary conditions are imposed on $\Gamma_D = \partial \Omega \setminus \Gamma_N$.

\begin{figure}[!tb]
	\centering
	\includegraphics[width=0.45\textwidth]{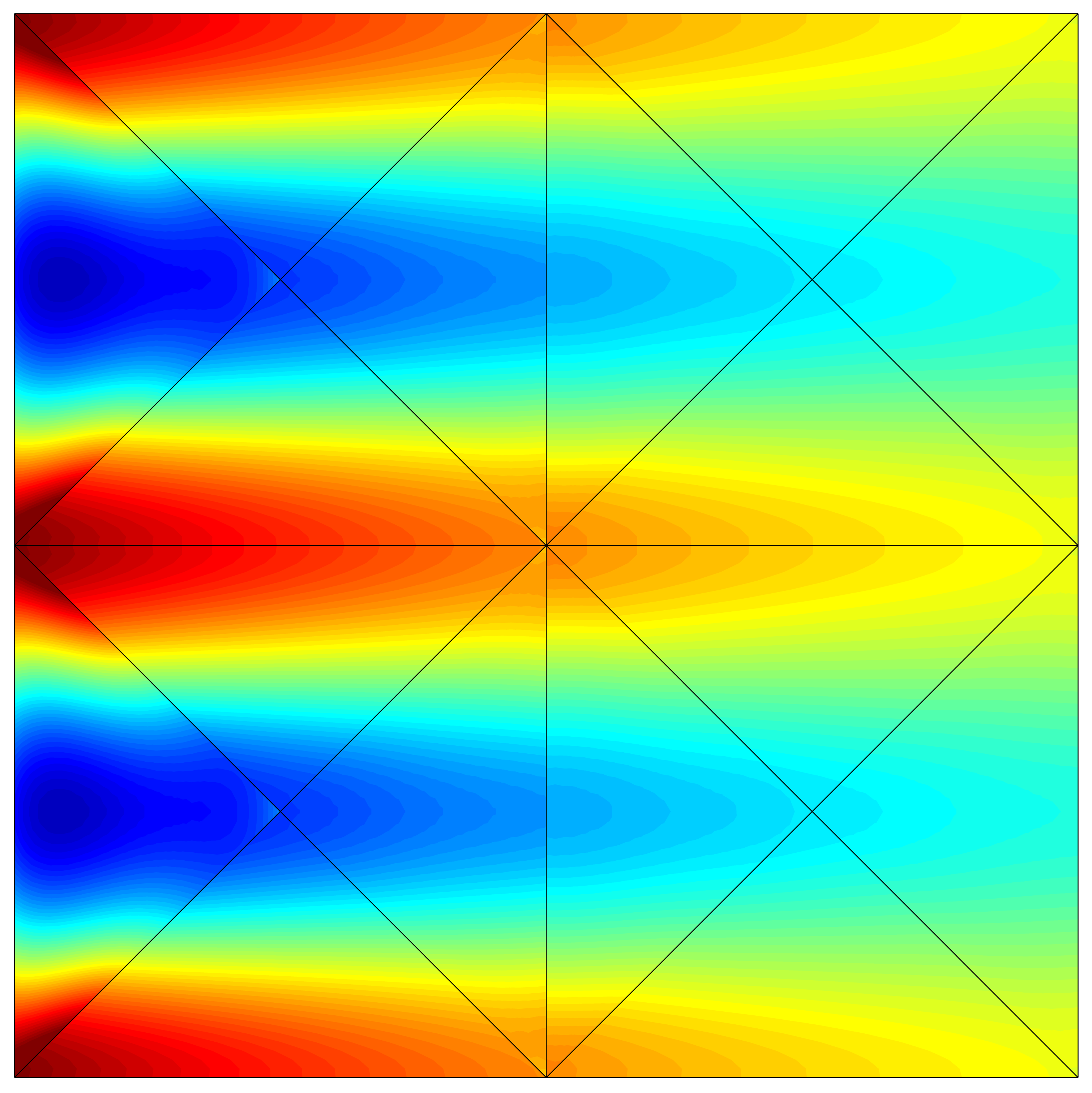}
	\includegraphics[width=0.45\textwidth]{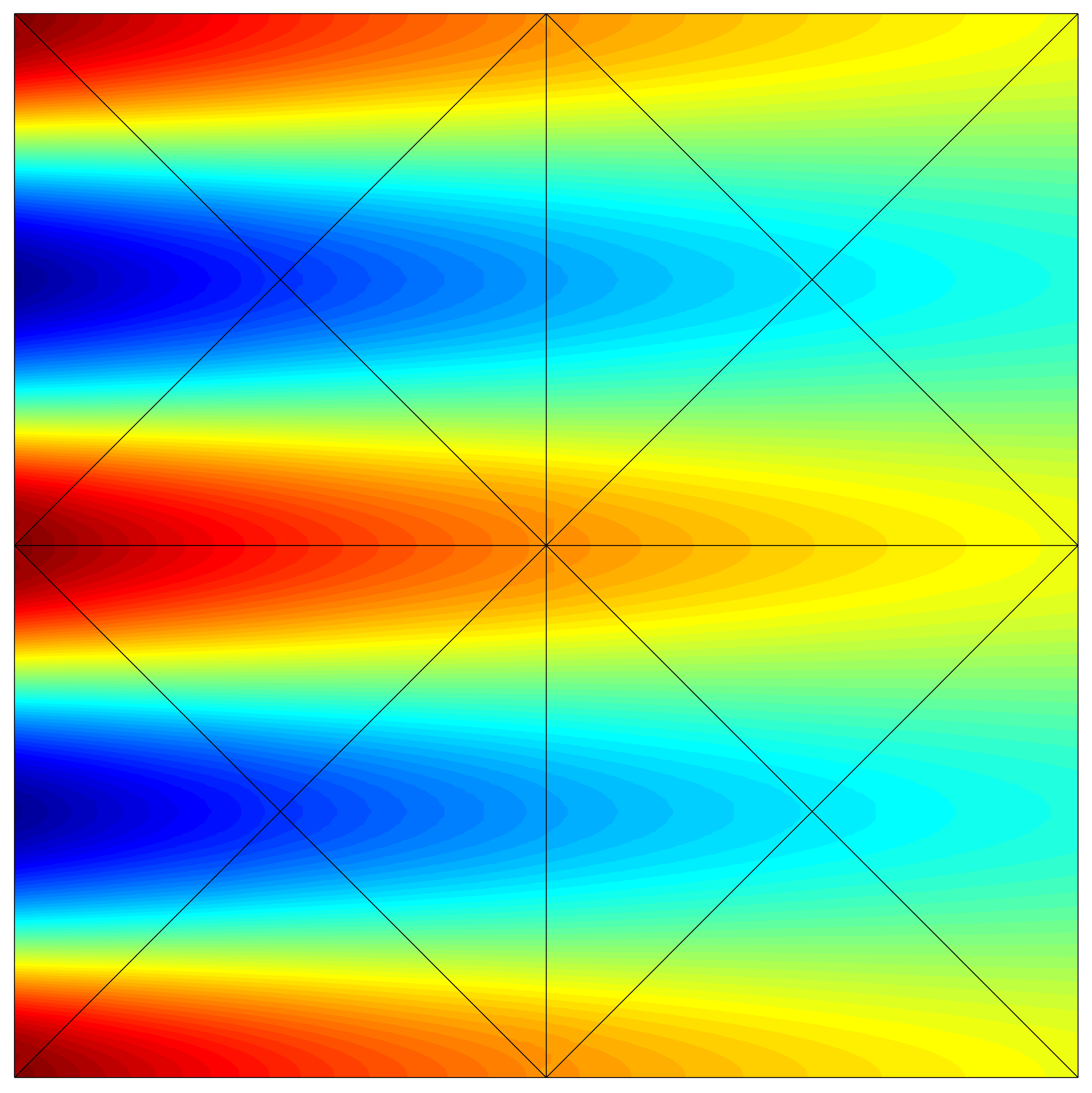}
	\caption{Oseen equations - Kovasznay flow: magnitude of the velocity field computed with the HDG-Voigt method using a degree of approximation $k=4$ (left) and postprocessed velocity (right).}
	\label{fig:OssenKowaznay_sol}
\end{figure}

Figure\ \ref{fig:OssenKowaznay_sol} shows the magnitude of the velocity computed in the first triangular mesh with a degree of approximation $k=4$ and the postprocessed velocity computed by using the strategy described in Section\ \ref{sc:PostProcess}.
The results clearly illustrate the increased accuracy provided by the element-by-element postprocess of the velocity field.

\begin{figure}[!tb]
	\centering
	\includegraphics[width=0.49\textwidth]{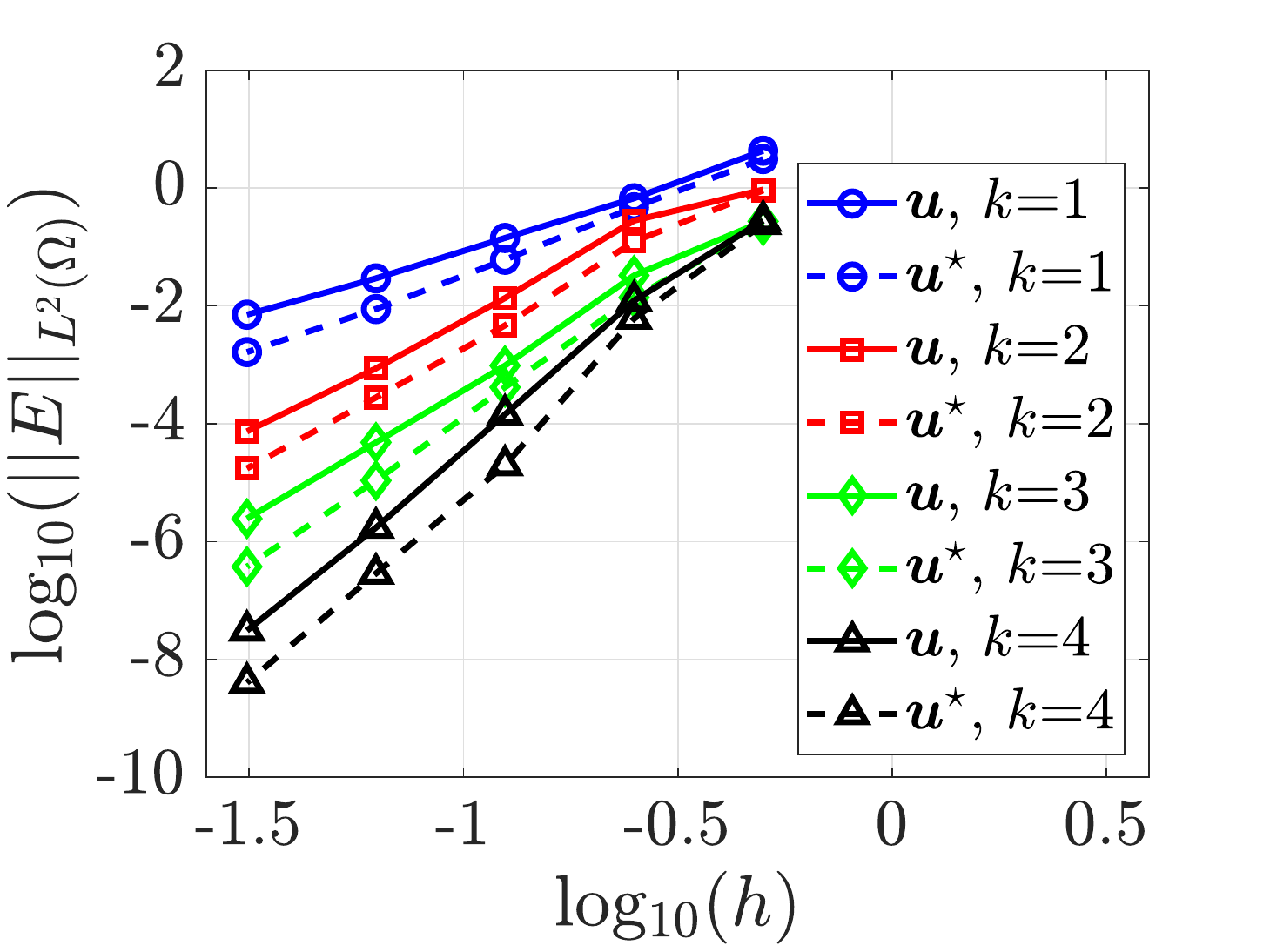}
	\includegraphics[width=0.49\textwidth]{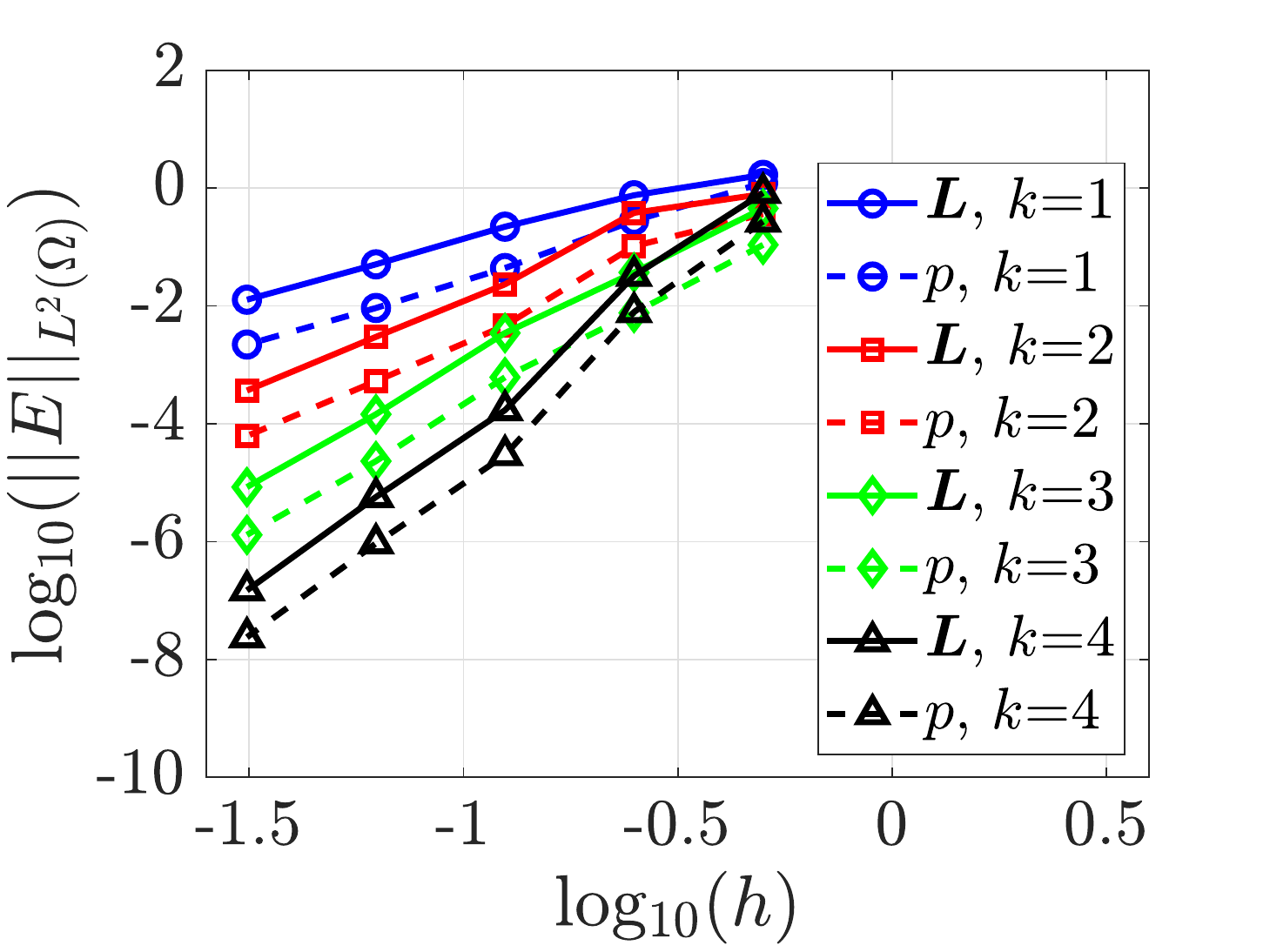}
	\caption{Oseen equations - Kovasznay flow: error of the velocity ($\bu$), postprocessed velocity ($\bu^\star$), scaled strain-rate tensor ($\bL$) and pressure ($p$) in the \eltwo($\Omega$) norm as a function of the characteristic element size $h$.}
	\label{fig:OssenKowaznay_hConv}
\end{figure}

To validate the implementation of the HDG-Voigt formulation, a mesh convergence study is performed as done in the previous example. Figure\ \ref{fig:OssenKowaznay_hConv} shows the relative error in the \eltwo($\Omega$) norm as a function of the characteristic element size $h$ for the velocity ($\bu$), postprocessed velocity ($\bu^\star$), scaled strain-rate tensor ($\bL$) and pressure ($p$).
The optimal rate of convergence is again observed for all the variables. In all the examples, the stabilization for diffusion is defined as $\tau^d=\kappa \nu/\ell$, whereas for convection it holds $\tau^a = \beta \max_{\bx \in \mathcal{N}_{\Omega}} \norm{\ba(\bx)}_2$ with parameters $\kappa=10$ and $\beta=0.02$.

\begin{figure}[!tb]
	\centering
	\includegraphics[width=0.49\textwidth]{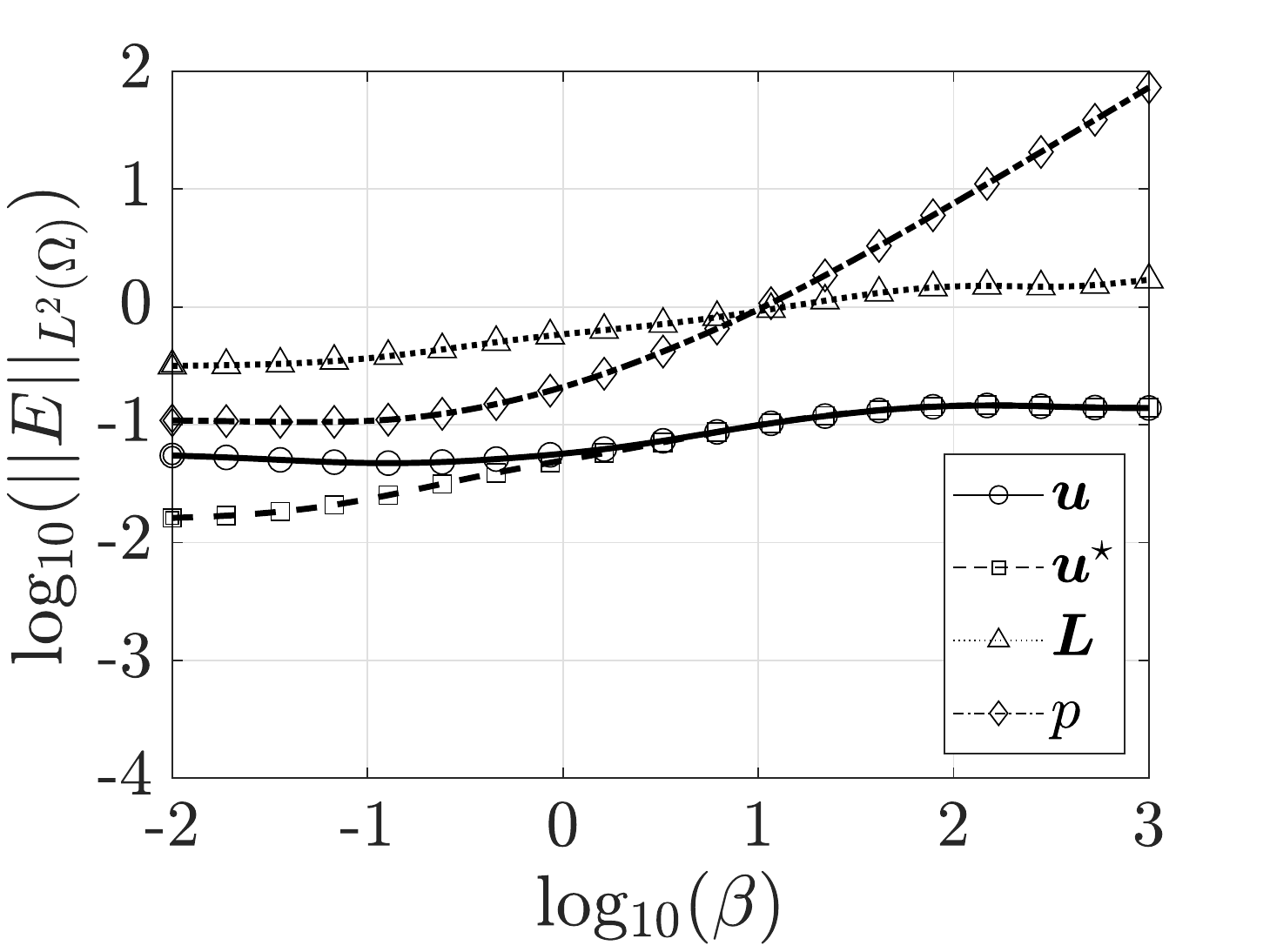}
	\includegraphics[width=0.49\textwidth]{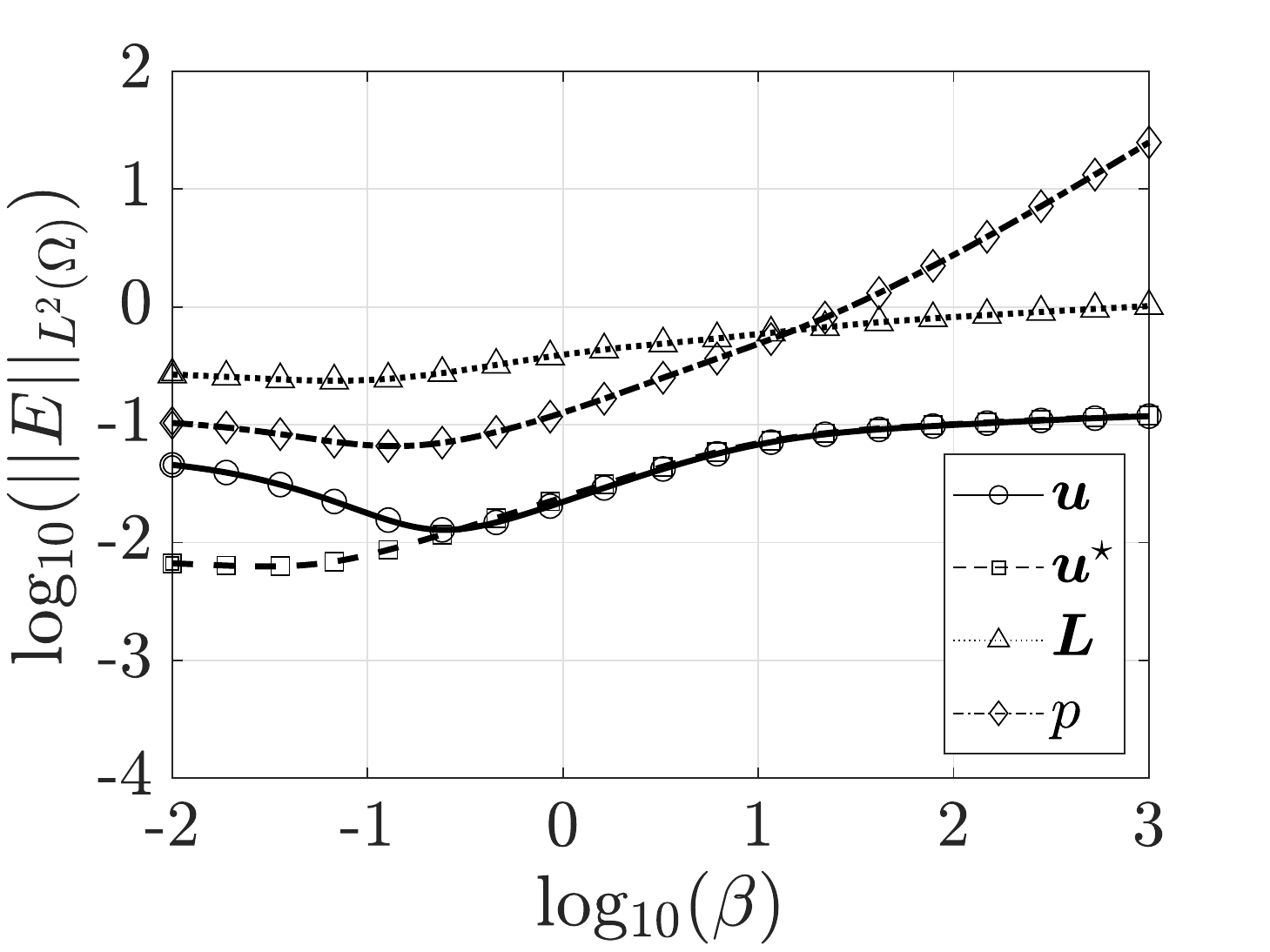}
	\caption{Oseen equations - Kovasznay flow: error of the velocity ($\bu$), postprocessed velocity ($\bu^\star$), scaled strain-rate tensor ($\bL$) and pressure ($p$) in the \eltwo($\Omega$) norm as a function of the parameter $\beta$ used to define $\tau^a = \beta \norm{\ba}_\infty$. Third mesh with $k=1$ (left) and second mesh with $k=2$ (right).}
	\label{fig:OssenKowaznay_tau}
\end{figure}

Finally, the influence of the stabilization parameter $\beta$ is studied numerically. The value of the stabilization parameter $\kappa$ is selected as 10. Figure\ \ref{fig:OssenKowaznay_tau} shows the relative error in the \eltwo($\Omega$) norm, as a function of the parameter $\beta$ used to define the stabilization parameter $\tau^a = \beta \max_{\bx \in \mathcal{N}_{\Omega}} \norm{\ba(\bx)}_2$. The error for the velocity, postprocessed velocity, scaled strain-rate tensor and pressure is considered.
The results are displayed for the same two levels of mesh refinement and degrees of approximation utilized for the sensitivity analysis in the Stokes flow problem. It can be observed that there is an optimal value of the stabilization parameter that provides the maximum accuracy for the velocity and the postprocessed velocity, corresponding to $\beta \in [0.01,0.1]$. The pressure is the variable that shows a higher dependence in terms of the selected value of the stabilization parameter $\beta$. In particular, large values of $\beta$ lead to sizeable errors in the pressure.

%________________________________________________________________________
\subsection{Navier-Stokes flow}\label{sc:NS-NumEx}

The solution of the nonlinear incompressible Navier-Stokes equations is considered. The HDG-Voigt formulation follows the rationale presented in Section\ \ref{sc:Oseen} with $\ba$ being $\bu$ and $\bha$ being $\bhu$. 
The resulting nonlinear local and global problems are solved by using a standard Newton-Raphson linearization. 
It is worth noticing that by substituting $\bu$ for $\ba$ and $\bhu$ for $\bha$ in \eqref{eq:tauConv}, the resulting stabilization parameter $\tau^a$ is now a function of the unknown primal, $\tau^a = \tau^a(\bu)$ or hybrid, $\tau^a = \tau^a(\bhu)$, variable. In order to avoid introducing an additional nonlinearity in the problem, the stabilization parameter $\tau^a$ is thus evaluated in the previous iteration of the Newton-Raphson algorithm.
For the following numerical experiments, at step $r+1$ of the Newton-Raphson method, the definition of $\tau^a = \beta \max_{\bx \in \mathcal{N}_{\Omega}} \norm{\bu^r(\bx)}_2$ is considered, $\bu^r$ being the velocity field computed at iteration $r$ of the nonlinear procedure.

%________________________________________________________________________
\subsubsection{Navier-Stokes Kovasznay flow}\label{sc:NS-Kovasznay}
First, the Kovasznay flow described in Section\ \ref{sc:O-NumEx} is considered with the same boundary conditions and using the same set of meshes used to solve the Ossen equations.  
\begin{figure}[!tb]
	\centering
	\includegraphics[width=0.45\textwidth]{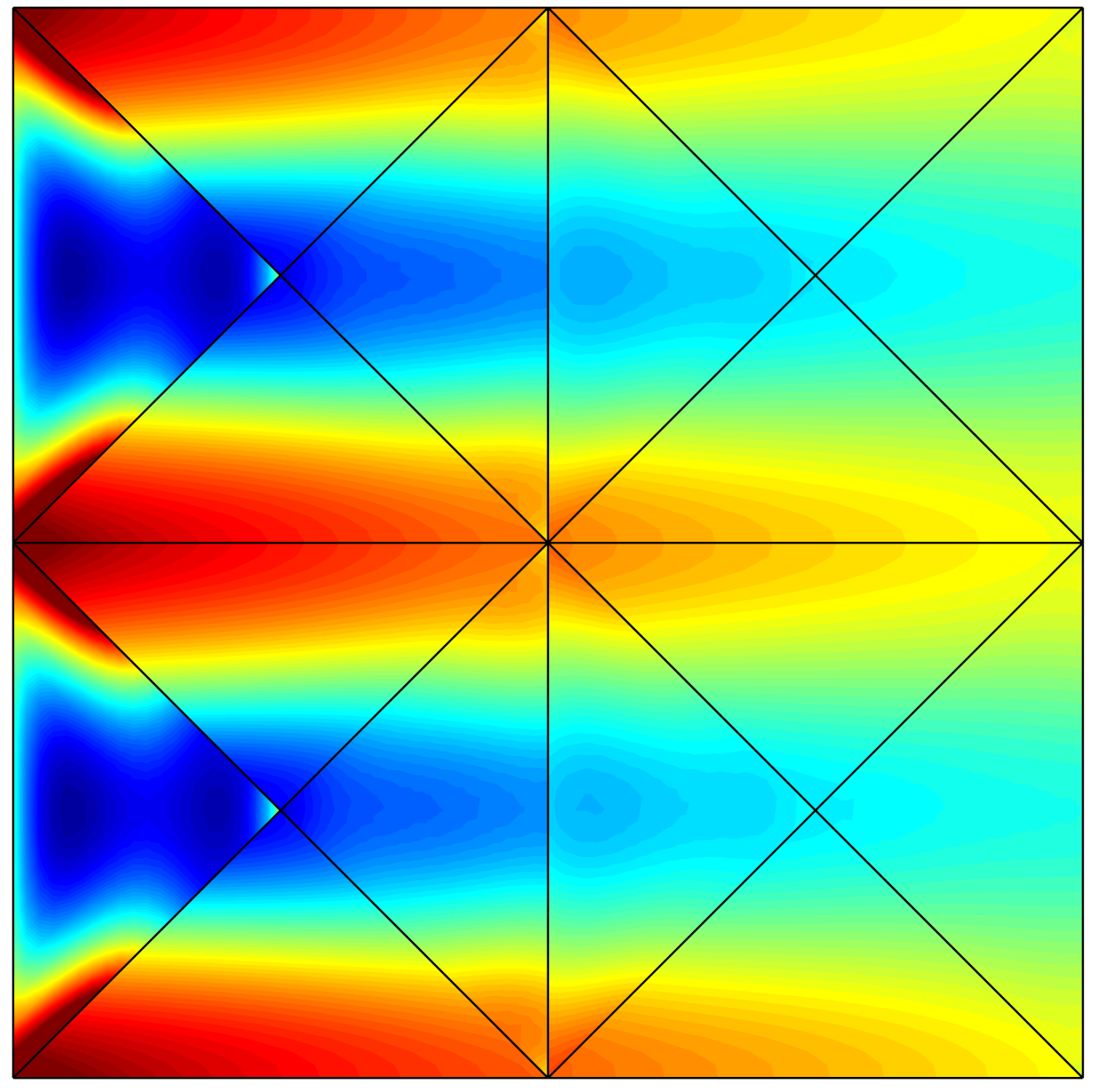}
	\includegraphics[width=0.45\textwidth]{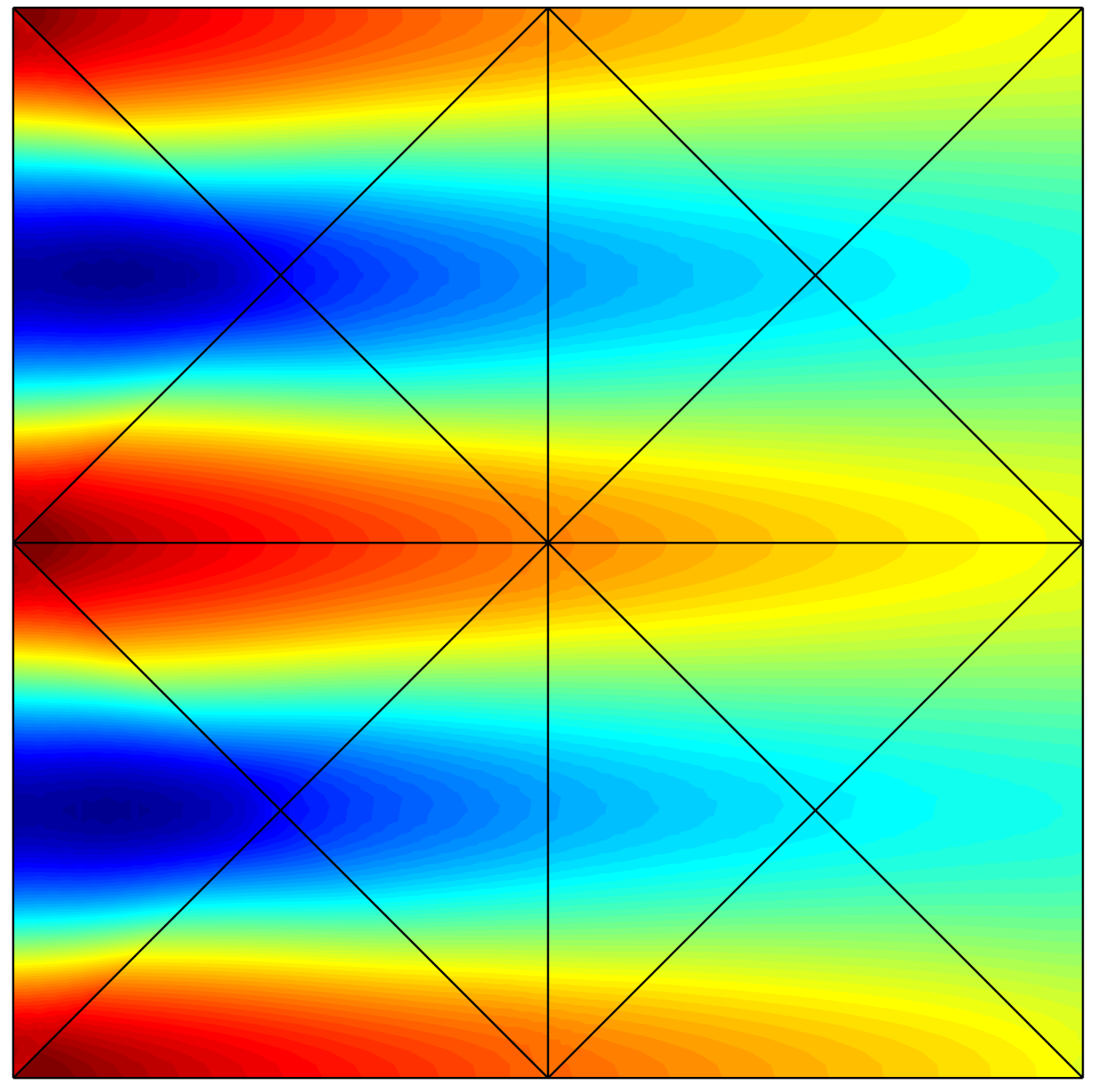}
	\caption{Navier-Stokes equations - Kovasznay flow: magnitude of the velocity field computed with the HDG-Voigt method using a degree of approximation $k=4$ (left) and postprocessed velocity (right).}
	\label{fig:NS-Kowaznay_sol}
\end{figure}
Figure\ \ref{fig:NS-Kowaznay_sol} shows the magnitude of the velocity computed in the first triangular mesh with a degree of approximation $k=4$ and the postprocessed velocity computed by using the strategy described in Section\ \ref{sc:PostProcess}. 
A visual comparison of the results between Navier-Stokes and Ossen solutions, computed using the same coarse mesh, illustrates the extra level of difficulty induced by the nonlinearity of the Navier-Stokes equations.

\begin{figure}[!tb]
	\centering
	\includegraphics[width=0.49\textwidth]{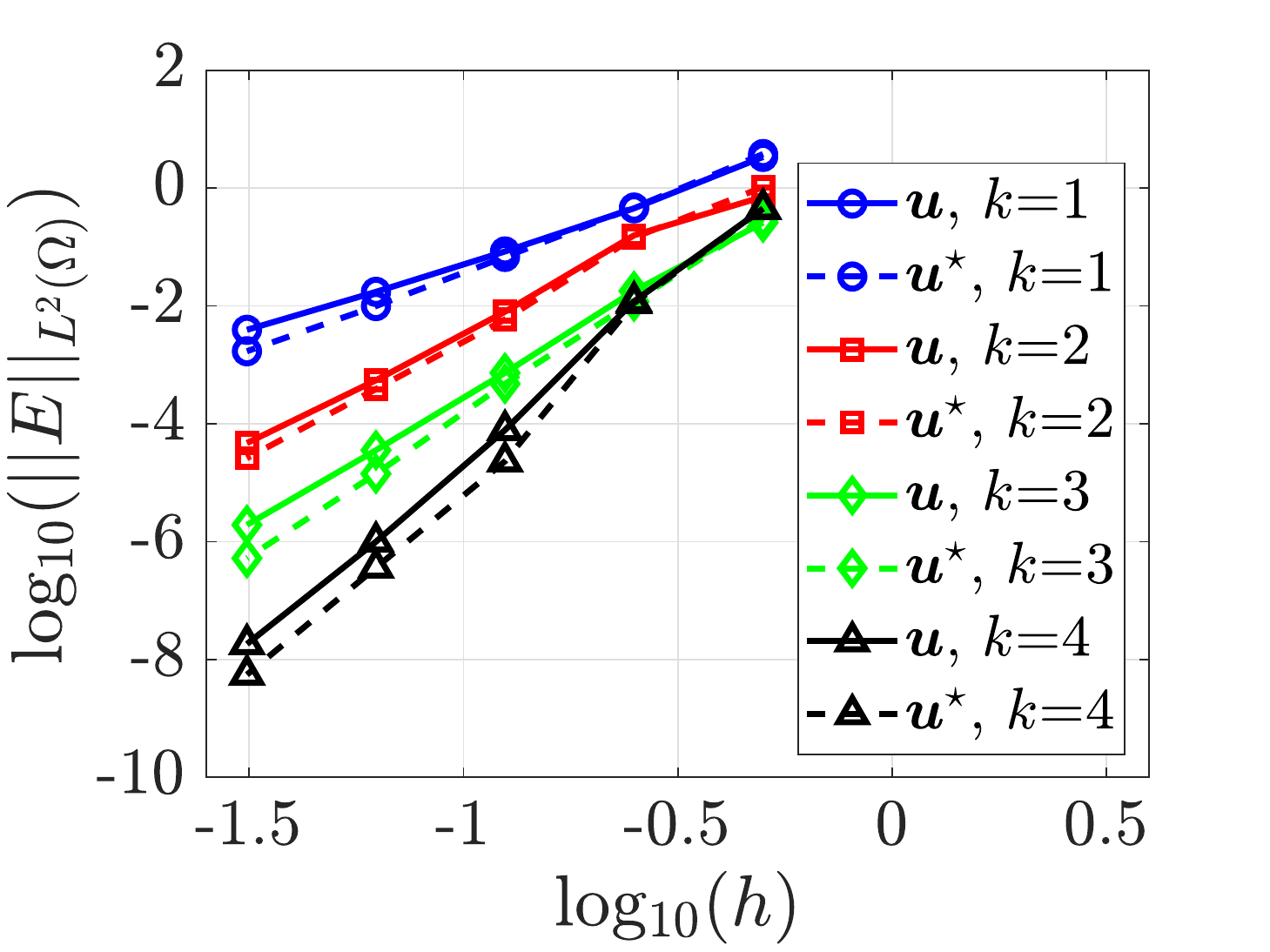}
	\includegraphics[width=0.49\textwidth]{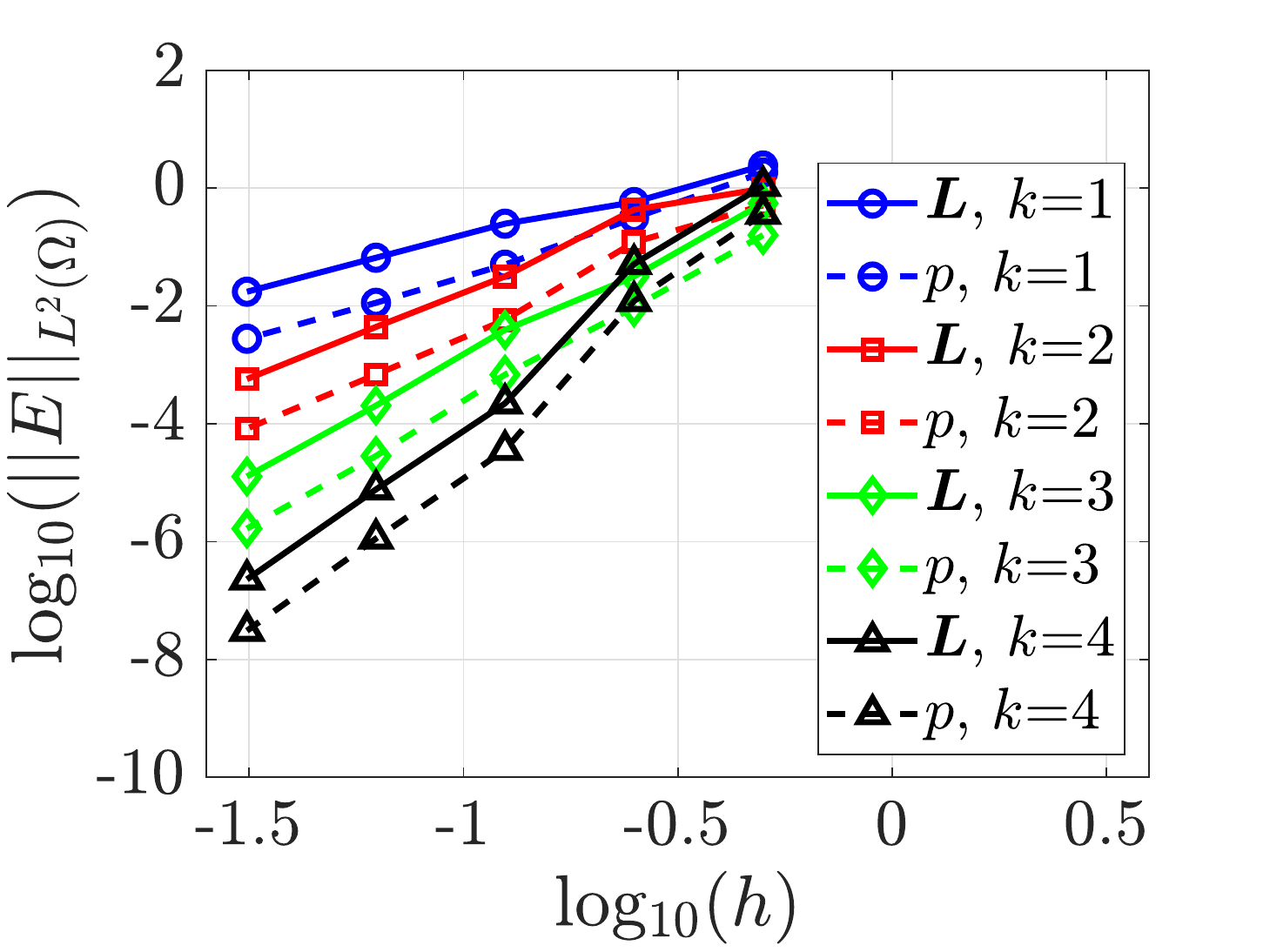}
	\caption{Navier-Stokes equations - Kovasznay flow: error of the velocity ($\bu$), postprocessed velocity ($\bu^\star$), scaled strain-rate tensor ($\bL$) and pressure ($p$) in the \eltwo($\Omega$) norm as a function of the characteristic element size $h$.}
	\label{fig:NS-Kowaznay_hConv}
\end{figure}

To validate the implementation of the HDG-Voigt formulation, a mesh convergence study is performed as done in the previous examples. Figure\ \ref{fig:NS-Kowaznay_hConv} shows the relative error in the \eltwo($\Omega$) norm as a function of the characteristic element size $h$ for the velocity ($\bu$), postprocessed velocity ($\bu^\star$), scaled strain-rate tensor ($\bL$) and pressure ($p$).
The optimal rate of convergence is again observed for all the variables, with a slightly higher rate for the primal variables when $k=4$. In all the examples, the stabilization parameters are taken as $\tau^d=\kappa \nu/\ell$ and $\tau^a = \beta \max_{\bx \in \mathcal{N}_{\Omega}} \norm{\bu^r(\bx)}_2$, with $\kappa=10$ and $\beta=0.1$, $\bu^r$ being the velocity field computed in the previous Newton-Raphson iteration.

%________________________________________________________________________
\subsubsection{Navier-Stokes Poiseuille flow}\label{sc:NS-Poiseuille}
The following example considers another well-known problem with analytical solution to illustrate the effect of the boundary condition imposed on an outflow part of the boundary, as discussed in Section\ \ref{sc:ProbSta}.  The Poiseuille flow in a rectangular channel $\Omega=[0,10]\times[0,1]$ is considered, with Dirichlet boundary conditions on the left, top and bottom part of the boundary and different boundary conditions, described in Remark\ \ref{rm:Outflow}, on the right-end of the boundary. 
The exact solution is given by
\begin{equation}
  \bu(\bx) = \begin{Bmatrix} 
                  4 V x_2(1-x_2)  \\[0.25ex] 
                  0
                  \end{Bmatrix} , \quad
                  p(\bx) = -8\nu V x_1 + C ,
\end{equation}
where $V$ is the maximum value of the velocity profile, achieved at the centerline of the channel, that is for $x_2 = 1/2$, and $C = 80 \nu V$.

\begin{figure}[!tb]
	\centering
	\includegraphics[width=\textwidth]{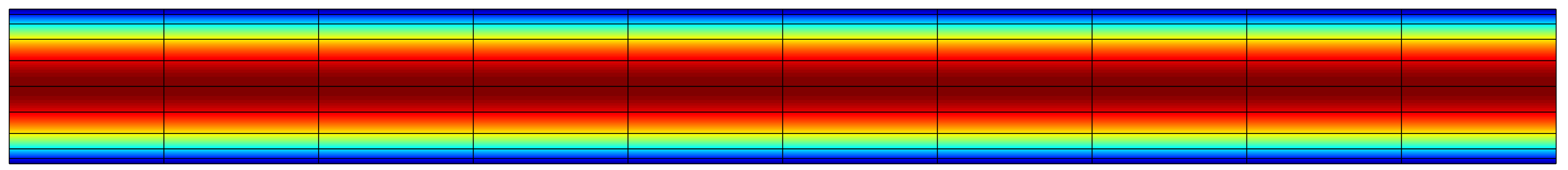}
	\caption{Navier-Stokes equations - Poiseuille flow: magnitude of the velocity field computed with the HDG-Voigt method using a degree of approximation $k=2$ with the outflow boundary condition in\ \eqref{eq:outCorrect}.}
	\label{fig:NS-Poiseuille_solBC3}
\end{figure}

The solution computed with a structured quadrilateral mesh with $10 \times 10$ biquadratic elements is shown in Figure\ \ref{fig:NS-Poiseuille_solBC3}, by imposing the outflow boundary condition given in\ \eqref{eq:outCorrect}. 
As expected, the computed solution reproduces the exact solution (with machine accuracy) as the latter belongs to the polynomial space used to define the functional approximation. 
The error of the velocity measured in the $\mathcal{L}^\infty(\Omega)$ norm is $3.6 \times 10^{-13}$. 
In contrast, when the homogeneous Neumann boundary condition given in\ \eqref{eq:outNeumann} is utilized, the error of the velocity measured in the $\mathcal{L}^\infty(\Omega)$ norm is $1.11$. Finally, when the traction-free boundary condition in\ \eqref{eq:outTraction} is imposed on the right part of the boundary, the error of the velocity measured in the $\mathcal{L}^\infty(\Omega)$ norm is $1.29$. 

\begin{figure}[!tb]
	\centering
	\includegraphics[width=0.325\textwidth]{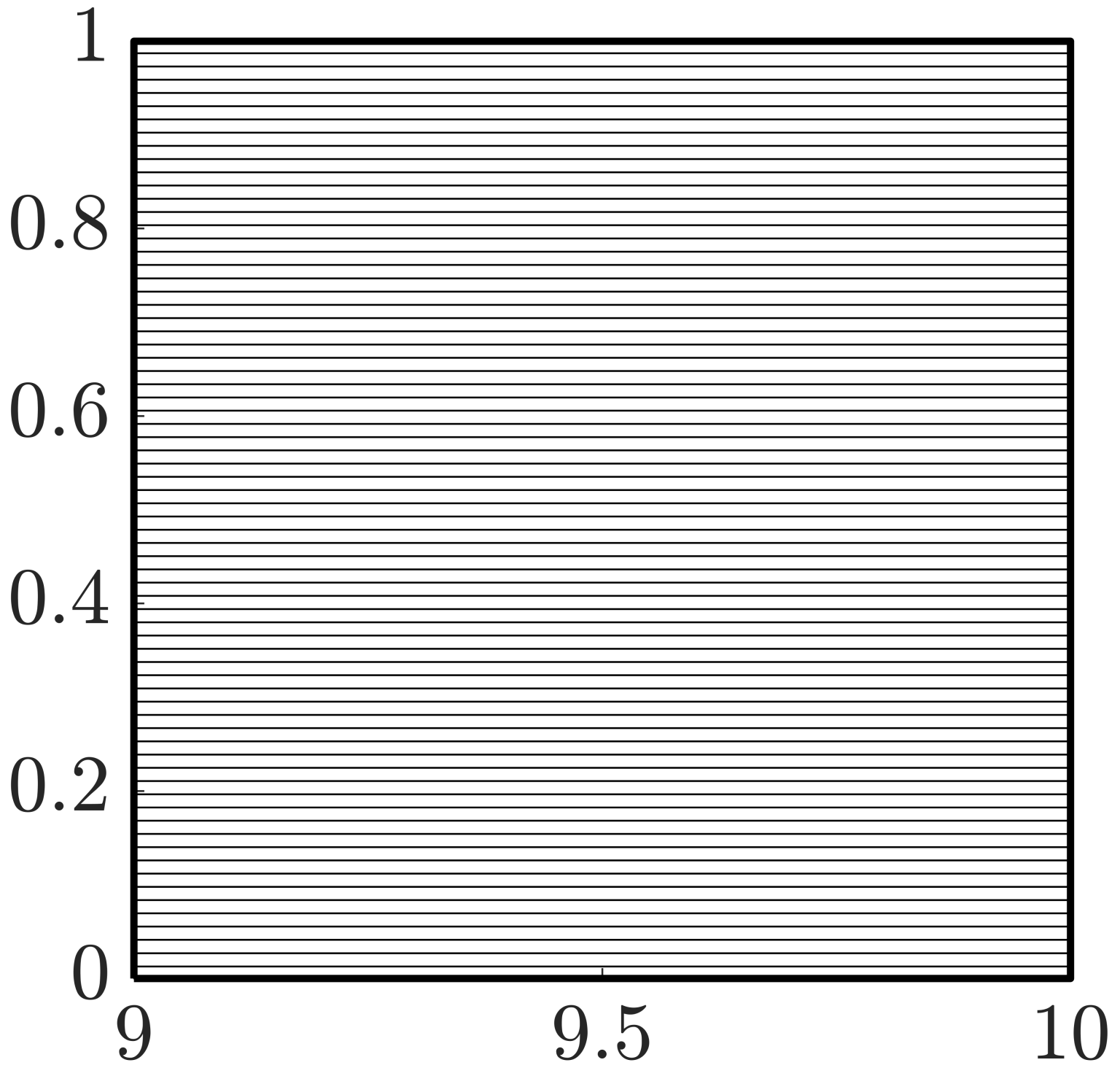}
	\includegraphics[width=0.325\textwidth]{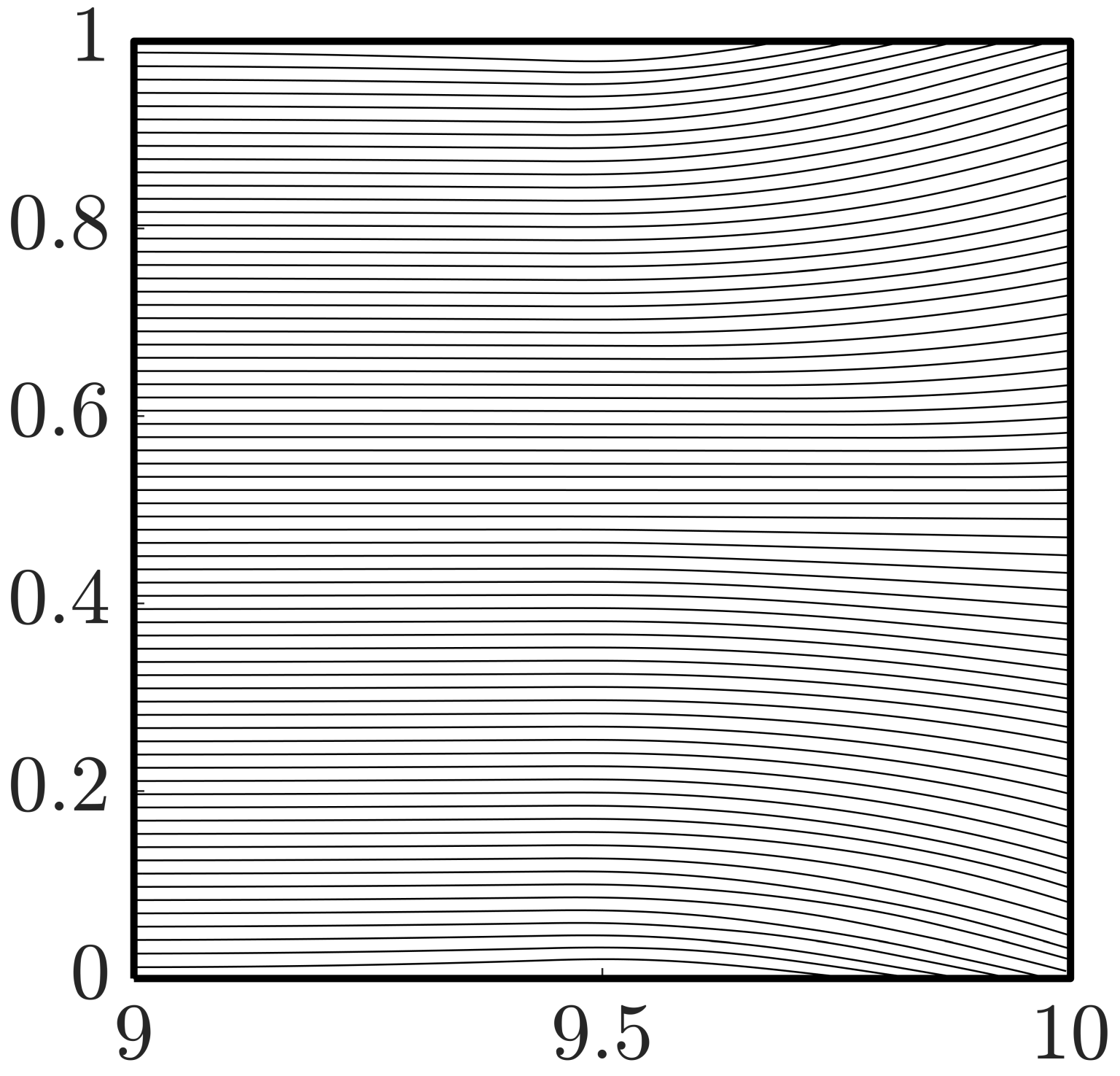}
	\includegraphics[width=0.325\textwidth]{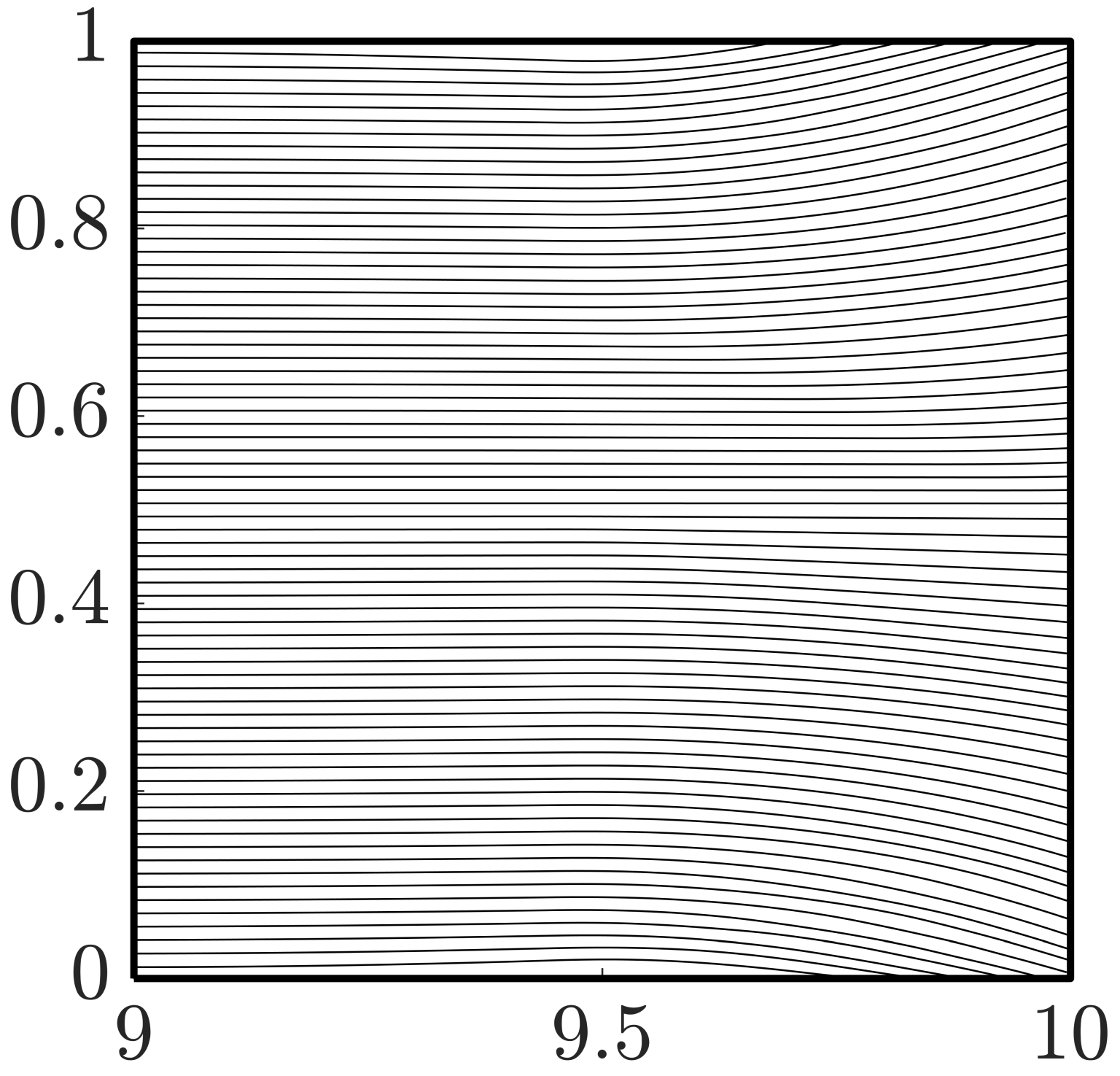}
	\caption{Navier-Stokes equations - Poiseuille flow: magnitude of the velocity field and streamlines $k=2$ with the outflow boundary condition (left) homogeneous Neumann boundary condition in\ \eqref{eq:outNeumann} (middle) and traction-free boundary condition in\ \eqref{eq:outTraction} (right).}
	\label{fig:NS-Poiseuille_solBC}
\end{figure}

A detailed view of the solution near the outflow boundary for the three different boundary conditions considered is shown in Figure\ \ref{fig:NS-Poiseuille_solBC}. 
The figure displays the streamlines. It can be clearly observed that, using the outflow boundary condition, the streamlines are parallel to the $x$-axis and the motion of a fluid between two parallel infinite plates is correctly reproduced. On the contrary, with the homogeneous Neumann boundary condition and the traction-free boundary condition the isolines present non-physical artifacts.

%________________________________________________________________________
\subsubsection{Backward facing step}\label{sc:NS-BFS}
The last example considers another well-known test case for the incompressible Navier-Stokes equations, the so-called backward facing step. This problem is traditionally employed to test the ability of a numerical scheme to capture the recirculation zones and position of the reattachment point\ \cite{armaly1983experimental,erturk2008numerical}. 

\begin{figure}[!tb]
	\centering
	\includegraphics[width=\textwidth]{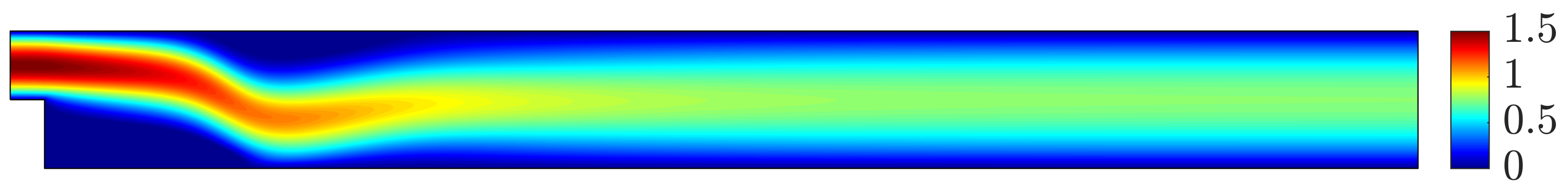}
	\includegraphics[width=\textwidth]{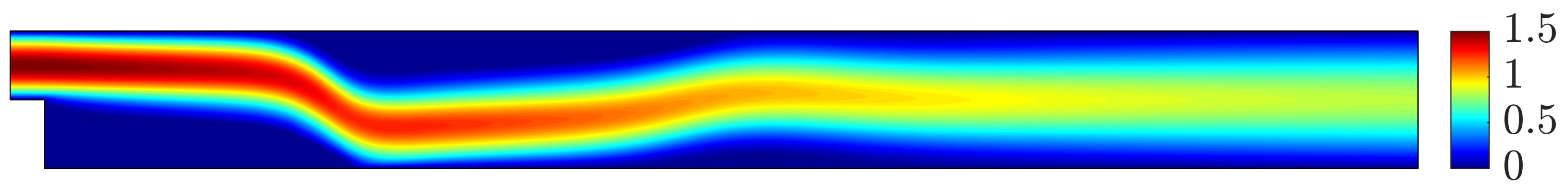}
	\includegraphics[width=\textwidth]{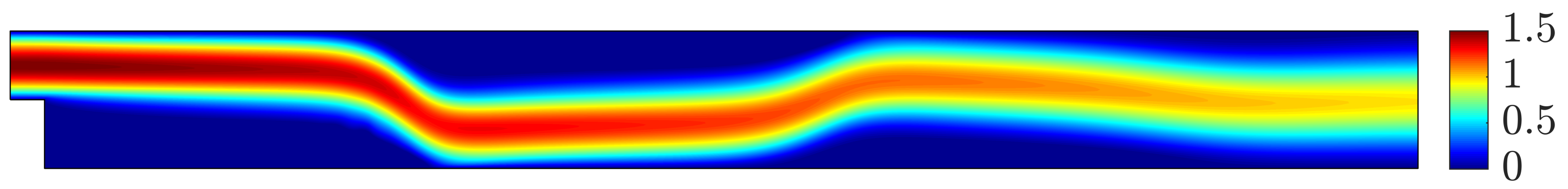}
	\caption{Navier-Stokes equations - Backward facing step: magnitude of the velocity field for $Re=800$ (top), $Re=1,900$ (middle) and $Re=3,000$ (bottom).}
	\label{fig:NS-BFS-Velo}
\end{figure}

Figure\ \ref{fig:NS-BFS-Velo} shows the magnitude of the velocity for three different Reynolds numbers, namely $Re=800$, $Re=1,900$ and $Re=3,000$. 
The results illustrate the higher complexity induced by an increase of the Reynolds number. 

\begin{figure}[!tb]
	\centering
	\includegraphics[width=\textwidth]{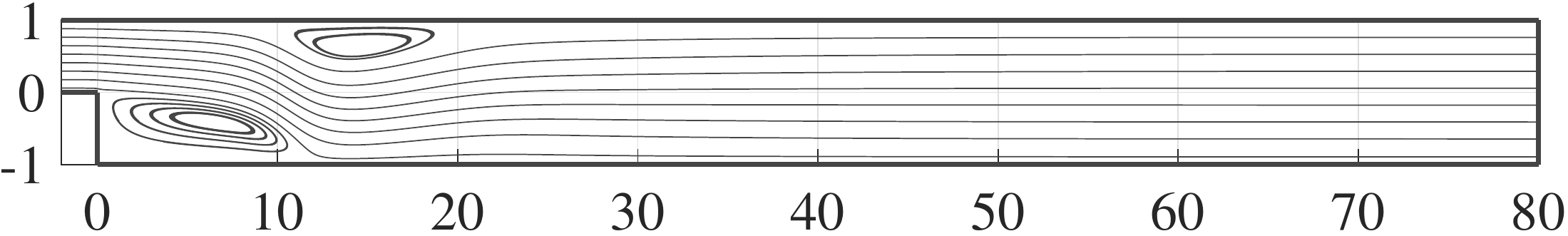}
	\includegraphics[width=\textwidth]{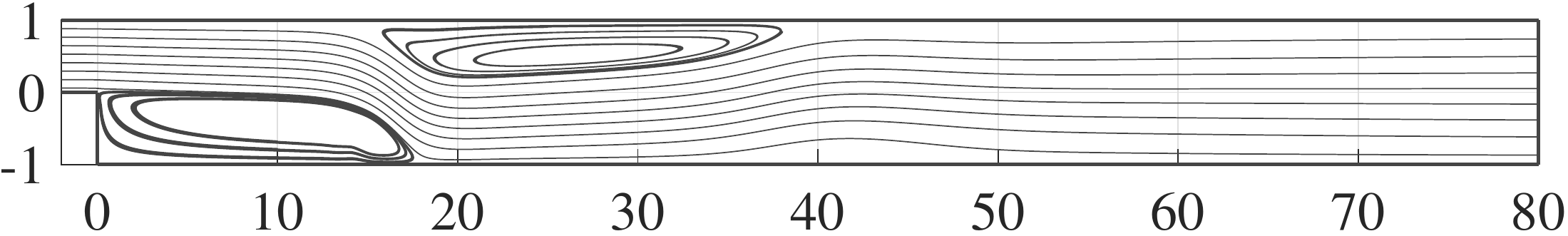}
	\includegraphics[width=\textwidth]{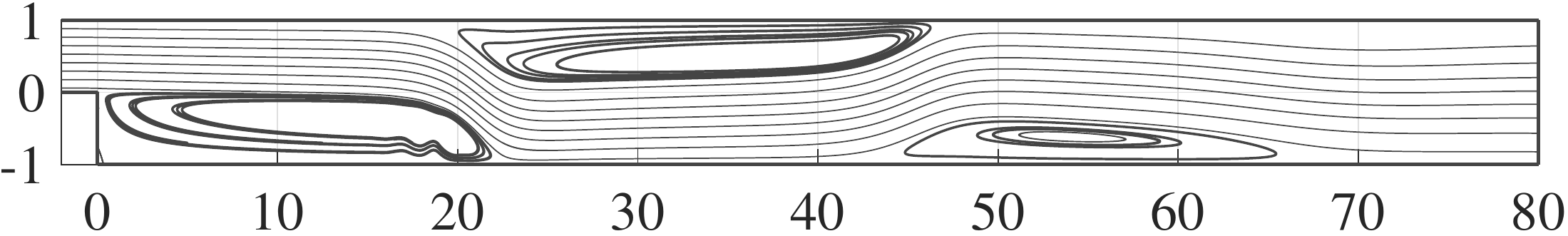}
	\caption{Navier-Stokes equations - Backward facing step: streamlines for $Re=800$ (top), $Re=1,900$ (middle) and $Re=3,000$ (bottom).}
	\label{fig:NS-BFS-Stream}
\end{figure}

To better observe the complexity of the flow and the different recirculation regions, Figure\ \ref{fig:NS-BFS-Stream} shows the streamlines of the velocity field for $Re=800$, $Re=1,900$ and $Re=3,000$. 
The results show the change in the number of recirculation regions as well as the change in the position of such regions as the Reynolds number is increased. 

\begin{figure}[!tb]
	\centering
	\includegraphics[width=\textwidth]{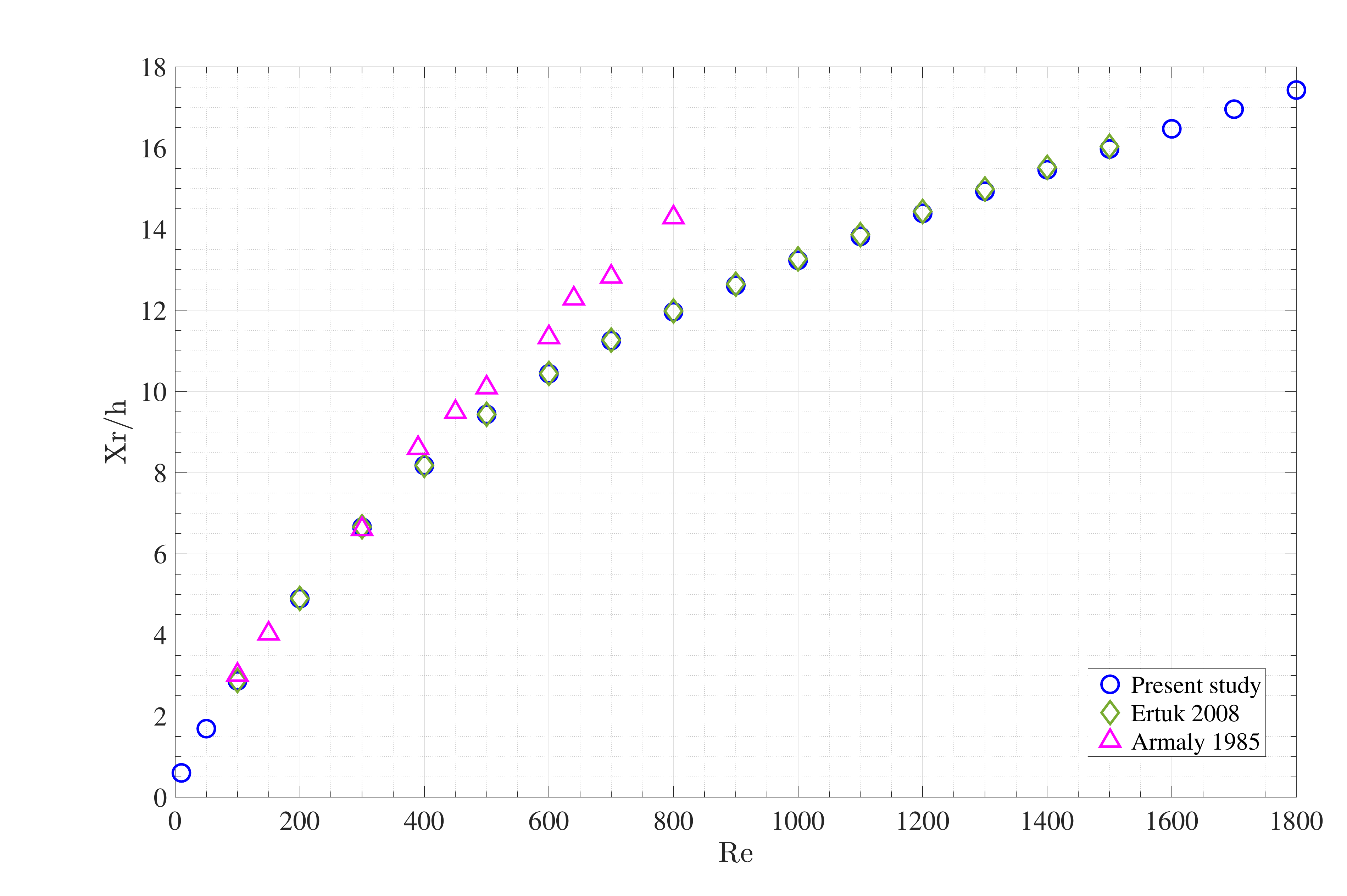}
	\caption{Navier-Stokes equations - Backward facing step: Position of the reattachment point as a function of the Reynolds number and comparison with published results.}
	\label{fig:NS-BFS_Reattachment}
\end{figure}

Finally, Figure\ \ref{fig:NS-BFS_Reattachment} shows a the position of the reattachment point as a function of the Reynolds number. The results obtained using the presented HDG-Voigt formulation are compared with the results in\ \cite{armaly1983experimental} and \cite{erturk2008numerical}, showing an excellent aggrement with the more recent results of\ \cite{erturk2008numerical}.

%________________________________________________________________________
\paragraph*{Acknowledgements}
This work is partially supported by the European Union's Horizon 2020 research and innovation programme under the Marie Sk\l odowska-Curie actions (Grant No. 675919 and 764636) and the Spanish Ministry of Economy and Competitiveness (Grant No. DPI2017-85139-C2-2-R). The first and third author also gratefully acknowledge the financial support provided by Generalitat de Catalunya (Grant No. 2017-SGR-1278).
%________________________________________________________________________

%________________________________________________________________________
\appendix
%________________________________________________________________________
\section{Appendix: Saddle-point structure of the global problem}\label{ap:Symmetry}

In this Appendix, the symmetry of the rectangular off-diagonal blocks in\ \eqref{eq:globalProblemSystemFinal} is demonstrated.
First, rewrite \eqref{eq:globalProblemSystemFinal} as
\begin{equation}\label{eq:globalNonSym}
  \begin{bmatrix}\widehat{\mat{K}} & \mat{H} \\
                   \mat{G}^T           & \mat{0}  \end{bmatrix}
  \begin{Bmatrix}\vect{\hu} \\ \bm{\rho} \end{Bmatrix} 
  =
  \begin{Bmatrix}\vect{\hat{f}}_{\hu} \\ \vect{\hat{f}}_{\rho} \end{Bmatrix},
\end{equation}
where the block $\mat{H}$ is obtained by the solution of the local problem in\ \eqref{eq:localProblemSystem} and has the following form
\begin{equation}\label{eq:eqH}
  \mat{H} = \Assem_{e=1}^{\numel}
     \left[\begin{array}{@{}c@{\,}c@{\,}c@{\,}c@{}} \mat{A}_{L \hu}^T & \mat{A}_{\hu u} & \mat{A}_{p\hu}^T & \mat{0}\end{array}\right]_{e}
     \left[\begin{array}{@{}c@{\,}c@{\,}c@{\,}c@{}}
       \mat{A}_{LL}    & \mat{A}_{Lu} & \mat{0}              & \mat{0}                   \\
       \mat{A}_{Lu}^T & \mat{A}_{uu} & \mat{A}_{pu}^T & \mat{0}                   \\
       \mat{0}             & \mat{A}_{pu} & \mat{0}              & \mat{a}_{\rho p}^T \\
       \mat{0}             & \mat{0}          & \mat{a}_{\rho p}& 0                   \\
     \end{array}\right]_{\! e}^{-1}
   \begin{Bmatrix}
     \vect{0} \\
     \vect{0} \\
     \vect{0} \\
     1
    \end{Bmatrix}_{\! e} .
\end{equation}

In order for the system in\ \eqref{eq:globalNonSym} to have a saddle-point structure, it needs to be proved that $\mat{H} = \mat{G}$.
For the sake of readability, rewrite the matrix of the local problem in\ \eqref{eq:localProblemSystem} using the block structure
\begin{equation}\label{eq:LocalBlockMatrix}
\Ke := 
      \begin{bmatrix}
        \Be    & \Ce \\
        \Ce^T & \De
      \end{bmatrix}
\end{equation}
where the blocks are defined as
\begin{equation*}%\label{eq:Blocks}
\Be :=
      \begin{bmatrix}
        \mat{A}_{LL}    & \mat{A}_{Lu} \\
        \mat{A}_{Lu}^T & \mat{A}_{uu}
      \end{bmatrix} 
      , \quad
\Ce :=      
      \begin{bmatrix}
       \mat{0}              & \mat{0} \\
       \mat{A}_{pu}^T & \mat{0}                   
      \end{bmatrix}
      , \quad
\De :=            
      \begin{bmatrix}
        \mat{0}              & \mat{a}_{\rho p}^T \\
        \mat{a}_{\rho p}& \mat{0}
      \end{bmatrix} .
\end{equation*}
\begin{proposition}
For each element $\Omega_e$, it holds
\begin{equation}\label{eq:locMatrix}
\inv{\Ke}
    \begin{Bmatrix}
      \vect{0} \\
      \vect{0} \\
      \vect{0} \\
      1
%      \abs{\partial\Omega_e}
    \end{Bmatrix}_{\! e}
= 
%\abs{\partial\Omega_e}
      \begin{bmatrix}
        - \begin{bmatrix} \inv{\Be} \end{bmatrix}_{12} \mat{A}_{pu}^T \begin{bmatrix} \inv{\Se} \end{bmatrix}_{12}\\[0.75ex]
        - \begin{bmatrix} \inv{\Be} \end{bmatrix}_{22} \mat{A}_{pu}^T \begin{bmatrix} \inv{\Se} \end{bmatrix}_{12}\\[0.75ex]
        \begin{bmatrix} \inv{\Se} \end{bmatrix}_{12}\\
       0
      \end{bmatrix}_{\! e} ,
\end{equation}
where 
\begin{equation*}%\label{eq:blockResult}
\begin{aligned}
\begin{bmatrix} \inv{\Be} \end{bmatrix}_{12} & = 
- \inv{\mat{A}_{LL}} \mat{A}_{Lu} \inv{(\mat{A}_{uu} - \mat{A}_{Lu}^T \inv{\mat{A}_{LL}} \mat{A}_{Lu})}  
, \\
\begin{bmatrix} \inv{\Be} \end{bmatrix}_{22} & = 
\inv{(\mat{A}_{uu} - \mat{A}_{Lu}^T \inv{\mat{A}_{LL}} \mat{A}_{Lu})} 
, \\
\begin{bmatrix} \inv{\Se} \end{bmatrix}_{12} & = 
 \pinv{\Bigl( \mat{a}_{\rho p} \bigl( \mat{I} - \pinv{( \mat{A}_{pu} \inv{(\mat{A}_{uu} - \mat{A}_{Lu}^T \inv{\mat{A}_{LL}} \mat{A}_{Lu})} \mat{A}_{pu}^T )} 
 \\[-2ex] &\hspace{53pt}
                                                      (\mat{A}_{pu} \inv{(\mat{A}_{uu} - \mat{A}_{Lu}^T \inv{\mat{A}_{LL}} \mat{A}_{Lu})} \mat{A}_{pu}^T ) \bigr) \Bigr)} 
.
\end{aligned}
\end{equation*}
\begin{proof}
The inverse of the block matrix in\ \eqref{eq:LocalBlockMatrix}, written using Schur-Banachiewicz form, see\ \cite[Section 2.17]{Bernstein-book}, is
\begin{equation*}%\label{eq:InverseBlockMatrix}
\inv{\Ke} {:=} 
      \left[\begin{array}{@{}c@{\;}c@{}}
        \inv{\Be} ( \mat{I} {+} \Ce \inv{(\De{-}\Ce^T\inv{\Be}\Ce)}\Ce^T\inv{\Be} )   & {-}\inv{\Be}\Ce \inv{(\De{-}\Ce^T\inv{\Be}\Ce)} \\
        {-}\inv{(\De{-}\Ce^T\inv{\Be}\Ce)}\Ce^T\inv{\Be} & \inv{(\De{-}\Ce^T\inv{\Be}\Ce)}
      \end{array}\right] ,
\end{equation*}
where the block $(2,2)$ is the inverse of the Schur complement
\begin{equation}\label{eq:SchurComplementGlobal}
\Se := \De-\Ce^T\inv{\Be}\Ce =
      \begin{bmatrix}
        - \mat{A}_{pu} \begin{bmatrix} \inv{\Be} \end{bmatrix}_{22} \mat{A}_{pu}^T  & \mat{a}_{\rho p}^T  \\
        \mat{a}_{\rho p} & 0
      \end{bmatrix} 
\end{equation}
of block $\Be$ of the matrix $\Ke$.
Moreover, the block $(1,2)$ of the inverse matrix $\inv{\Ke}$ has the form
\begin{equation}\label{eq:eqK12}
  \begin{aligned}
  \begin{bmatrix} \inv{\Ke} \end{bmatrix}_{12} &= - \inv{\Be} \Ce \inv{\Se} \\
%= - 
%      \begin{bmatrix}
%        \begin{bmatrix} \inv{\Be} \end{bmatrix}_{11} & \begin{bmatrix} \inv{\Be} \end{bmatrix}_{12} \\
%        \begin{bmatrix} \inv{\Be} \end{bmatrix}_{12}^T & \begin{bmatrix} \inv{\Be} \end{bmatrix}_{22}
%      \end{bmatrix}
%      \begin{bmatrix}
%        \mat{0} & \mat{0} \\
%        \mat{A}_{pu}^T & \mat{0}
%      \end{bmatrix}
%      \begin{bmatrix}
%        \begin{bmatrix} \inv{\Se} \end{bmatrix}_{11}    & \begin{bmatrix} \inv{\Se} \end{bmatrix}_{12} \\
%        \begin{bmatrix} \inv{\Se} \end{bmatrix}_{12}^T & \begin{bmatrix} \inv{\Se} \end{bmatrix}_{22}
%      \end{bmatrix}
      &= -\begin{bmatrix}
        \begin{bmatrix} \inv{\Be} \end{bmatrix}_{12} \mat{A}_{pu}^T \begin{bmatrix} \inv{\Se} \end{bmatrix}_{11}    & \begin{bmatrix} \inv{\Be} \end{bmatrix}_{12} \mat{A}_{pu}^T \begin{bmatrix} \inv{\Se} \end{bmatrix}_{12} \\[0.75ex]
        \begin{bmatrix} \inv{\Be} \end{bmatrix}_{22} \mat{A}_{pu}^T \begin{bmatrix} \inv{\Se} \end{bmatrix}_{11} & \begin{bmatrix} \inv{\Be} \end{bmatrix}_{22} \mat{A}_{pu}^T \begin{bmatrix} \inv{\Se} \end{bmatrix}_{12}
      \end{bmatrix} .
  \end{aligned}
\end{equation}

To fully determine the inverse matrix $\inv{\Ke}$, the blocks of $\inv{\Be}$ and $\inv{\Se}$ need to be computed. 
%
%Following the rationale utilized above, the Schur-Banachiewicz form of $\inv{\Be}$ is
%
%\begin{equation}\label{eq:InverseB}
%\inv{\Be} := 
%      \begin{bmatrix}
%        \inv{\mat{A}_{LL}} ( \mat{I} + \mat{A}_{Lu} \inv{(\mat{A}_{uu} - \mat{A}_{Lu}^T \inv{\mat{A}_{LL}} \mat{A}_{Lu})} \mat{A}_{Lu}^T \inv{\mat{A}_{LL}} )  & - \inv{\mat{A}_{LL}} \mat{A}_{Lu} \inv{(\mat{A}_{uu} - \mat{A}_{Lu}^T \inv{\mat{A}_{LL}} \mat{A}_{Lu})}  \\
%        - \inv{(\mat{A}_{uu} - \mat{A}_{Lu}^T \inv{\mat{A}_{LL}} \mat{A}_{Lu})} \mat{A}_{Lu}^T \inv{\mat{A}_{LL}}  & \inv{(\mat{A}_{uu} - \mat{A}_{Lu}^T \inv{\mat{A}_{LL}} \mat{A}_{Lu})} 
%      \end{bmatrix} .
%\end{equation}
%
Following the Schur-Banachiewicz rationale utilized above, the blocks of $\inv{\Be}$ have the form
\begin{equation}\label{eq:InverseB}
\begin{aligned}
\begin{bmatrix} \inv{\Be} \end{bmatrix}_{11} & := 
\inv{\mat{A}_{LL}} ( \mat{I} + \mat{A}_{Lu} \inv{(\mat{A}_{uu} - \mat{A}_{Lu}^T \inv{\mat{A}_{LL}} \mat{A}_{Lu})} \mat{A}_{Lu}^T \inv{\mat{A}_{LL}} ) 
, \\
\begin{bmatrix} \inv{\Be} \end{bmatrix}_{12} & := 
- \inv{\mat{A}_{LL}} \mat{A}_{Lu} \inv{(\mat{A}_{uu} - \mat{A}_{Lu}^T \inv{\mat{A}_{LL}} \mat{A}_{Lu})}  
, \\
\begin{bmatrix} \inv{\Be} \end{bmatrix}_{22} & := 
\inv{(\mat{A}_{uu} - \mat{A}_{Lu}^T \inv{\mat{A}_{LL}} \mat{A}_{Lu})} 
,
\end{aligned}
\end{equation}
and, from the symmetry of $\Be$, it follows that $\begin{bmatrix} \inv{\Be} \end{bmatrix}_{21} = \begin{bmatrix} \inv{\Be} \end{bmatrix}_{12}^T$.

Plugging the expression of $\begin{bmatrix} \inv{\Be} \end{bmatrix}_{22} $, see\ \eqref{eq:InverseB}, into the definition of $\Se$ in\ \eqref{eq:SchurComplementGlobal}, it follows that the block $(1,1)$ of such matrix is
\begin{equation}\label{eq:eqS11}
\begin{bmatrix} \Se \end{bmatrix}_{11} := - \mat{A}_{pu} \inv{(\mat{A}_{uu} - \mat{A}_{Lu}^T \inv{\mat{A}_{LL}} \mat{A}_{Lu})} \mat{A}_{pu}^T .
\end{equation}
It is straightforward to observe that this matrix is the Schur complement of block
$$
      \begin{bmatrix}
        \mat{A}_{LL}    & \mat{A}_{Lu} \\
        \mat{A}_{Lu}^T & \mat{A}_{uu}
      \end{bmatrix} 
$$
of the matrix
$$
      \begin{bmatrix}
        \mat{A}_{LL}    & \mat{A}_{Lu} & \mat{0} \\
        \mat{A}_{Lu}^T & \mat{A}_{uu} & \mat{A}_{pu}^T \\
        \mat{0} & \mat{A}_{pu} & \mat{0}
      \end{bmatrix} ,
$$
which is singular, since it is obtained from the discretization of an incompressible flow problem with purely Dirichlet boundary conditions. Hence, to compute the blocks of $\inv{\Se}$, the framework of the generalized inverse of a partitioned matrix is exploited, see\ \cite{Miao-91}, leading to
\begin{equation}\label{eq:InverseS}
\begin{aligned}
  \begin{bmatrix} \inv{\Se} \end{bmatrix}_{11} & :=
      \Bigl( \mat{I} -  \pinv{ \bigl( \mat{a}_{\rho p} ( \mat{I} - \pinv{\begin{bmatrix} \Se \end{bmatrix}_{11}} \begin{bmatrix} \Se \end{bmatrix}_{11} ) \bigr) } 
                              \mat{a}_{\rho p} \Bigr) 
      \\[-1ex] &\hspace{80pt}
      \pinv{\begin{bmatrix} \Se \end{bmatrix}_{11}} 
      \Bigl( \mat{I} -  \mat{a}_{\rho p}^T \pinv{ \bigl( ( \mat{I} - \begin{bmatrix} \Se \end{bmatrix}_{11} \pinv{\begin{bmatrix} \Se \end{bmatrix}_{11}} ) 
                              \mat{a}_{\rho p}^T \bigr) } \Bigr) 
, \\
  \begin{bmatrix} \inv{\Se} \end{bmatrix}_{12} & :=
  \pinv{ \bigl( \mat{a}_{\rho p} ( \mat{I} - \pinv{\begin{bmatrix} \Se \end{bmatrix}_{11}} \begin{bmatrix} \Se \end{bmatrix}_{11} ) \bigr) }
, \\
  \begin{bmatrix} \inv{\Se} \end{bmatrix}_{21} & :=
  \pinv{ \bigl( ( \mat{I} - \begin{bmatrix} \Se \end{bmatrix}_{11} \pinv{\begin{bmatrix} \Se \end{bmatrix}_{11}} ) \mat{a}_{\rho p}^T \bigr) }
, \\
\begin{bmatrix} \inv{\Se} \end{bmatrix}_{22} & :=
0
\end{aligned}
\end{equation}
where the Moore-Penrose pseudoinverse $\pinv{\begin{bmatrix} \Se \end{bmatrix}_{11}}$ of the singular matrix $\begin{bmatrix} \Se \end{bmatrix}_{11}$ has the form 
\begin{equation}\label{eq:InverseS11}
  \pinv{ \begin{bmatrix} \Se \end{bmatrix}_{11} } 
    := - \pinv{ \bigl( \mat{A}_{pu} \inv{(\mat{A}_{uu} - \mat{A}_{Lu}^T \inv{\mat{A}_{LL}} \mat{A}_{Lu})} \mat{A}_{pu}^T \bigr) }.
\end{equation}

From\ \eqref{eq:eqH}, it is straightforward to observe that only the blocks in the last column of the inverse matrix are involved in the definition of the product in\ \eqref{eq:locMatrix}. The result\ \eqref{eq:locMatrix} follows directly from\ \eqref{eq:eqK12},\ \eqref{eq:InverseS} and\ \eqref{eq:InverseS11}.
\end{proof}
\end{proposition}

\begin{proposition}
Given $\mat{H}$ and $\mat{G}^T$ from\ \eqref{eq:eqH} and\ \eqref{eq:eqGT} respectively, it holds that $\mat{H} = \mat{G}$.
\begin{proof}
From\ \eqref{eq:eqH} and\ \eqref{eq:locMatrix}, it follows that 
\begin{equation}\label{eq:eqGe}
  \mat{H}_e 
%= \begin{bmatrix} \mat{A}_{L \hu}^T & \mat{A}_{\hu u} & \mat{A}_{p\hu}^T & \mat{0}\end{bmatrix}_{e} 
%\inv{\Ke}
%    \begin{Bmatrix}
%      \vect{0} \\
%      \vect{0} \\
%      \vect{0} \\
%      1
%%      \abs{\partial\Omega_e}
%    \end{Bmatrix}_{\! e}
  = - \bigl(\mat{A}_{L \hu}^T \begin{bmatrix} \inv{\Be} \end{bmatrix}_{12} + \mat{A}_{\hu u} \begin{bmatrix} \inv{\Be} \end{bmatrix}_{22} \bigr) 
            \mat{A}_{pu}^T \begin{bmatrix} \inv{\Se} \end{bmatrix}_{12} + \mat{A}_{p\hu}^T \begin{bmatrix} \inv{\Se} \end{bmatrix}_{12} ,
\end{equation}
for each element $\Omega_e$.

First, recall that the matrix $\mat{I} - \pinv{\begin{bmatrix} \Se \end{bmatrix}_{11}} \begin{bmatrix} \Se \end{bmatrix}_{11}$ defines an orthogonal projector onto the kernel of $\begin{bmatrix} \Se \end{bmatrix}_{11}$, see\ \cite[Section 6.1]{Bernstein-book}.
As observed in the previous proposition, see\ \eqref{eq:eqS11}, $\begin{bmatrix} \Se \end{bmatrix}_{11}$ is the Schur complement of the velocity block of the matrix obtained from the discretization of an incompressible flow problem with purely Dirichlet boundary conditions. Thus, the kernel of $\begin{bmatrix} \Se \end{bmatrix}_{11}$ contains all constant vectors representing the mean value of pressure. It follows that
\begin{equation*}%\label{eq:KernelElement}
  \mat{a}_{\rho p} \bigl( \mat{I} - \pinv{\begin{bmatrix} \Se \end{bmatrix}_{11}} \begin{bmatrix} \Se \end{bmatrix}_{11} \bigr) 
  = \frac{1}{\nen} \vect{1}^T 
\end{equation*}
is the constant vector obtained as the average of $1$ over the $\nen$ element nodes of $\Omega_e$ and, consequently, $\begin{bmatrix} \inv{\Se} \end{bmatrix}_{12} = \vect{1}$. Moreover, since the kernel of the matrix $\mat{A}_{pu}^T$ also includes all constant vectors, $\mat{A}_{pu}^T \begin{bmatrix} \inv{\Se} \end{bmatrix}_{12} = \vect{0}$.
Hence, from\ \eqref{eq:eqGe}, it follows that
\begin{equation*}%\label{eq:eqGnew}
  \mat{H} 
 = \begin{bmatrix} [\mat{H}_1] & [\mat{H}_2] & \cdots & [\mat{H}_{\numel}] \end{bmatrix}  
 = \begin{bmatrix} [\mat{A}_{p\hu}^T]_1 \vect{1} & [\mat{A}_{p\hu}^T]_2 \vect{1} & \cdots & [\mat{A}_{p\hu}^T]_{\numel} \vect{1} \end{bmatrix}
\end{equation*}
which proves the statement.
\end{proof}
\end{proposition}

When convection phenomena are neglected (Stokes flow), $\mat{A}_{\hu u} = \mat{A}_{u \hu}^T$ and the symmetry of $\widehat{\mat{K}}$ and the global matrix in\ \eqref{eq:globalProblemSystemFinal} follows straightforwardly.
For general incompressible flow problems, the matrix $\widehat{\mat{K}}$ is not symmetric but the off-diagonal blocks $\mat{G}$ and $\mat{G}^T$ are one the transpose of the other and the resulting global matrix maintains the above displayed saddle-point structure\ \cite{BenziGolub05}.

%________________________________________________________________________
\section{Appendix: Implementation details}\label{ap:Implementation}

In this Appendix, the matrices and vectors appearing in the discrete form of the HDG-Voigt approximation of the Oseen equations are detailed.
The elemental variables $\bu$, $p$ and $\bL$ are defined in a reference element $\widetilde{\Omega}(\bxi), \ \bxi = (\xi_1, \ldots, \xi_{\nsd})$ whereas the face variable $\bhu$, is defined on a reference face $\widetilde{\Gamma}(\bet), \ \bet = (\eta_1, \ldots, \eta_{\nsd-1})$ as
\begin{align*} 
  \bu(\bxi) &\simeq \sum_{j=1}^{\nen} \nodaluV_j N_j(\bxi), &
  p(\bxi) &\simeq \sum_{j=1}^{\nen} \nodalpV_j N_j(\bxi), \\
  \bL(\bxi) &\simeq \sum_{j=1}^{\nen} \nodalLV_j N_j(\bxi), &
  \bhu(\bet) &\simeq \sum_{j=1}^{\nfn} \nodaluhV_j \hat{N}_j(\bet),
\end{align*} 
where $\nodaluV_j, \nodalpV_j, \nodalLV_j$ and $\nodaluhV_j $ are the nodal values of the approximation, $\nen$ and $\nfn$ the number of nodes in the element and face, respectively and $N_j$ and $\hat{N}_j$ the polynomial shape functions in the reference element and face, respectively.

An isoparametric formulation is considered and the following transformation is used to map reference and local coordinates
\begin{equation*}
\bx(\bxi) = \sum_{k=1}^{\nen} \bx_k N_k(\bxi),
\end{equation*}
where the vector $\{\bx_k\}_{k=1,\ldots,\nen}$ denotes the elemental nodal coordinates. 

Following\ \cite{RS-SGKH:18}, the matrices $\gradS$ and  $\bN$ in\ \eqref{eq:symmGrad} and\ \eqref{eq:normalVoigt}, respectively, are expressed in compact form as
\begin{equation*}
\gradS = \sum_{k=1}^{\nsd} \mat{F}_k \frac{\partial }{\partial x_k}, \qquad  \bN = \sum_{k=1}^{\nsd} \mat{F}_k n_k ,
\end{equation*}
where the matrices $\mat{F}_k$ are defined as
\begin{equation*}
\mat{F}_1 = 
\begin{bmatrix}
1 & 0 & 0 \\
0 & 0 & 1
\end{bmatrix}^T
, \quad
\mat{F}_2 = 
\begin{bmatrix}
0 & 0 & 1 \\
0 & 1 & 0
\end{bmatrix}^T
\end{equation*}
in two dimensions and 
\begin{equation*}
  \mat{F}_1 = \left[\begin{array}{@{}c@{\;}c@{\;}c@{\;}c@{\;}c@{\;}c@{}}
                      1 & 0 & 0 & 0 & 0 & 0 \\
                      0 & 0 & 0 & 1 & 0 & 0 \\
                      0 & 0 & 0 & 0 & 1 & 0 \end{array}\right]^T
  , \;
  \mat{F}_2 = \left[\begin{array}{@{}c@{\;}c@{\;}c@{\;}c@{\;}c@{\;}c@{}}
                      0 & 0 & 0 & 1 & 0 & 0 \\
                      0 & 1 & 0 & 0 & 0 & 0 \\
                      0 & 0 & 0 & 0 & 0 & 1 \end{array}\right]^T
  , \;
  \mat{F}_3 = \left[\begin{array}{@{}c@{\;}c@{\;}c@{\;}c@{\;}c@{\;}c@{}}
                      0 & 0 & 0 & 0 & 1 & 0 \\
                      0 & 0 & 0 & 0 & 0 & 1 \\
                      0 & 0 & 1 & 0 & 0 & 0 \end{array}\right]^T
\end{equation*}
in three dimensions.
Moreover, from the definition of $\bE$ in\ \eqref{eq:EDVoigt}, it holds $\bN^T \bE = \bn$ and $\gradS^T\bE = \grad$ for the gradient operator applied to a scalar function.
The following compact forms of the shape functions and their derivatives are introduced
\begin{align*}
  \Nmat & := \begin{bmatrix} N_1\Insd & N_2\Insd & \dots & N_{\nen}\Insd \end{bmatrix}^T , \\
  \NmatHat &:= \begin{bmatrix} \hat{N}_1\Insd & \hat{N}_2\Insd & \dots & \hat{N}_{\nfn}\Insd \end{bmatrix}^T , \\
  \Nmat^{\tau} & := \begin{bmatrix} N_1\tau \Insd & N_2\tau \Insd & \dots & N_{\nen}\tau \Insd \end{bmatrix}^T , \\
  \NmatHat^{\tau} &:= \begin{bmatrix} \hat{N}_1\tau \Insd & \hat{N}_2 \tau \Insd & \dots & \hat{N}_{\nfn} \tau \Insd \end{bmatrix}^T , \\
  \Nmat^n & := \begin{bmatrix} N_1\bn & N_2\bn & \dots & N_{\nen}\bn \end{bmatrix}^T , \\
  \Mmat &:= \begin{bmatrix} N_1\Imsd & N_2\Imsd & \dots & N_{\nen}\Imsd \end{bmatrix}^T , \\
  \Nmat^a & := \begin{bmatrix} N_1(\tau - \bha {\cdot} \bn) \Insd & N_2 (\tau - \bha {\cdot} \bn) \Insd & \dots & N_{\nen} (\tau - \bha {\cdot} \bn) \Insd \end{bmatrix}^T , \\
  \NmatHat^a & := 
  \begin{bmatrix} \hat{N}_1 (\tau - \bha {\cdot} \bn) \Insd & \hat{N}_2 (\tau - \bha {\cdot} \bn) \Insd & \dots & \hat{N}_{\nfn} (\tau - \bha {\cdot} \bn) \Insd \end{bmatrix}^T , \\
  \Qmat & := \begin{bmatrix} (\bJ^{-1} \grad \! N_1)^T & (\bJ^{-1} \grad \! N_2)^T& \dots &  (\bJ^{-1} \grad \! N_{\nen})^T\end{bmatrix}^T , \\
\Qmat^a & :=
\begin{bmatrix} 
\bm{a} \cdot (\bJ^{-1} \grad \! N_1) & \bm{a} \cdot (\bJ^{-1} \grad \! N_2)  & \dots & \bm{a} \cdot (\bJ^{-1} \grad \! N_{\nen})
\end{bmatrix}^T ,
\end{align*}
where $\bn$ is the outward unit normal vector to a face, $\ba$ is the convection field evaluated in the reference element, and $\bha$ is the convection field evaluated on the reference face. Moreover, for each spatial dimension, that is, for $k=1,\ldots,\nsd$, define
\begin{align*}
  \Nmat^{\mathrm{D}}_k & {:=}\!
  \begin{bmatrix} N_1 n_k \mat{F}^T_k \bDHalf & N_2 n_k \mat{F}^T_k \bDHalf & \dots &  N_{\nfn} n_k \mat{F}^T_k \bDHalf \end{bmatrix}^T, \\
  \NmatHat^{\mathrm{D}}_k & {:=}\! 
  \begin{bmatrix} \hat{N}_1 n_k \mat{F}^T_k \bDHalf & \hat{N}_2 n_k \mat{F}^T_k \bDHalf & \dots &  \hat{N}_{\nfn} n_k \mat{F}^T_k \bDHalf \end{bmatrix}^T, \\
  \Qmat^{\mathrm{D}}_k & {:=}\! \left[\begin{array}{@{}c@{\;}c@{\;}c@{\;}c@{}}
  [ \bJ^{-1} \grad \! N_1 ]_{k}  \mat{F}^T_k \bDHalf & [ \bJ^{-1} \grad \! N_2 ]_{k}  \mat{F}^T_k \bDHalf & \dots & [ \bJ^{-1} \grad \! N_{\nen} ]_{k}  \mat{F}^T_k \bDHalf \end{array}\right]^T , 
\end{align*}
where $n_k$ is the $k$-th components of the outward unit normal vector $\bn$ to a face  and $\bJ$ is the Jacobian of the isoparametric transformation.

The discretization of\ \eqref{eq:O-Loc-L} leads to the following matrices and vector
\begin{align*}
  [\mat{A}_{LL} ]_e & = -\sumge \Mmat(\bxige) \Mmat^T\!(\bxige) \abs{\bJ(\bxige)} \wge , \\[-1ex]
  [\mat{A}_{Lu}]_e & = \sum_{k=1}^{\nsd} \sumge \Qmat^{\mathrm{D}}_k(\bxige) \Nmat^T\!(\bxige) \abs{\bJ(\bxige)} \wge , \\[-1ex]
  [\mat{A}_{L \hat{u}}]_e & = \sum_{f=1}^{\numfa^e} 
                                              \biggl( \sum_{k=1}^{\nsd} \sumgf \Nmat^{\mathrm{D}}_k (\bxigf) \NmatHat^T\!(\bxigf) \abs{\bJ(\bxigf)} \wgf \biggr) 
                                              \bigl(1 - \chi_{\Gamma_D}(f) \bigr), \\[-1ex]
  [\vect{f}_{L}]_e & = \sum_{f=1}^{\numfa^e} 
                                   \biggl( \sum_{k=1}^{\nsd} \sumgf \Nmat^{\mathrm{D}}_k(\bxigf) \bu_D \!\bigl(\bx(\bxigf) \bigr) \abs{\bJ(\bxigf)} \wgf \biggr) 
                                   \chi_{\Gamma_D}(f) ,
\end{align*}
where $\numfa^e$ is the number of faces $\Gamma_{e,f}, \ f = 1,\dots,\numfa^e$, of the element $\Omega_e$, $\bxige$ and $\wge$ (resp., $\bxigf$ and $\wgf$ ) are the $\nipe$ (resp., $\nipf$) integration points and weights defined on the reference element (resp., face) and $\chi_{\Gamma_D}$ is the indicator function of the boundary $\Gamma_D$, namely
\begin{equation*}
\chi_{\Gamma_D}(f) = \begin{cases}%\Bigg\{\begin{array}{ll}
  1 & \text{ if } \Gamma_{e,f} \cap \Gamma_D \neq \emptyset  \\
  0 & \text{ otherwise}
  \end{cases}%\end{array}.
\end{equation*}
Similarly, from the discretization of\ \eqref{eq:O-Loc-u} the following matrices and vectors are obtained
\begin{align*}
   [\mat{A}_{uu}]_e & = - \sumge \Qmat^a(\bxige) \Nmat^T\!(\bxige) \abs{\bJ(\bxige)} \wge 
   \\[-2.5ex]&\hspace{90pt}
                                    + \sum_{f=1}^{\numfa^e}  \sumgf \Nmat^{\tau }\!(\bxigf) \Nmat^T\!(\bxigf) \abs{\bJ(\bxigf)} \wgf, \\[-1ex]
  [\mat{A}_{u \hat{u}}]_e & = \sum_{f=1}^{\numfa^e} 
                                             \biggl(  \sumgf \Nmat^a\!(\bxigf) \NmatHat^T\!(\bxigf) \abs{\bJ(\bxigf)} \wgf \biggr)  
                                             \bigl(1 - \chi_{\Gamma_D}(f) \bigr), \\[-1ex]
  [\vect{f}_{u}]_e & = \sumge \Nmat(\bxige)\bm{s}\!\bigl(\bx(\bxige) \bigr) \abs{\bJ(\bxige)} \wge 
   \\[-2.5ex]&\hspace{90pt}
                              + \sum_{f=1}^{\numfa^e} \biggl(  \sumgf \Nmat^a\!(\bxigf) \bu_D\!\bigl(\bx(\bxigf) \bigr) \abs{\bJ(\bxigf)} \wgf \biggr) \chi_{\Gamma_D}(f) .
\end{align*}
%
%where $\bn_{e,f}$ is the unit outward normal to the face $\Gamma_{e,f}$ of element $\Omega_e$ and $\bC$ is a diagonal matrix of dimension $\nen \nsd \times \nen \nsd$ accounting for the information of the convection field $\ba$. More precisely, $\bC_{ii}\!\left(\bx(\bxige) \right) := a_{\mathrm{p}(g(i))}\!\left(\bx(\bxige) \right)$, where $g(i) := \nsd - i \mathbin{\%} \nsd$, $\mathbin{\%}$ being the remainder of the division $i/\nsd$, and $\mathrm{p}$ denotes the permutation index
%%
%\begin{equation*}
%\mathrm{p} := 
%\begin{cases}
%[1 \;\, 2]^T , & \nsd =2 \\
%[2 \;\, 1 \;\, 3]^T , & \nsd =3
%\end{cases} .
%\end{equation*}
%
The discrete forms of the incompressibility constraint in\ \eqref{eq:O-Loc-p} and the restriction in\ \eqref{eq:O-Loc-rho} feature the following matrices and vector
\begin{align*}
  [\mat{A}_{pu}]_e & = \sumge \Nmat(\bxige) \Qmat^T\!(\bxige)  \abs{\bJ(\bxige)} \wge, \\[-1ex]
  [\mat{A}_{p \hat{u}}]_e & = \sum_{f=1}^{\numfa^e} 
                                             \biggl(  \sumgf \Nmat^n(\bxigf) \NmatHat^T\!(\bxigf) \abs{\bJ(\bxigf)} \wgf \biggr)
                                             \bigl(1 - \chi_{\Gamma_D}(f) \bigr) , \\[-1ex]
  [\vect{f}_{p}]_e & = \sum_{f=1}^{\numfa^e} 
                                 \biggl(  \sumgf \Nmat^n(\bxigf) \bu_D\!\bigl(\bx(\bxigf) \bigr) \abs{\bJ(\bxigf)} \wgf \biggr) \chi_{\Gamma_D}(f) , \\[-1ex]
  [\mat{a}_{\rho p}]_e & = \sum_{f=1}^{\numfa^e} \sumgf \Nmat(\bxigf) \vect{1} \abs{\bJ(\bxigf)} \wgf  ,
\end{align*}
Finally, the matrices and vectors resulting from the discretization of the global problem in\ \eqref{eq:O-Glob-Transmission} are
\begin{align*}
  [\mat{A}_{\hat{u} \hat{u}}]_e & =  
	- \sum_{f=1}^{\numfa^e} \biggl( \sumgf \NmatHat^{\tau}\!(\bxigf) \NmatHat^T(\bxigf) \abs{\bJ(\bxigf)} \wgf \biggr) \chi_{\Gamma}(f) \\
   \\[-3em]&\hspace{90pt}
	- \sum_{f=1}^{\numfa^e} \biggl( \sumgf \NmatHat^a\!(\bxigf) \NmatHat^T(\bxigf) \abs{\bJ(\bxigf)} \wgf \biggr) \chi_{\Gamma_N}(f) , \\[-1ex]
  [\mat{A}_{\hat{u} u}]_e & = \sum_{f=1}^{\numfa^e} 
                                              \biggl( \sumgf \NmatHat^{\tau}\!(\bxigf) \Nmat^T\!(\bxigf) \abs{\bJ(\bxigf)} \wgf \biggr)  
                                              \bigl(1 - \chi_{\Gamma_D}(f) \bigr) , \\[-1ex]
  [\vect{f}_{\hat{u}}]_e & = - \sum_{f=1}^{\numfa^e} \biggl( \sumgf \NmatHat(\bxigf) \bt\!\bigl(\bx(\bxigf) \bigr) \abs{\bJ(\bxigf)} \wgf \biggr) \chi_{\Gamma_N}(f) , 
\end{align*}
where $\chi_{\Gamma}$ and $\chi_{\Gamma_N}$ are the indicator functions of the internal skeleton $\Gamma$ and the Neumann boundary $\Gamma_N$, respectively.

%______________________________________________________________________________
%\bibliographystyle{amsplain}
%\bibliographystyle{abbrv}
\bibliographystyle{unsrt}
\bibliography{references}

\end{document}